\input amstex
\documentstyle{amsppt}
\loadmsbm

\nologo

\TagsOnRight

\NoBlackBoxes

\define\acc{\operatorname{acc}}

\define\supp{\operatorname{supp}}
\define\Lip{\operatorname{Lip}}
\define\diam{\operatorname{diam}}

\def\floor{\mathbin{\hbox{\vrule height1.2ex width0.8pt depth0pt
        \kern-0.8pt \vrule height0.8pt width1.2ex depth0pt}}}

\font\letter=cmss10 

\font\normalsmall=cmss10 scaled 500

\font\normal=cmss10 scaled 700

\font\normalbig=cmss10

\define\ini{\operatorname{i}}

\define\termi{\operatorname{t}}

\define\smallsmallG{\text{\normalsmall G}}
\define\smallsmallvertexi{\text{\normalsmall i}}
\define\smallsmallvertexj{\text{\normalsmall j}}

\define\smalledge{\text{\normal e}}
\define\smallE{\text{\normal E}}
\define\smallV{\text{\normal V}}
\define\smallG{\text{\normal G}}
\define\smallvertexi{\text{\normal i}}
\define\smallvertexj{\text{\normal j}}

\define\edge{\text{\normalbig e}}

\define\E{\text{\normalbig E}}
\define\V{\text{\normalbig V}}
\define\G{\text{\normalbig G}}
\define\vertexi{\text{\normalbig i}}
\define\vertexj{\text{\normalbig j}}

\define\Normal{\text{\normal N}}

\define\smallNormal{\text{\normalsmall N}}

\define\Haus{\text{\normal H}}

\define\scon{\text{\normal con}}

\define\sdyncon{\text{\normal dyn-con}}

\define\dyn{\text{\normal dyn}}

\define\vector{\text{\normal vec}}

\define\Hol{\text{\normal H\"ol}}

\define\co{\text{\normal co}}

\define\erg{\text{\normal erg}}

\define\specradsmall{\text{\normal spec-rad}}

\define\subdif{\text{\normal sub}}

\define\abs{\text{\normal ab}}

\define\radius{\text{\normal rad}}

\define\distance{\text{\letter d}}

\define\LDistance{\text{\letter L}}

\font\tenscr=callig15 scaled 800
\font\sevenscr=callig15
\font\fivescr=callig15
\skewchar\tenscr='177 \skewchar\sevenscr='177 \skewchar\fivescr='177
\newfam\scrfam \textfont\scrfam=\tenscr \scriptfont\scrfam=\sevenscr
\scriptscriptfont\scrfam=\fivescr
\def\scr#1{{\fam\scrfam#1}}

\font\tenscri=suet14
\font\sevenscri=suet14
\font\fivescri=suet14
\skewchar\tenscri='177 \skewchar\sevenscri='177 \skewchar\fivescri='177
\newfam\scrifam \textfont\scrifam=\tenscri \scriptfont\scrifam=\sevenscri
\scriptscriptfont\scrifam=\fivescri
\def\scri#1{{\fam\scrifam#1}}

\hsize = 6.1 true in
\vsize = 9.33 true in

%
%

\topmatter
\title
\centerline{Fine and coarse 
multifractal zeta-functions:}
On the multifractal
formalism for
multifractal 
zeta-functions
\endtitle
\endtopmatter

\centerline{\smc L\. Olsen}
\centerline{Department of Mathematics}
\centerline{University of St\. Andrews}
\centerline{St\. Andrews, Fife KY16 9SS, Scotland}
\centerline{e-mail: {\tt lo\@st-and.ac.uk}}

\topmatter
\abstract{Multifractal analysis
refers to the study of 
the local 
properties of measures 
and functions,
and consists of two parts:
the fine multifractal theory 
and the coarse multifractal theory.
The fine and the coarse theory are linked by a 
web of 
conjectures known collectively as the Multifractal Formalism.
Very roughly speaking the
Multifractal Formalism
says that the multifractal spectrum  
from
fine theory 
equals 
the Legendre transform of the
Renyi dimensions 
from
 the coarse theory.

Recently 
{\it fine}
multifractal zeta-functions,
i\.e\.
multifractal zeta-functions
designed to 
produce detailed information
about the fine multifractal theory,
have been introduced and investigated.
The purpose of this work is to
complement and expand this study by
introducing and investigating
{\it coarse}
multifractal zeta-functions,
i\.e\.
multifractal zeta-functions
designed to 
produce
information 
about the coarse multifractal theory,
and, in particular,
to
establish a
{\it Multifractal Fortmalism for Zeta-Functions}
 linking
fine multifractal zeta-functions and coarse multifractal zeta-functions
via the Legendre transform.

Several applications are given, including applications to
multifractal analysis of 
graph-directed self-conformal 
 measures 
 and multifractal analysis of ergodic Birkhoff averages of continuous functions
 on graph-directed self-conformal sets.

 In judiciously chosen examples 
the {\it Multifractal Fortmalism for Zeta-Functions}
reduces to
known 
results linking 
different types of classical
multifractal spectra
to 
the Legendre transform of
certain
Renyi dimensions,
and in other examples
the {\it Multifractal Fortmalism for Zeta-Functions}
provides 
new 
results
linking
multifractal spectra
to 
the Legendre transform of
various
Renyi dimensions.
 In particular,
 our results leads to new 
 representations for 
 many
 different types of  multifractal spectra of
 ergodic Birkhoff averages of continuous functions.
}
\endabstract
\endtopmatter

%
%
%
%
%
%
%

\footnote""
{
\!\!\!\!\!\!\!\!
}

\leftheadtext{L\. Olsen}

\rightheadtext{On the multifractal
formalism for
zeta-functions}


\heading{1. Introduction.}\endheading

Mutifractal analysis
emerged in the 1990's as a powerful tool for analysing the local behaviour of 
measures with widely varying intensity.
Multifractal analysis consists of two distinct  parts:
the 
{\it fine} part and the 
{\it coarse} part.
the fine multifractal theory 
and the coarse multifractal theory.
The fine and the coarse theory are linked, via the Legendre transform,
by a 
web of 
conjectures known collectively as the Multifractal Formalism.

Recently 
{\it fine}
multifractal zeta-functions,
i\.e\.
multifractal zeta-functions
designed to 
produce detailed information
about the fine multifractal theory,
have been introduced and investigated.
The purpose of this work is to
complement and expand this study by
introducing and investigating
{\it coarse}
multifractal zeta-functions,
i\.e\.
multifractal zeta-functions
designed to 
produce
information 
about the coarse multifractal theory,
and, in particular,
to
establish a
{\it Multifractal Fortmalism for Zeta-Functions}
 linking
fine multifractal zeta-functions and coarse multifractal zeta-functions
via the Legendre transform.

\medskip

{\bf Multifractal analysis.}
Multifractal analysis consists of two distinct  parts:
the 
{\it fine} part and the 
{\it coarse} part.
These two parts are linked together by the so-called Multifractal Formalism.

\medskip

{\it Fine multifractal analysis.}
The first main ingredient in multifractal analysis is the  
multifractal spectrum.
The multifractal spectrum is defined as follows.
For a Borel measure $\mu$ on $\Bbb R^{d}$ 
and
a positive number $\alpha$,
let us consider
the set  of
those points
$x$ in $\Bbb R^{d}$ for which the measure
$\mu(B(x,r))$ of the ball
$B(x,r)$ with center $x$ and radius $r$ behaves like
$r^{\alpha}$ for small $r$,
i\.e\. the set
 $$
   \Bigg\{
   x\in K
   \,\Bigg|\,
   \lim_{r\searrow 0}
   \frac{\log\mu(B(x,r))}{\log r}
   =
   \alpha
   \Bigg\}\,.
   \tag1.1
 $$
If the intensity of the measure $\mu$ varies very widely, it may
happen that the sets in (1.1)
display a
fractal-like character for a range of values of $\alpha$. In this case
it is natural to study
the Hausdorff dimensions of the sets in (1.1)
 as $\alpha$
varies.
We therefore define the  fine multifractal spectrum of $\mu$
by
 $$
 f_{\mu}(\alpha)
 =
 \dim_{\Haus}
   \Bigg\{
   x\in K
   \,\Bigg|\,
   \lim_{r\searrow 0}
   \frac{\log\mu(B(x,r))}{\log r}
   =
   \alpha
   \Bigg\}\,,
 \tag1.2
$$
where $\dim_{\Haus}$ denotes the Hausdorff dimension.
 It is clear that the fine multifractal spectrum is 
 defined using  local 
 (i\.e\. {\it fine}) properties of the measure, namely, the local dimension.

 \medskip
 
 {\it Coarse multifractal analysis.}
The second main ingredient in multifractal analysis  is the  
Renyi dimensions.
The Renyi
dimensions are defined 
using global (i\.e\. {\it coarse}) 
properties
of the measure,
namely,
they quantify
the varying intensity of a measure by analyzing its moments at different scales. 
Formally, 
for $q\in\Bbb R$,
the $q$'th
Renyi dimensions 
$\tau_{\mu}(q)$
of $\mu$ is defined by
 $$
 \tau_{\mu}(q)
 =
 \lim_{r\searrow 0}
 \frac{\dsize \log\int\limits_{\,\,\,K}\mu(B(x,r))^{q-1}\,d\mu(x)}{-\log r}\,,
 \tag1.3
 $$
provided the limit exists.

\medskip

{\it The Multifractal Formalism: linking  fine and coarse multifractal analysis.}
One of the main problems
in multifractal analysis is to 
 understand the
 multifractal spectrum and the Renyi dimensions,
 and
their relationship with each other.
Indeed,
based on a remarkable insight
together with a clever heuristic argument,
theoretical physicists
Halsey et al\. [HaJeKaPrSh]
suggested
 in the 1980's
that the multifractal spectrum
can be computed
from the Renyi dimensions
using a principle known
as the Multifractal Formalism.
More precisely,
the 
Multifractal Formalism
predicts  that
the
multifractal spectrum 
equals the 
Legendre transform of the 
Renyi dimensions.
Before stating this formally,
we remind the reader
that if
$X$ is an inner product space with inner product $\langle\cdot |\cdot\rangle$
and
$f:X\to\Bbb R$ is a real valued function,
then the Legendre transform $f^{*}:\Bbb R\to[-\infty,\infty]$
of $f$ is defined by
 $$
 f^{*}(x)
 =
 \inf_{y}(\langle x|y\rangle+f(y))\,.
 \tag1.4
 $$
We can now state a
mathematically precise version of
the 
Multifractal Formalism.

\bigskip

\proclaim{Definition. The Multifractal Formalism for Measures}
Let  $\mu$ be a Borel measure on $\Bbb R^{d}$
with support $K$.
We will say that
$\mu$
satisfies 
the 
Multifractal Formalism
if the
following conditions hold:
\roster
\item"(i)"
The limit in (1.3) exists for all $q$.
\item"(ii)"
The function $\tau_{\mu}$ is convex and differentiable.
\item"(iii)"
$\tau_{\mu}(0)=\dim_{\Haus}K$.
\item"(iv)"
The
multifractal spectrum
of $\mu$
equals the 
Legendre transform
of
the Renyi dimensions, i\.e\.
 $$
 f_{\mu}(\alpha)
 =
\tau_{\mu}^{*}(\alpha)
\,\,\,\,
\text{for all $\alpha$.}
 $$
 \endroster
\endproclaim

\bigskip

\noindent
Conditions (ii) and (iii) are usually not included
in the 
Multifractal Formalism;
however, since these properties
play an important role later,
we have decided to include them.
During the past 20 years
there has been an enormous interest 
in
verifying the
Multifractal Formalism 
and
computing the multifractal spectra of measures
in 
the mathematical literature.
In the mid 1990's
Cawley \& Mauldin [CaMa] and Arbeiter \& Patzschke [ArPa]
verified the Multifractal Formalism for self-similar measures 
satisfying the so-called  Open Set Condition,
and within the last
20 years the multifractal spectra of various classes of measures
in Euclidean space $\Bbb R^{d}$ 
exhibiting some degree of self-similarity have been computed 
rigorously, cf\. 
the textbooks [Fa2,Pe2]
and the references therein.
On the other hand, it is also known that most 
measures
(both in the sense of Baire category and in the sense of 
\lq\lq shyness")
do not satisfy the 
Multifractal Formalism, see, for example, [Bay1,Bay2,Ol3,Ol4].
While most 
measures
do not satisfy the 
Multifractal Formalism,
it is well-known and not difficult to 
show that
the multifractal spectrum is always bounded above by the Legendre transform of the
Renyi dimension.
Since this result plays an important role later,
we have decided to
state it formally.

\bigskip

\proclaim{Proposition A [LauNg,Ol1,Pe1]}
Let  $\mu$ be a Borel measure on $\Bbb R^{d}$
and assume that
the limit in (1.3) exists
for all $q$.
Then
 $$
 f_{\mu}(\alpha)
 \le
 \tau_{\mu}^{*}(\alpha)
 \,\,\,\,
\text{for all $\alpha$.}
 $$
\endproclaim

\bigskip

More generally,
if $K\subseteq\Bbb R^{d}$ and $\Phi:K\to\Bbb R$ is a function,
then the  fine 
multifractal spectrum
of $\Phi$
is defined by
 $$
 f_{\Phi}(\alpha)
 =
 \dim_{\Haus}
 \Big\{
 x\in K
 \,\Big|\,
 \Phi(x)=\alpha
 \Big\}\,.
 \tag1.5
 $$
Of course, if
$\Phi(x)
=
\lim_{r\searrow 0}\frac{\log\mu(B(x,r))}{\log r}$,
then
the multifractal spectrum of $\Phi$ equals the 
multifractal spectrum of $\mu$ from (1.2).
Other interesting and important examples are obtained by,
for example,
 letting $\Phi(x)$
equal the ergodic Birkhoff average of a continuous function  at $x$,
or by letting
$\Phi(x)$
equal the local entropy
of a dynamical system at $x$.
In this more general setting
it is also natural to 
attempt to find a 
\lq\lq coarse moment scaling"
function
$\tau_{\Phi}:\Bbb R\to\Bbb R$
such that
$f_{\Phi}(\alpha)=\tau_{\Phi}^{*}(\alpha)$
for all $\alpha$,
and during the past 10 years
an extensive theory
analysing  this problem
have been developed,
 see, for example, the texts by Pesin [Pe2] and
Barrierra [Bar].

\bigskip

{\bf Multifractal zeta-functions.}
Dynamical zeta-functions were introduced by Artin \& Mazur
in the mid 1960's [ArMa]
based on an analogy with the number theoretical zeta-functions
associated with a function field
over a finite ring.
Subsequently Ruelle [Rue1,Rue2]
associated
zeta-functions
to 
certain statistical mechanical models in one dimensions.
During the past 35 years 
many parallels have been 
drawn between
number theory 
zeta-functions, 
dynamical zeta-functions,
and statistical mechanics zeta-functions.
However, much more recently and
motivated by the 
powerful techniques
provided by 
the use of the Artin-Mazur zeta-functions in number theory
and the use of the Ruelle
zeta-functions in dynamical systems,
Lapidus and collaborators
(see the intriguing books by Lapidus \& van Frankenhuysen [Lap-vF1,Lap-VF2]
and the references therein)
have recently
introduced and pioneered
to use of
zeta-functions in fractal geometry.
Inspired by this development,
within the past 4--5 years
several authors have 
paralleled this development
by
introducing
zeta-functions
into  multifractal geometry.

\medskip

{\it Fine multifractal zeta-functions.}
%
%
%
Indeed, in 2009,
Lapidus and collaborators
[LapRo,LapLe-VeRo] introduced 
various intriguing 
 multifractal zeta-functions 
designed to provide information about the multifractal spectrum
$f_{\mu}(\alpha)$
of 
self-similar measures $\mu$,
and
a number of connections with multifractal spectra were
suggested
and in some cases proved;
for example, in simplified cases 
the multifractal spectrum of a self-similar measure could be recovered from a 
zeta-function.
The 
key idea in [LapRo,LapLe-VeRo] is both simple and attractive:
while traditional zeta-functions are defined  by 
\lq\lq summing over all data",
the
multifractal zeta-functions
in
[LapRo,LapLe-VeRo] 
are defined by only
\lq\lq summing over data that are multifractally relevant".
%
%
Ideas similar to those in
[LapRo,LapLe-VeRo]
have  much more recently
been revisited and
investigated 
in
[Bak,MiOl1,MiOl2,Ol7] where
the authors
introduce related
multifractal zeta-functions 
tailored to study
the
multifractal spectra
of self-conformal measures
and a number of connections with
very general types of multifractal spectra were established.
Indeed,  in [MiOl2,Ol7] 
we proposed  a 
theory of
{\it fine} 
multifractal zeta-functions
paralleling
the existing
theory
 of 
dynamical zeta-functions 
 introduced and developed by Ruelle [Rue1,Rue2]
and others [Bal1,Bal2,ParPo1,ParPo2].
In particular,
in the setting of graph-directed
self-conformal constructions,
we introduced 
a family
of
{\it fine}
 multifractal zeta-functions
designed
to
 provide
precise information of very general classes of multifractal spectra, 
including, for example, the multifractal 
spectra of self-conformal measures and the multifractal spectra of ergodic 
Birkhoff averages 
of continuous functions on self-conformal sets.
More 
precisely, 
we fix a
metric space $X$ and
a continuous map
$U:\Cal P(\Sigma_{\smallG}^{\Bbb N})\to X$ 
where
$\Sigma_{\smallG}^{\Bbb N}$
denotes the
shift space 
modelling the underlying 
graph-directed
self-conformal construction $\G$
(the precise definition of  $\Sigma_{\smallG}^{\Bbb N}$ will be given in Section 2)
and
$\Cal P(\Sigma_{\smallG}^{\Bbb N})$
denotes the
family of probability measures
on $\Sigma_{\smallG}^{\Bbb N}$.
For each 
subset $C$ of $X$
and each
 continuous function
 $\varphi:\Sigma_{\smallG}^{\Bbb N}\to\Bbb R$,
we now associated
a
{\it fine}
dynamical  multifractal zeta-function
$\zeta_{C}^{\dyn,U}(\varphi;z)$
defined by
 $$
\zeta_{C}^{\dyn,U}(\varphi;z)
=
   \sum_{n}
   \,\,
   \frac{z^{n}}{n}
 \left(
   \sum
  \Sb
 \bold i\in\Sigma_{\smallsmallG}^{n}\\
  {}\\
  UL_{n}[\bold i]\subseteq C
  \endSb
  \,\,
  \sup_{\bold u\in[\bold i]}
  \,\,
 \exp
 \,\,
 \sum_{k=0}^{n-1}\varphi (S^{k}\bold u)
\right)\,,
\tag1.6
$$
where
$z$ is a complex variable;
in (1.6) we have used the following
notation, namely,
$S$ denotes the shift map,
$L_{n}$ denotes the $n$'th empirical measure
(i\.e\.
$L_{n}\bold i=\frac{1}{n}\sum_{k=0}^{n-1}\delta_{S^{k}\bold i}$ for 
$\bold i\in\Sigma_{\smallG}^{\Bbb N}$),
$\Sigma_{\smallG}^{n}$ denotes the 
family of all 
permissible
words of length $n$
and if $\bold i\in\Sigma_{\smallG}^{n}$,
then $[\bold i]$ denotes the cylinder generated by $\bold i$;
the precise definitions will be given in Section 2.
One of main results
from [MiOl2,Ol7]
shows that if $\alpha$ is a point in $X$
and $\varphi<0$,
then there is a unique
real number 
$\,\,\scri f(\alpha)$
such that
 $$
 \lim_{r\searrow 0}
 \sigma_{\radius}
 \big(
 \,
 \zeta_{B(\alpha,r)}^{\dyn,U}(\,\,\,\scri f(\alpha)\,\varphi;\cdot)
 \,
 \big)
 =
 1\,,
 $$
where
 $\sigma_{\radius}$ denotes the radius of convergence;
furthermore,
$\,\,\scri f(\alpha)$ 
satisfies a variational principle
and
for judicious 
choices of
$X$, $U$ and $\varphi$,
the number
$\,\,\scri f(\alpha)$ equals 
important fine multifractal spectra,
see, [MiOl2,Ol7] for a large number of
examples 
illustrating this.

\medskip

{\it Coarse multifractal zeta-functions.}
Because of the 
importance of the Renyi dimensions, it is natural to
construct a
theory of
{\it coarse} 
multifractal zeta-functions
paralleling
the existing
theory
 of 
fractal zeta-functions introduced and pioneered
by Lapidus \& van Frankenhuysen [Lap-vF1,Lap-VF2].
This has recently been done
by
 Levy-Vehel \& Mendivil
[Le-VeMe]
who
 introduced
 multifractal zeta-functions
tailored to
provide information about the 
Renyi dimensions 
$\tau_{\mu}(q)$
of self-similar measures $\mu$.
Ideas related to those in [Le-VeMe]
have also been investigated
in [Ol5,Ol6]
where related
multifractal zeta-functions 
designed
 to study
the 
Renyi dimensions
and the (closely related) multifractal  Minkowski volume
of self-conformal measures
are introduced and investigated.

While Renyi dimensions play an important role in multifractal analysis,
the Multifractal Formalism 
underpinnning the 
link between the Renyi dimensions and the multifractal spectrum
is 
fundamental for
a fuller understanding
of
 the Renyi dimensions
 as well as
the 
genuinely
mysterious
relationship 
between these 
and the multifractal spectrum.
In recognition
 of this viewpoint,
 and in order to place
 the somewhat ad-hoc 
 multifractal zeta-functions from
 [Le-VeMe,Ol5,Ol6]
in to a broader context,
it is natural to ask if
the 
\lq\lq fine multifractal spectrum"
$\,\,\scri f(\alpha)$
defined using the 
{\it fine} multifractal zeta-functions
satisfies a multifractal formalism
based on
a natural family of {\it coarse} multifractal zeta-functions.
More precisely,
assuming $X$ is an inner product space 
(indeed, 
since
the classical Multifractal Formalism
involves the Legendre transform,
it seems natural to 
work in a setting which
allows Legendre transforms,
and the 
most natural setting for this in
within the context of inner product spaces),
we pose the following question: 
is it possible to associated
a 
natural
{\it coarse} 
multifractal zeta-function
$\zeta_{q}^{\co,U}(\varphi;s)$,
where
 $s$ is a complex variable,
 to each 
 continuous function
 $\varphi:\Sigma_{\smallG}^{\Bbb N}\to\Bbb R$
 and each point $q\in X$, 
 such that
if we let
 $$
 \tau(q)
 =
 \sigma_{\abs}
 \big(
 \,
 \zeta_{q}^{\co,U}(\varphi;\cdot)
 \,
 \big)
$$
where
$ \sigma_{\abs}$ denotes 
the abscissa of convergence, then
 $$
\,\,\scri f(\alpha)
=
\tau^{*}(\alpha)
\,\,\,\,
\text{for all (or some) $\alpha$?}
$$
Adopting this viewpoint,
the purpose 
of this paper is to 
introduce and 
develop
a meaningful
theory of
{\it coarse}
multifractal zeta-functions in the setting of
graph-directed self-conformal constructions
and
establish 
a multifractal formalism for zeta-functions relating 
the
{\it fine}
multifractal zeta-function in (1.6)
and the
{\it fine}
multifractal zeta-function in (1.7) below via 
the Legendre transform.
The 
{\it coarse} 
multifractal zeta-function
$\zeta_{q}^{\co,U}(\varphi;s)$
is defined 
by
 $$
 \align
 \zeta_{q}^{\co,U}(\varphi;s)
&=
\sum_{n}
\,\,
 \sum_{\bold i\in\Sigma_{\smallsmallG}^{n}}
 \,\,
 \exp
 \Bigg(
 \Bigg(
\sup_{\bold u\in[\bold i]}\big\langle q\big |UL_{n}\bold u\, \big\rangle
+s
\Bigg)
 \sup_{\bold u\in[\bold i]}
 \sum_{k=0}^{n-1}\varphi(S^{k}\bold u)
 \Bigg)
 \tag1.7
 \endalign 
$$
where $s$ is complex 
variable and where we have used the same notation as in (1.6).

We illustrate 
the definition of
$\zeta_{q}^{\co,U}(\varphi;s)$ by considering the following simple example
involving self-similar measures.
Self-similar measures form a special case of the 
family of 
graph-directed self-conformal measures and are defined as follows.
Let
$S_{1},\ldots,S_{N}:\Bbb R^{d}\to\Bbb R^{d}$ be contraction similarities
and let $r_{i}$ denote the contraction ratio of $S_{i}$, i\.e\.
$|S_{i}x-S_{i}y|=r_{i}|x-y|$
for all $x,y\in\Bbb R^{d}$.
Also, fix a probability vector  $(p_{1},\ldots,p_{N})$.
The self-similar set $K$ and the self-similar measure $\mu$
associated with the list
$(S_{1},\ldots,S_{N},p_{1},\ldots,p_{N})$
 is the unique set and the unique probability measure satisfying
 $$
 K
 =
 \bigcup_{i}S_{i}(K)\,,\,\,\,\,
 \mu
 =
 \sum_{i}p_{i}\mu\circ S_{i}^{-1}\,.
 \tag1.8
 $$
In this particular case, 
the shift space $\Sigma_{\smallG}^{\Bbb N}$
is given by
$\Sigma_{\smallG}^{\Bbb N}=\{1,\ldots,N\}^{\Bbb N}$.
We now  define maps
 $\varphi,\psi:\{1,\ldots,N\}^{\Bbb N}\to\Bbb R$ and
$U:\Cal P(\{1,\ldots,N\}^{\Bbb N})\to\Bbb R$
by
 $$
 \align
 \varphi(\bold i)
&=
 \log r_{i_{1}}\,,\\
  \psi(\bold i)
&=
 \log p_{i_{1}}
 \endalign
 $$
for $\bold i=i_{1}i_{2}\ldots\in \{1,\ldots,N\}^{\Bbb N}$
and
 $$
 U\mu
 =
 \frac{\int\psi\,d\mu}{\int\varphi\,d\mu}
 $$
for $\mu\in \Cal P(\{1,\ldots,N\}^{\Bbb N})$.
For this choice of $\varphi$ and $U$, a simple calculation shows that
(the reader is referred to Section 5 for details
and for several other examples) 
$$
 \align
 \zeta_{q}^{\co,U}(\varphi;s)
&=
\sum_{n}
\,\,
 \sum_{\bold i=i_{1}\ldots i_{n}\in\{1,\ldots,N\}^{n}}
 \,\,
 \exp
 \Bigg(
 \Bigg(
\sup_{\bold u\in[\bold i]}q \,U(L_{n}\bold u)
+s
\Bigg)
 \sup_{\bold u\in[\bold i]}
 \sum_{k=0}^{n-1}\varphi(S^{k}\bold u)
 \Bigg)\\
&= 
\sum_{n}
\,\,
 \sum_{\bold i=i_{1}\ldots i_{n}\in\{1,\ldots,N\}^{n}}
 \,\,
 \exp
 \Bigg(
 \Bigg(
q\,\frac{\sum_{k=0}^{n-1}\log p_{i_{k}}}{\sum_{k=0}^{n-1}\log r_{i_{k}}}+s
\Bigg)
 \sum_{k=0}^{n-1}\log r_{i_{k}}
 \Bigg)\\
&= 
\sum_{n}
\,\,
 \sum_{\bold i=i_{1}\ldots i_{n}\in\{1,\ldots,N\}^{n}}
 \,\,
p_{i_{1}}^{q}\cdots p_{i_{n}}^{q}\,r_{i_{1}}^{s}\cdots r_{i_{n}}^{s}\\ 
&=
\sum_{n}
\,\,
\Bigg(
 \sum_{i=1}^{N}
p_{i}^{q} r_{i}^{s}
\Bigg)^{n}\\ 
&=
\frac
{
\sum_{i=1}^{N}
p_{i}^{q} r_{i}^{s}
}
{1- 
 \sum_{i=1}^{N}
p_{i}^{q} r_{i}^{s}
}
\tag1.9
 \endalign 
$$
provided the geometric series 
$\sum_{n}
(
 \sum_{i=1}^{N}
p_{i}^{q} r_{i}^{s}
)^{n}$
converges,
i\.e\.
provided
$|\sum_{i=1}^{N}p_{i}^{q} r_{i}^{s}|<1$.
The zeta-function in (1.9), i\.e\. the zeta-function
 $$
\zeta_{q}(s)
=
 \frac
{
\sum_{i=1}^{N}
p_{i}^{q} r_{i}^{s}
}
{1-
 \sum_{i=1}^{N}
p_{i}^{q} r_{i}^{s}
}\,,
\tag1.10
$$
and similar zeta-functions have, in fact, been investigated
intensively during the past 10 years
by Lapidus et at
[[Lap-vF1,Lap-VF2]
(for $q=0$)
and others
[Le-VeMe,Ol5,Ol6].
For example,
it is well-known 
that
if $K$ and $\mu$ denote the self-similar set and the self-similar measure satisfying (1.8),
then the  
Renyi dimension $\tau_{\mu}(q)$ of $\mu$ (see (1.3))
equals the 
abscissa of convergence of 
the zeta-function in (1.10), see [Le-VeMe,Ol5,Ol6]
(and [Lap-vF1,Lap-VF2] for the case $q=0$).
In addition,
explicit formulas for the
Minkowski volume of the fractal $K$ and multifractal $\mu$
can be found
from the poles and the residues of the 
zeta-function in (1.10); see
Lapidus \& van Frankenhuysen's books [Lap-vF1,Lap-VF2]
for a detailed discussion of this in the fractal case
(i\.e\. the the case $q=0$)
and Olsen [Ol5,Ol6]
for a discussion of the multifractal case
(i\.e\. the case $q\in \Bbb R$).

The 
{\it coarse}
multifractal zeta-functions
introduced in this paper 
may therefore be view as
a continuation
of the work initiated in 
[Le-VeMe,Ol5,Ol6]
placing the zeta-functions from
[Le-VeMe,Ol5,Ol6]
in to a broader context
provider a fuller understanding
and
allowing
further applications 
(see, in particular, Section 5
where applications to 
multifractal analysis of graph-directed self-conformal measures 
and
multifractal analysis of various ergodic spectra are considered).

\medskip

{\it The Multifractal Formalism for zeta-functions:
linking fine and coarse multi fractal zeta-functions.}
For judicious 
choices of
$X$, $U$ and $\varphi$,
the abscissa of convergence
 $$
 \tau(q)
 =
 \sigma_{\abs}
 \big(
 \,
 \zeta_{q}^{\co,U}(\varphi;\cdot)
 \,
 \big)
$$
of the coarse multifractal zeta-function 
$\zeta_{q}^{\co,U}(\varphi;\cdot)$
has many of the properties
that
any meaningful
\lq\lq coarse moment scaling" function
is expected to have,
including the following
(see the statement of the Multifractal Formalism for Measures and Proposition A):
\roster
\item"$\bullet$"
the function
$\tau$ is convex and 
(under some mild assumptions) differentiable;
\item"$\bullet$"
the
\lq\lq fine multifractal spectrum"
 $\,\,\scri f(\alpha)$
 is {\it always}
bounded above by the 
Legendre transform of
$\tau$;
\item"$\bullet$"
if $\tau$ is differentiable, then
 $\,\,\scri f(\alpha)$ equals the 
 Legendre transform of
$\tau$, i\.e\.
the Multifractal Formalism is satisfied.
\endroster
This is the content of the 
Theorem 1.1 below;
we emphasise that Theorem 1.1
is special  case of 
the 
three main results in this paper, namely,
Theorem 4.2, Theorem 4.3 and Theorem 4.4,
where more general results are presented.

\bigskip

\proclaim{Theorem 1.1. The multifractal formalism for multifractal zeta-functions; a special case}
Let $X$ be an inner product space and let $U:\Cal P(\Sigma_{\smallG}^{\Bbb N})\to X$ be 
continuous with respect to the weak topology.
Fix a 
H\"older continuous function  $\varphi:\Sigma_{\smallG}^{\Bbb N}\to\Bbb R$ 
with
$\varphi<0$.
For
$C\subseteq X$ and $q\in X$,
we let 
$ \zeta_{C}^{\dyn,U}(\varphi;\cdot)$
and
$ \zeta_{q}^{\co,U}(\varphi;\cdot)$
denote the fine and coarse multifractal zeta-functions in (1.6) and (1.7),
respecively,
and let
 $$
 \tau(q)
 =
 \sigma_{\abs}
 \big(
 \,
 \zeta_{q}^{\co,U}(\varphi;\cdot)
 \,
 \big)\,.
$$
\roster
\item"(1)"
The function $\tau$ is convex
and if there is a 
H\"older continuous function  $\psi:\Sigma_{\smallG}^{\Bbb N}\to\Bbb R$ 
such that
$U\mu=\frac{\int\psi\,d\mu}{\int\varphi\,d\mu}$
for all $\mu$,
then $\tau$ is real analytic.
\item"(2)"
Fix $\alpha$ be a  point in $X$
and let
$\,\,\scri f(\alpha)$
be the unique real number
such that
 $$
 \align
 \limsup_{r\searrow 0}
 \sigma_{\radius}
 \big(
 \,
 \zeta_{B(\alpha,r)}^{\dyn,U}(\,\,\scri f(\alpha)\,\varphi;\cdot)
 \,
 \big)
 &=
 1\,.
 \endalign
 $$ 
Then
 $$
  \scri f(\alpha)
  \le
  \tau^{*}(\alpha)\,.
  $$
If, in addition, $X=\Bbb R^{M}$
and 
there is $q\in\Bbb R^{M}$ such that $\tau$ is differentiable 
at $q$ with
$\alpha=-\nabla\tau(q)$, then
$$
  \scri f(\alpha)
  =
  \tau^{*}(\alpha)\,.
  \tag1.11
  $$
(Observe that since $\tau$ is convex,
we conclude that $\tau$
is differentiable almost everywhere,
and the conclusion in (1.11)
is therefore 
satisfied for 
\lq\lq many"
points $\alpha$.)  
\endroster

\endproclaim

  \bigskip

\noindent
For judicious 
choices of $X$, $U$ and $\varphi$
equation (1.11)
reduces to
known Legendre transform representations of
different types of classical
multifractal spectra,
and for other
choices of $X$, $U$ and $\varphi$
equation (1.11)
provides
new Legendre transform representations of
multifractal spectra;
for example,
for certain choices of $X$, $U$ and $\varphi$
equation (1.11) equals the
well-know 
Legendre transform representation of the
multifractal spectrum of a graph directed self-conformal measure
(see the examples in Section 5.1),
and for 
other
choices of $X$, $U$ and $\varphi$
equation (1.11) provides 
new
Legendre transform representations of the
multifractal spectrum of 
different types of 
 multifractal spectra of ergodic averages of continuous functions on
 graph-directed
 self-conformal sets
 (see the examples in Section 5.2).

%
%
%

The remaining part
of the paper is organised as follows.
In Section 2 we briefly recall the 
definitions
of self-conformal constructions
and in Section 3 we recall the definition of the fine 
dynamical multifractal zeta-functions
introduced in 
[MiOl2,Ol7].
In Section 4
we provide the definition of 
the coarse multifractal zeta-functions 
 and we 
 state our main 
  results.
 Section 5 contains a number of examples
 illustrating the main results.
 Finally, the proofs are presented in Sections 6--9.

\newpage

\bigskip

%
%
%

\centerline{\smc 2. The setting:}

\centerline{\smc
Graph-directed
self-conformal sets and 
graph-directed
self-conformal measures.}

\medskip

{\bf Notation from symbolic dynamics.}
We first recall the 
 notation
and terminology  from symbolic dynamics 
that will be used in this paper.
Fix
a finite directed multigraph
$\G=(\V,\E)$
where
$\V$ denotes the set of vertices of $\G$ and
$\E$ denotes the set of edges of $\G$.
For an edge $\edge\in\E$,
we 
write
$\ini(\edge)$
for the initial vertex of $\edge$ 
and we 
write
$\termi(\edge)$ for the terminal vertex of $\edge$.
For $\vertexi,\vertexj\in\V$, write
 $$
 \aligned
 \E_{\smallvertexi}
&=
\Big\{\edge\in\E\,\Big|\,\ini(\edge)=\vertexi\Big\}\,,\\
 \E_{\smallvertexi,\smallvertexj}
&=
\Big\{\edge\in\E\,\Big|\,\ini(\edge)=\vertexi\,\,\text{and}\,\,\termi(\edge)=\vertexj\Big\}\,;
\endaligned
\tag2.1
 $$  
i\.e\.
 $\E_{\smallvertexi}$ is the family of all edges starting at $\vertexi$;
 and
 $ \E_{\smallvertexi,\smallvertexj}$
 is the 
 is the family of all edges starting at $\vertexi$
 and ending at $\vertexj$.
Also, for a positive integer $n$,
we
 write
 $$
 \aligned
 \Sigma_{\smallG}^{n}
&=
\Big\{
\edge_{1}\ldots\edge_{n}
\,\Big|\,
\edge_{i}\in\E\,\,\text{for $1\le i\le n$,}\\
&\qquad\qquad
   \quad\,\,\,\,\,\,
\termi(\edge_{1})=\ini(e_{2})\,,\\
&\qquad\qquad
   \quad\,\,\,\,\,\,
\termi(e_{i-1})=\ini(\edge_{i})
\,\,\text{and}\,\,
\termi(\edge_{i})=\ini(e_{i+1})
\,\,\text{for $1< i < n$,}\\
&\qquad\qquad
   \quad\,\,\,\,\,\,
\termi(e_{n-1})=\ini(\edge_{n})
\Big\}\\
 \Sigma_{\smallG}^{*}
&=
\bigcup_{m}\,\Sigma_{\smallG}^{m}\,,\\
 \Sigma_{\smallG}^{\Bbb N}
&=
\Big\{
\edge_{1}\edge_{2}\ldots\,
\,\Big|\,
\edge_{i}\in\E\,\,\text{for $1\le i$,}\\
&\qquad\qquad
   \quad\,\,\,\,\,\,
\termi(\edge_{1})=\ini(e_{2})\,,\\
&\qquad\qquad
   \quad\,\,\,\,\,\,
\termi(e_{i-1})=\ini(\edge_{i})
\,\,\text{and}\,\,
\termi(\edge_{i})=\ini(e_{i+1})
\,\,\text{for $1< i $}
\Big\}\,;
 \endaligned
 \tag2.2
 $$
i\.e\. $\Sigma_{\smallG}^{n}$ is the family of all
finite strings
$\bold i=\edge_{1}\ldots \edge_{n}$
consisting of finite paths in $\G$
of length $n$;
$\Sigma_{\smallG}^{*}$ is the family of all finite strings
$\bold i=\edge_{1}\ldots \edge_{m}$
with $m\in\Bbb N$ consisting of finite paths in $\G$;
and 
$\Sigma_{\smallG}^{\Bbb N}$ is the family of all
infinite
strings
$\bold i=\edge_{1}\edge_{2}\ldots $
consisting of infinite paths in $\G$.
For 
a finite string 
$\bold i=\edge_{1}\ldots \edge_{n}\in\Sigma_{\smallG}^{n} $, we write
 $$
 \ini(\bold i)
 =
 \ini(\edge_{1})\,,\,\,\,\,
 \termi(\bold i)
 =
 \termi(\edge_{n})\,,
 \tag2.3
 $$
and for an infinite string
$\bold i=\edge_{1} \edge_{2}\ldots\in\Sigma_{\smallG}^{\Bbb N} $, we write
 $$
 \ini(\bold i)
 =
 \ini(\edge_{1})\,.
 \tag2.4
 $$
Next,
for an infinite string 
$\bold i=\edge_{1}\edge_{2}\ldots\in\Sigma_{\smallG}^{\Bbb N} $
and a positive integer $n$, we will write
$\bold i|n
=
\edge_{1}\ldots \edge_{n}$.
In addition, for
a positive integer $n$
and
a finite string 
$\bold i=\edge_{1}\ldots \edge_{n}\in\Sigma_{\smallG}^{n} $
with length equal to  $n$,
 we will write
$|\bold i|
=
n$, and we let $[\bold i]$ denote the cylinder 
generated by $\bold i$, i\.e\.
  $$
  [\bold i]
 =
 \Big\{
 \bold j\in\Sigma_{\smallG}^{\Bbb N}
 \,\Big|\,
 \bold j|n=\bold i
 \Big\}\,.
 \tag2.5
 $$
Finally, let $S:\Sigma_{\smallG}^{\Bbb N}\to\Sigma_{\smallG}^{\Bbb N}$ denote the shift map, i\.e\.
 $$
 S(\edge_{1}\edge_{2}\ldots)
 =
 \edge_{2}\edge_{3}\ldots\,.
 $$

\bigskip

{\bf Graph-directed
self-conformal sets and 
graph-directed 
self-conformal measures.}
Next, we recall the definition of
graph-directed
 self-conformal 
sets and measures. 
A graph-directed conformal iterated function system with probabilities 
is a list
 $$
 \Big(\,
 \V,\,
 \E,\,
 (V_{\smallvertexi})_{\smallvertexi\in\smallV},\,
 (X_{\smallvertexi})_{\smallvertexi\in\smallV},\,
 (S_{\smalledge})_{\smalledge\in\smallE},\,
 (p_{\smalledge})_{\smalledge\in\smallE}\,
 \Big)
 $$
where

\bigskip

\roster
\item"$\bullet$"
For each $\vertexi\in\V$
we have:
$V_{\smallvertexi}$ is an open, connected subset of $\Bbb R^{d}$.
\item"$\bullet$" 
For each $\vertexi\in\V$
we have:
$X_{\smallvertexi}\subseteq V_{\smallvertexi}$ is a compact set with
$X_{\smallvertexi}^{\circ\,-}=X_{\smallvertexi}$.
\item"$\bullet$" 
For each $\vertexi,\vertexj\in\V$
and
$\edge\in E_{\smallvertexi,\smallvertexj}$
we have:
$S_{\smalledge}:V_{\smallvertexj}\to V_{\smallvertexi}$
is a contractive
$C^{1+\gamma}$ diffeomorphism with
$0<\gamma<1$
such that
$S_{\smalledge}(X_{\smallvertexj})\subseteq X_{\smallvertexi}$.
\item"$\bullet$"
The Conformality Condition.
For $\vertexi,\vertexj\in\V$,
$\edge\in E_{\smallvertexi,\smallvertexj}$ and $x\in V_{\smallvertexj}$, let
$
(DS_{\smalledge})(x):\Bbb R^{d}\to\Bbb R^{d}
$
denote the derivative of $S_{\smalledge}$ at $x$. 
For each $\vertexi,\vertexj\in\V$
and
$\edge\in E_{\smallvertexi,\smallvertexj}$
we have:
$(DS_{\smalledge})(x)$ is a contractive similarity map, i\.e\.
there exists
$s_{\smalledge}(x)\in(0,1)$ such that
$|(DS_{\smalledge})(x)u-(DS_{\smalledge})(x)v|
 =
 s_{\smalledge}(x)|u-v|$
for all $u,v\in\Bbb R^{d}$.
\item"$\bullet$" 
For each $\vertexi\in\V$
we have:
$(p_{\smalledge})_{\smalledge\in\smallE_{\smallsmallvertexi}}$ is a probability vector.
\endroster

\bigskip

\noindent
It follows from [Hu] that there exists a unique
list $(K_{\smallvertexi})_{\smallvertexi\in\smallV}$ of
non-empty compact sets $K_{\smallvertexi}\subseteq X_{\smallvertexi}$ such that
 $$
 K_{\smallvertexi}
 \,=\,
 \bigcup_{\smalledge\in \smallE_{\smallsmallvertexi}}\,
 S_{\smalledge}K_{\termi(\smalledge)}\,,
 \tag2.6
 $$
and a unique
list $(\mu_{\smallvertexi})_{\smallvertexi\in\smallV}$ 
of probability measures with
$\supp\mu_{\smallvertexi}=K_{\smallvertexi}$ such that
 $$
 \mu_{\smallvertexi}
 \,=\,
 \sum_{\smalledge\in E_{\smallsmallvertexi}}\,
 p_{\smalledge}\,\mu_{\termi(\smalledge)}\circ S_{\smalledge}^{-1}\,.
 \tag2.7
 $$
The sets $(K_{\smallvertexi})_{\smallvertexi\in\smallV}$
and measures $(\mu_{\smallvertexi})_{\smallvertexi\in\smallV}$
are called the self-conformal sets and self-conformal
measures associated with the list (1.21), respectively.
We will frequently assume that the so-called Open Set condition (OSC) 
 is satisfied. The OSC is defined as follows:

\bigskip

\roster
\item"$\bullet$"
The Open Set Condition:
There exists a list
$(U_{\smallvertexi})_{\smallvertexi\in\smallV}$ of open non-empty and bounded sets
$U_{\smallvertexi}\subseteq X_{\smallvertexi}$ with 
$S_{\smalledge}\big(\,U_{\smallvertexj}\,\big)
 \subseteq
 U_{\smallvertexi}$
for all $\vertexi,\vertexj\in\V$ and all $\edge\in\E_{\smallvertexi,\smallvertexj}$ such that
 $S_{\smalledge_{1}}\big(\,U_{\termi(\smalledge_{1})}\,\big)
 \,\cap\,
 S_{\smalledge_{2}}\big(\,U_{\termi(\smalledge_{2})}\,\big)
 =
 \varnothing$
for all $\vertexi\in\V$ and all $\edge_{1},\edge_{2}\in\E_{\smallvertexi}$ with 
$\edge_{1}\not=\edge_{2}$.
\endroster 
 
\bigskip

\noindent
For $\bold i=\edge_{1}\ldots\edge_{n}\in\Sigma_{\smallG}^{n}$, we write
 $$
 \aligned
 S_{\bold i}
&=
 S_{\smalledge_{1}}\cdots S_{\smalledge_{n}}\,,\\
 K_{\bold i}
&=
 S_{\smalledge_{1}}\cdots S_{\smalledge_{n}}
 \big(K_{\termi(\smalledge_{n})}\big)\,,
 \endaligned
 \tag2.8
 $$ 
and
we define the projection
$\pi:\Sigma_{\smallG}^{\Bbb N}\to \Bbb R^{d}$ by
 $$
 \big\{\,\pi(\bold i)\,\big\}
 =
 \bigcap_{n}
 \,
 K_{\bold i|n}
 \tag2.9
 $$
for $\bold i=\edge_{1}\edge_{2}\ldots\in\Sigma_{\smallG}^{\Bbb N}$.
Finally, we define
$\Lambda:\Sigma_{\smallG}^{\Bbb N}\to\Bbb R$ by
 $$
 \align
 \Lambda(\bold i)
 \,
&=\,
 \log
 \big|
 \,
 \big(DS_{\smalledge_{1}}\big)\big(\pi_{\termi(\smalledge_{1})}(S\bold i)\big)
 \,
 \big|
 \tag2.10
 \endalign
 $$
for $\bold i=\edge_{1}\edge_{2}\ldots\in\Sigma_{\smallG}^{\Bbb N}$;
loosely speaking the 
map $\Lambda$ represents the local change of scale as one goes from
$\pi_{\termi(\smalledge_{1})}(S\bold i)$ to
$\pi_{\ini(\smalledge_{1})}(\bold i)$.

\bigskip


\heading{3. Fine dynamical multifractal zeta-functions.}\endheading

The purpose of this section is to
recall the 
fine dynamical multifractal zeta-functions
introduced in [MiOl2,Ol7].
Throughout 
 this section, 
and in the remaining parts of the paper, we will used the
following notation.
 Namely,
if $(a_{n})_{n}$ is a sequence of 
complex numbers
and if $f$ is the power
series
defined by $f(z)=\sum_{n}a_{n}z^{n}$ for $z\in\Bbb C$,
then we will denote the 
radius
of convergence of $f$ by $\sigma_{\radius}(f)$, i\.e\. we 
write
 $$
 \sigma_{\radius}(f)
 =
 \text{
 \lq\lq the radius of convergence of $f$"
 }\,.
 $$

 The 
fine dynamical multifractal zeta-functions
 in [MiOl2,Ol7] are
 motivated by the 
 notion of pressure
 from the thermodynamic formalism
  and
 the
 dynamical zeta-functions introduced by
 Ruelle [Rue1,Rue2]; see, also 
 [Bal1,\allowlinebreak
 Bal2,\allowlinebreak
 ParPo1,\allowlinebreak
 ParPo2].
Because of this we
now briefly
recall the definition 
of pressure and dynamical zeta-function.
Let
$\varphi:\Sigma_{\smallG}^{\Bbb N}\to\Bbb R$
be a
continuous function.
 The
pressure
of 
 $\varphi$
 is defined 
by
 $$
 \align
 P(\varphi)
 &=
\lim_{n}
  \,\,
  \frac{1}{n}
 \,\,
 \log
   \sum
  \Sb
\bold i\in\Sigma_{\smallsmallG}^{n}
  \endSb
  \,\,
  \sup_{\bold u\in[\bold i]}
  \,\,
   \exp
 \,\,
 \sum_{k=0}^{n-1}\varphi S^{k}\bold u\,,
 \tag3.1
  \endalign
 $$
 see
 [Bo]
or
[ParPo2];
we note that
it is well-known that the limit in (3.1) exists.
Also,
the dynamical zeta-function of $\varphi$ is defined by
$$
  \zeta^{\dyn}(\varphi;z)
  =
   \sum_{n}
   \,\,
   \frac{z^{n}}{n}
 \left(
   \sum
  \Sb
\bold i\in\Sigma_{\smallsmallG}^{n}
  \endSb
  \,\,
  \sup_{\bold u\in[\bold i]}
  \,\,
 \exp
 \,\,
 \sum_{k=0}^{n-1}\varphi S^{k}\bold u
\right)
\tag3.2
 $$
for those complex numbers $z$ for which the series converge, see
[ParPo2].

We now define fine dynamical multifractal zeta-functions
from [MiOl2,Ol7].
However, we first introduce the following notation.
We denote the family of Borel probability measures on 
$\Sigma_{\smallG}^{\Bbb N}$
and the 
 family of shift invariant Borel probability measures on 
$\Sigma_{\smallG}^{\Bbb N}$
by $\Cal P(\Sigma_{\smallG}^{\Bbb N})$
and
$\Cal P_{S}(\Sigma_{\smallG}^{\Bbb N})$, respectively, i\.e\.
we write
 $$
 \align
 \Cal P(\Sigma_{\smallG}^{\Bbb N})
&=
 \Big\{
 \mu\,\Big|\,
 \text{
 $\mu$ is a Borel probability measures on 
$\Sigma_{\smallG}^{\Bbb N}$
 }
 \Big\}\,,
 \tag3.3\\
 \Cal P_{S}(\Sigma_{\smallG}^{\Bbb N})
&=
 \Big\{
 \mu\,\Big|\,
 \text{
 $\mu$ is a shift invariant Borel probability measures on 
$\Sigma_{\smallG}^{\Bbb N}$
 }
 \Big\}\,;
 \tag3.4\\
 \endalign
 $$
we will always
equip
$\Cal P(\Sigma_{\smallG}^{\Bbb N})$ 
and
$\Cal P_{S}(\Sigma_{\smallG}^{\Bbb N})$ with the weak topologies.
We now fix a metric space $X$ and a
continuous map
$U:\Cal P(\Sigma_{\smallG}^{\Bbb N})\to X$.
The multifractal zeta-function
framework
developed 
in [MiOl2,Ol7]
depends on the space $X$ and the map $U$;
judicious
choices  of $X$ and $U$ 
will provide
 important 
examples, including,
multifractal spectra of 
graph-directed self-conformal measures
(see [MiOl2,Ol7] and Section 5.1)
and
a variety of
multifractal spectra 
of
ergodic averages
of continuous
functions on graph-directed self-conformal sets
(see [MiOl2,Ol7] and Section 5.2).
Next,
 for a positive integer $n$,
 let
$L_{n}:\Sigma_{\smallG}^{\Bbb N}\to\Cal P(\Sigma_{\smallG}^{\Bbb N})$ be defined 
by
 $$
 L_{n}\bold i
 =
 \frac{1}{n}\sum_{k=0}^{n-1}\delta_{S^{k}\bold i}\,;
 \tag3.5
 $$
 recall, that $S:\Sigma_{\smallG}^{\Bbb N}\to\Sigma_{\smallG}^{\Bbb N}$ denotes the shift map.
 We can now define the 
multifractal pressure and zeta-function
associated with the space $X$ and the map $U$.

 \bigskip
 
 \proclaim{Definition.
 The multifractal pressure 
 $ \underline P_{C}^{U}(\varphi)$
 and
 $ \overline P_{C}^{U}(\varphi)$
associated with the space $X$ and the map $U$}
Let $X$ be a metric space and let $U:\Cal P(\Sigma_{\smallG}^{\Bbb N})\to X$ be 
continuous with respect to the weak topology.
Fix a continuous function
$\varphi:\Sigma_{\smallG}^{\Bbb N}\to\Bbb R$.
    For $C\subseteq X$,
we define the lower and upper
mutifractal pressure
of 
 $\varphi$
associated with the space $X$ and the map $U$ by
 $$
 \aligned
  \underline P_{C}^{U}(\varphi)
 &=
 \,
  \liminf_{n}
 \,
  \,\,
  \frac{1}{n}
 \,\,
 \log
   \sum
  \Sb
  \bold i\in\Sigma_{\smallsmallG}^{n}\\
  {}\\
  UL_{n}[\bold i]\subseteq C
  \endSb
   \,\,
  \sup_{\bold u\in[\bold i]}
  \,\,
 \exp
 \,\,
 \sum_{k=0}^{n-1}\varphi S^{k}\bold u\,,\\
   \overline P_{C}^{U}(\varphi)
 &=
  \limsup_{n}
  \,\,
  \frac{1}{n}
 \,\,
 \log
   \sum
  \Sb
 \bold i\in\Sigma_{\smallsmallG}^{n}\\
  {}\\
  UL_{n}[\bold i]\subseteq C
  \endSb
 \,\,
  \sup_{\bold u\in[\bold i]}
  \,\,
 \exp
 \,\,
 \sum_{k=0}^{n-1}\varphi S^{k}\bold u\,.
 \endaligned
 \tag3.6
 $$

\endproclaim

 \bigskip
 
 \proclaim{Definition.
 The dynamical multifractal zeta-function 
 $\zeta_{C}^{\dyn,U}(\varphi;\cdot)$
 associated with the space $X$ and the map $U$}
 Let $X$ be a metric space and let $U:\Cal P(\Sigma_{\smallG}^{\Bbb N})\to X$ be 
continuous with respect to the weak topology.
Fix a continuous function
$\varphi:\Sigma_{\smallG}^{\Bbb N}\to\Bbb R$.
    For $C\subseteq X$,
we define the 
dynamical
multifractal
 zeta-function $\zeta_{C}^{\dyn,U}(\varphi;\cdot)$
associated with the space $X$ and the map $U$ by
 $$
  \zeta_{C}^{\dyn,U}(\varphi;z)
  =
   \sum_{n}
   \,\,
   \frac{z^{n}}{n}
 \left(
   \sum
  \Sb
 \bold i\in\Sigma_{\smallsmallG}^{n}\\
  {}\\
  UL_{n}[\bold i]\subseteq C
  \endSb
  \,\,
  \sup_{\bold u\in[\bold i]}
  \,\,
 \exp
 \,\,
 \sum_{k=0}^{n-1}\varphi S^{k}\bold u
\right)
\tag3.7
 $$
for those complex numbers $z$ for which the series converges.

\endproclaim

\bigskip

\noindent
{\bf Remark.}
Comparing the
definition of the pressure (3.1)
(the dynamical zeta-function (3.2))
and
the
definition of the multifractal pressure (3.6)
(the dynamical multifractal zeta-function (3.7)),
it is clear that 
the 
multifractal pressure
(the dynamical multifractal zeta-function)
is obtained by only summing over those 
strings 
$\bold i\in\Sigma_{\smallG}^{n}$
that are multifractally relevant, i\.e\.
those 
strings 
$\bold i\in\Sigma_{\smallG}^{n}$
for which
$UL_{n}[\bold i]\subseteq C$.

\bigskip

\noindent
{\bf Remark.}
It is clear that if $C=X$, then the 
multifractal 
\lq\lq constraint"
$UL_{n}[\bold i]\subseteq C$
is vacuously
satisfied,
and that,
in this case,
the 
multifractal pressure and 
dynamical multifractal zeta-function
reduce to the usual pressure and the usual
dynamical zeta-function, i\.e\.
 $
 \underline P_{X}^{U}(\varphi)
 =
 \overline P_{X}^{U}(\varphi)
 =
  P(\varphi)
 $
and
 $
 \zeta_{X}^{\dyn,U}(\varphi;\cdot)
 =
 \zeta^{\dyn}(\varphi;\cdot)
 $.

\bigskip

Recall, the following classical and beautiful
result known.
If
$(K_{\smallvertexi})_{\smallvertexi\in\smallV}$
denotes the
the
graph-directed
 self-conformal sets
 (2.6)
 and the OSC is satisfied, then $\dim_{\Haus}K_{\smallvertexi}=s$ 
 for all $\vertexi$
 where
 $s$ is the unique real number
such that
 $P(s\Lambda)=0$
 (alternatively,
 $s$ 
  is the unique real number such that
  $\sigma_{\radius}\big(\,\zeta^{\dyn}(s\Lambda;\cdot)\,\big)
  =
  1)$;
here  
 $\Lambda$
is the scaling function defined by (2.10).
The formula $P(s\Lambda)=0$ is known as Bowen's formula.
The next proposition shows
that
 if $\varphi:\Sigma_{\smallG}^{\Bbb N}\to\Bbb R$ is a continuous function
 with $\varphi<0$, then
there is are unique real numbers
 $\,\,\scri f(C)$ and
 $\,\,\scri F\,(C)$
solving the 
natural multifractal 
versions of Bowen's equations.

\bigskip

\proclaim{Proposition B [MiOl2,Ol7]. 
Solutions to mutifractal Bowen equations}
Let $X$ be a metric space and let $U:\Cal P(\Sigma_{\smallG}^{\Bbb N})\to X$ be 
continuous with respect to the weak topology.
Let $C\subseteq X$ be a  subset of $X$.
Fix a continuous
function
$\varphi:\Sigma_{\smallG}^{\Bbb N}\to\Bbb R$
with
$\varphi<0$.
Then 
there are unique real numbers
$\,\,\scri f(C)$
and
$\scri F\,(C)$ 
such that
 $$
 \align
 \lim_{r\searrow 0}\overline P_{B(C,r)}^{U}(\,\,\,\scri f(C)\,\varphi)
 &=
 0\,,
 \tag3.8\\
 \overline P_{C}^{U}(\,\scri F\,(C)\,\varphi)
&=
0\,;
\tag3.9
 \endalign
 $$
alternatively,
$\,\,\scri f(C)$
and
$\scri F\,(C)$ 
 are the  unique real numbers 
 such that
 $$
 \align
 \lim_{r\searrow 0}
 \sigma_{\radius}
 \big(
 \,
 \zeta_{B(C,r)}^{\dyn,U}(\,\,\,\scri f(C)\,\varphi;\cdot)
 \,
 \big)
 &=
 1\,,\\
 \sigma_{\radius}
 \big(
 \,
 \zeta_{C}^{\dyn,U}(\,\scri F\,(C)\,\varphi;\cdot)
 \,
 \big)
&=
1\,.
 \endalign
 $$ 
 If $C=\{\alpha\}$, then we will
 write 
 $\,\,\scri f(\alpha)$
and
$\scri F\,(\alpha)$
for
$\,\,\scri f(C)$
and
$\scri F\,(C)$, respectively.
\endproclaim

\bigskip

%
The main results in [MiOl2,Ol7]
relate the solutions
 $\,\,\scri f(C)$ and
 $\,\,\scri F\,(C)$
of the multifractal Bowen
 equations (3.8) and (3.9)
to various multifractal spectra.
However, it is clear that if
$C=\{\alpha\}$, then
the
sum in (3.7) may be empty,
and the 
solution
 $\,\scri F\,(\alpha)$
 to the equation
 $\overline P_{C}^{U}(\,\scri F\,(\alpha)\,\varphi)
=
0$
is therefore
 $-\infty$, i\.e\.
$
\,\scri F\,(\alpha)=-\infty
$;
in fact, examples in [MiOl2,Ol7] shows that this may happen
for all but countably many $\alpha$.
Hence, 
if $C$ is 
\lq\lq too small", then it may happen that
the
pressure (3.6)
and  the zeta-function (3.7)
do not  encode
sufficient
information
allowing us to recover any  
meaningful
dynamical or geometric
 characteristics associated with $U$,
including, for example,
multifractal spectra.
Because of this problem, 
[MiOl2,Ol7] proceed in the 
following 
two equally naturally ways.
Either, we can 
consider
a family
of enlarged
\lq\lq target" sets 
shrinking to the original 
main \lq\lq target" $\{\alpha\}$;
this approach will be referred to as the
shrinking target approach.
Or, alternatively,
we can consider
a
fixed enlarge 
\lq\lq target" set
and regard this as our
original main \lq\lq target";
this approach will be referred to as the
fixed target approach.

\bigskip

{\bf The shrinking target setting.}
In the
shrinking target
setting,
[MiOl2,Ol7]
provide a
variational principle
for
the solution
$\,\,\scri f(C)$ to the multifractal Bowen equation (3.8);
thisn is Theorem C below.
Below we 
denote the entropy
of a measure $\mu\in\Cal P(\Sigma_{\smallG}^{\Bbb N})$
by $h(\mu)$.

\bigskip

\proclaim{Theorem C [MiOl2,Ol7].
The shrinking target
multifractal Bowen equation}
Let $X$ be a metric space and let $U:\Cal P(\Sigma_{\smallG}^{\Bbb N})\to X$ be 
continuous with respect to the weak topology.
Let $C\subseteq X$ be a  subset of $X$.
Fix a continuous function  $\varphi:\Sigma_{\smallG}^{\Bbb N}\to\Bbb R$ 
with
$\varphi<0$
and let
$\,\,\scri f(C)$
be the unique real number
such that
 $$
 \align
 \lim_{r\searrow 0}
 \sigma_{\radius}
 \big(
 \,
 \zeta_{B(C,r)}^{\dyn,U}(\,\,\scri f(C)\,\varphi;\cdot)
 \,
 \big)
 &=
 1\,.
 \endalign
 $$ 
Then
 $$
\scri f(C)
=
\sup
\Sb
\mu\in\Cal P_{S}(\Sigma_{\smallsmallG}^{\Bbb N})\\
{}\\
U\mu\in \overline C
\endSb
-
\frac{h(\mu)}{\int\varphi\,d\mu}\,.
$$ 
\endproclaim

\bigskip

{\bf The fixed target setting.}
If the set $C$ is 
\lq\lq too small",
then
it follows from the above discussion 
that the solution
 $\,\scri F\,(C)$
 to the equation
 $\overline P_{C}^{U}(\,\scri F\,(C)\,\varphi)
=
0$
does not encode any meaningful dynamical or geometric information about $U$.
However, if the set $C$ satisfies a 
non-degeneracy 
condition guaranteeing that
it is not
\lq\lq too small"
(namely condition (3.10) below),
then 
results can be obtained.
Indeed,
[MiOl2,Ol7] provides
a
variational principles for the 
 solution
$\,\,\scri F\,(C)$ to the multifractal Bowen equation (3.9)
 in the
fixed target
setting.;
this is Theorem D below.

\bigskip

\proclaim{Theorem D [MiOl2,Ol7].
The fixed target
multifractal Bowen equation}
Let $X$ be a normed vector space.
Let $\Gamma:\Cal P(\Sigma_{\smallG}^{\Bbb N})\to X$
be continuous and affine
and let
$\Delta:\Cal P(\Sigma_{\smallG}^{\Bbb N})\to \Bbb R$
be continuous and affine
with
$\Delta(\mu)\not=0$
for all $\mu\in\Cal P(\Sigma_{\smallG}^{\Bbb N})$.
Define 
$U:\Cal P(\Sigma_{\smallG}^{\Bbb N})\to X$
by
$U=\frac{\Gamma}{\Delta}$.
Let $C$ be a closed and convex subset of $X$ and assume that
 $$
 \overset{\,\circ}\to{C}
 \cap
 \,
 U\big(\,\Cal P_{S}(\Sigma_{\smallG}^{\Bbb N})\,\big)
 \not=
 \varnothing\,.
 \tag3.10
 $$
Let $\varphi:\Sigma_{\smallG}^{\Bbb N}\to\Bbb R$ be continuous
with
$\varphi<0$.
Let 
$\scri F\,(C)$
be the unique real number
such that
 $$
 \align
 \sigma_{\radius}
 \big(
 \,
 \zeta_{C}^{\dyn,U}(\,\,\scri F\,(C)\,\varphi;\cdot)
 \,
 \big)
 &=
 1\,.
 \endalign
 $$ 
Then
 $$
\scri F\,(C)
=
\sup
\Sb
\mu\in\Cal P_{S}(\Sigma_{\smallsmallG}^{\Bbb N})\\
{}\\
U\mu\in C
\endSb
-
\frac{h(\mu)}{\int\varphi\,d\mu}\,.
$$ 
\endproclaim 
  \bigskip

\newpage


\heading{4. Coarse multifractal zeta-functions.}\endheading

This is the main sections in the paper.
In this section we introduce the 
the coarse
multifractal zeta-functions
and show that
the fine dynamical multifractal
zeta-functions defined in Section 3 and the 
coarse
multifractal zeta-functions
defined below are linked by 
a
Multifractal Formalism.
Finally, in Section
5
we provide a number of examples, including
multifractal spectra of graph-directed self-conformal measures
and multifractal spectra of ergodic averages of
continuous functions on 
graph-directed self-conformal  sets.

Throughout this section, and in the remaining parts of the paper, we will use
the following notation.
If 
$(\lambda_{n})_{n}$
and $(a_{n})_{n}$ are sequences of real numbers
and if $f$ is the 
\lq\lq zeta"-function defined by
$f(s)=
\sum_{n}a_{n}e^{\lambda_{n}s}$
for $s\in\Bbb C$,
then we will denote the abscissa 
of convergence 
of $f$ by
$\sigma_{\abs}(f)$, i\.e\.
 $$
 \sigma_{\abs}(f)
 =
 \text{
 \lq\lq the abscissa of convergence of $f$"
 }\,.
 $$
We can give define the 
coarse multifractal zeta-function 
and the 
corresponding
Renyi dimension
 associated with the space $X$ and the map $U$.

\bigskip
 
 \proclaim{Definition.
 The coarse multifractal zeta-function 
 $\zeta_{C}^{\co,U}(\varphi;\cdot)$
 associated with the space $X$ and the map $U$}
 Let $X$ be an inner product space with inner product 
$\langle\cdot|\cdot\rangle$
 and let $U:\Cal P(\Sigma_{\smallG}^{\Bbb N})\to X$ be 
continuous with respect to the weak topology.
Fix a continuous function
$\varphi:\Sigma_{\smallG}^{\Bbb N}\to\Bbb R$.
For $q\in X$,
we define the coarse multifractal zeta-function 
 $\zeta_{q}^{\co,U}(\varphi;\cdot)$
associated with the space $X$ and the map $U$ 
by
 $$
 \align
 \zeta_{q}^{\co,U}(\varphi;s)
&=
 \sum_{\bold i}
 \,\,
 \exp
 \Bigg(
 \Bigg(
\sup_{\bold u\in[\bold i]}\big\langle q\big |UL_{|\bold i|}\bold u\, \big\rangle
+s
\Bigg)
 \sup_{\bold u\in[\bold i]}
 \sum_{k=0}^{|\bold i|-1}\varphi(S^{k}\bold u)
 \Bigg)\\
&{}\\ 
&=
\sum_{n}
\,\,
 \sum_{\bold i\in\Sigma_{\smallsmallG}^{n}}
 \,\,
 \exp
 \Bigg(
 \Bigg(
\sup_{\bold u\in[\bold i]}\big\langle q\big |UL_{n}\bold u\, \big\rangle
+s
\Bigg)
 \sup_{\bold u\in[\bold i]}
 \sum_{k=0}^{n-1}\varphi(S^{k}\bold u)
 \Bigg)
 \endalign 
$$
for those complex numbers $s$ for which the series converges.

We define the  zeta-function 
Renyi dimension
 $\tau(q)$
associated with the space $X$ and the map $U$
by
 $$
 \tau(q)
 =
 \sigma_{\abs}
 \big(
 \,
 \zeta_{q}^{\co,U}(\varphi;\cdot)
 \,
 \big)\,.
 \tag4.1
 $$

\endproclaim

  \bigskip

{\bf Properties of $\tau$.}
It turns out that the 
zeta-function Renyi dimension
satisfies 
a variational principle. This is the statement of Theorem 4.1 below.

\bigskip

\proclaim{Theorem 4.1. The variational principle for $\tau$}
Let $X$ be an inner product space
with
inner product
$\langle\cdot|\cdot\rangle$
and let $U:\Cal P(\Sigma_{\smallG}^{\Bbb N})\to X$ be 
continuous with respect to the weak topology.
Let $q\in X$.
Fix a H\"older continuous function $\varphi:\Sigma_{\smallG}^{\Bbb N}\to\Bbb R$ with $\varphi<0$,
and let the function $\tau$ be defined by (4.1).
\roster
\item"(1)"
We have
 $$
 \tau(q)
 =
\sup
\Sb
\mu\in\Cal P_{S}(\Sigma_{\smallsmallG}^{\Bbb N})
\endSb
\Bigg(
\,
-
\frac{h(\mu)}{\int\varphi\,d\mu}
-
\big\langle q\big|U\mu \big\rangle
\,
\Bigg)\,.
$$ 
\item"(2)"
We have
  $$
\sup
\Sb
\mu\in\Cal P_{S}(\Sigma_{\smallsmallG}^{\Bbb N})
\endSb
\Bigg(
\,
h(\mu)
+
\big(
\big\langle q\big|U\mu \big\rangle
+
\tau(q)
\big)
\int\varphi\,d\mu
\,
\Bigg)
=
0\,.
$$ 
\endroster
\endproclaim

\bigskip

\noindent
Theorem 4.1 is proved in Section 6 and Section 7
 using techniques from large 
deviation theory

\bigskip

\noindent
{\bf Remark.}
Recall that
if $\varphi:\Sigma_{\smallG}^{\Bbb N}\to\Bbb R$
is a continuous function, then $P(\varphi)$
denotes the pressure of $\varphi$ (see (3.1)).
Also, recall that it follows from the variational principle
that
 $$
 P(\varphi)
 =
 \sup
\Sb
\mu\in\Cal P_{S}(\Sigma_{\smallsmallG}^{\Bbb N})
\endSb
\Bigg(
\,
h(\mu)
+
\int\varphi\,d\mu
\,
\Bigg)\,.
 $$
We conclude from this and Theorem 4.1 that
$\tau(0)$
is the solution to the 
following
pressure equation, namely,
 $$
 P(\,\tau(0)\varphi\,)
 =
 \sup
\Sb
\mu\in\Cal P_{S}(\Sigma_{\smallsmallG}^{\Bbb N})
\endSb
\Bigg(
\,
h(\mu)
+
\tau(0)\int\varphi\,d\mu
\,
\Bigg)
=
0\,.
 $$
This viewpoint
suggests that
$\tau(q)$ may be viewed as the 
solution to the
\lq\lq non-linear" pressure
equation
given by
$$
\sup
\Sb
\mu\in\Cal P_{S}(\Sigma_{\smallsmallG}^{\Bbb N})
\endSb
\Bigg(
\,
h(\mu)
+
\big(
\big\langle q\big|U\mu \big\rangle
+
\tau(q)
\big)
\int\varphi\,d\mu
\,
\Bigg)
=
0\,;
\tag4.1
$$ 
observe that (4.1)
is a
\lq\lq non-linear" pressure
equation
since it
(in addition
to the usual \lq\lq linear" term
$\tau(q)
\int\varphi\,d\mu$
depending linearly on $\mu$)
contains
the \lq\lq non-linear"
term
$\big\langle q\big|U\mu \big\rangle
\int\varphi\,d\mu$
depending non-linearly on $\mu$.

\bigskip

The next result shows that
$\tau$
has the expected properties, including,
(1) convexity.
(2) differentiability (under certain additional conditions),
and
(3) $\tau(0)$ equals the Hausdorff dimension of the 
invariant set of the associated dynamical system.

\bigskip

\proclaim{Theorem 4.2. Properties of $\tau$}
Let $X$ be an inner product space and let $U:\Cal P(\Sigma_{\smallG}^{\Bbb N})\to X$ be 
continuous with respect to the weak topology.
Fix a continuous function  $\varphi:\Sigma_{\smallG}^{\Bbb N}\to\Bbb R$,
and let the function $\tau$ be defined by (4.1).
\roster
\item"(1)"
The function $\tau$ is convex.
\item"(2)"
If $\varphi$ is H\"older continuous with $\varphi<0$
and there is a H\"older continuous function
$\psi:\Sigma_{\smallG}^{\Bbb N}\to\Bbb R^{M}$
such that
 $$
 U\mu=
 \frac{\int\psi\,d\mu}{\int\varphi\,d\mu}
 $$
 for all $\mu\in\Cal P\big(\Sigma_{\smallG}^{\Bbb N}\big)$,
 then the function $\tau$ is real analytic.
 \item"(3)"
 Let
 $(K_{\smallvertexi})_{\smallvertexi\in\smallV}$ be the 
 graph-directed
 self-conformal sets defined by (2.6)
 and
let $\varphi=\Lambda$ where
 $\Lambda:\Sigma_{\smallG}^{\Bbb N}\to\Bbb R$
 is the scaling function defined by (2.10).
If the OSC is satisfied
and $\vertexi\in\V$, then
  $$
  \tau(0)
  =
  \dim_{\Haus}K_{\smallvertexi}
  $$ 
 for all
 inner product spaces
 $X$ and  
 all
 continuous maps $U:\Cal P(\Sigma_{\smallG}^{\Bbb N})\to X$. 

\endroster
\endproclaim

  \bigskip

  {\bf The multifractal formalism and $\tau$.}
  The next result shows that
the fine dynamical multifractal
zeta-functions defined in Section 3 and the 
coarse
multifractal zeta-functions are linked by 
a
Multifractal Formalism.
Recall, that
if
  $X$ is an inner product space with inner product $\langle\cdot|\cdot\rangle$
  and
  $f:X\to\Bbb R$ is a  function, then the
  Legendre transform
  $f^{*}:X\to[-\infty,\infty]$ of $f$ is defined by
   $$
   f^{*}(x)
   =
   \inf_{y}
   (
   \langle x|y\rangle+f(y)
   )\,.
   $$

  Our main results are divided into two
parts.
  In analogy with the classical
   Multifractal Formalism
  for measures
  (see Section 1)
the first part provides a natural
  zeta-function analogue of
  Proposition A, namely,
    the 
  \lq\lq zeta-function
  multifractal spectrum"
  $\,\,\scri f(C)$
  is
  {\it always}
  majorized by
  the Legendre transform of $\tau$;
  more precisely
  $\,\,\scri f(C)
  \le
  \sup_{\alpha\in \overline C}\tau^{*}(\alpha)$.
  This is the statement of Theorem 4.3.


   The second part
 shows 
 that
 if 
 $\tau$ is differentiable, then
  the 
  \lq\lq zeta-function
  multifractal spectrum"
  $\,\,\scri f(C)$
  equals
  the Legendre transform of $\tau$;
  more precisely
  $\,\,\scri f(C)
  =
  \sup_{\alpha\in \overline C}\tau^{*}(\alpha)$.
  In 
 other words,
  if
 $\tau$ is differentiable,
 then
 the
  Multifractal Formalism
  is satisfied.
  This is the statement of Theorem 4.4.
  We use the following
  notation is the statement of Theorem 4.4 below, namely,
 if
$f:\Bbb R^{M}\to\Bbb R$ is a convex function,
then
we write
$D_{\subdif}f(x)$ for the subdifferential of $f$ at $x$.

  \bigskip

\proclaim{Theorem 4.3.
The majorant theorem}
Let $X$ be an inner product space and let $U:\Cal P(\Sigma_{\smallG}^{\Bbb N})\to X$ be 
continuous with respect to the weak topology.
Fix a continuous function  $\varphi:\Sigma_{\smallG}^{\Bbb N}\to\Bbb R$ 
with
$\varphi<0$,
and let the function $\tau$ be defined by (4.1).
Let $C\subseteq X$ be a  subset of $X$
and let
$\,\,\scri f(C)$
be the unique real number
such that
 $$
 \align
 \limsup_{r\searrow 0}
 \sigma_{\radius}
 \big(
 \,
 \zeta_{B(C,r)}^{\dyn,U}(\,\,\scri f(C)\,\varphi;\cdot)
 \,
 \big)
 &=
 1\,.
 \endalign
 $$ 
Then
 $$
  \scri f(C)
  \le
  \sup_{\alpha\in \overline C}\tau^{*}(\alpha)\,.
  $$

\endproclaim

\bigskip

\proclaim{Theorem 4.4.
The Multifractal Formalism for multifractal zeta-functions}
Let $U:\Cal P(\Sigma_{\smallG}^{\Bbb N})\to \Bbb R^{M}$ be 
continuous with respect to the weak topology.
Fix a 
H\"older continuous function  $\varphi:\Sigma_{\smallG}^{\Bbb N}\to\Bbb R$ 
with
$\varphi<0$,
and let the function $\tau$ be defined by (4.1).
\roster
\item"(1)"
Let $C\subseteq \Bbb R^{M}$
and
let
$\,\,\scri f(C)$
be the unique real number
such that
 $$
 \align
 \limsup_{r\searrow 0}
 \sigma_{\radius}
 \big(
 \,
 \zeta_{B(C,r)}^{\dyn,U}(\,\,\scri f(C)\,\varphi;\cdot)
 \,
 \big)
 &=
 1\,.
 \endalign
 $$ 
 
 \noindent
 Assume that there is a subset $Q$ of $\Bbb R^{M}$ 
 such that
 $\tau$
 is differentiable
 at all
 $q\in Q$
 and
if
  $\alpha\in\overline C\cap U\big(\Cal P_{S}(\Sigma_{\smallG}^{\Bbb N})\big)$
and $\tau^{*}(\alpha)>-\infty$,
then
there is $q\in Q$
with
$\alpha\in-D_{\subdif} \tau (q)$.
 %
%
%
%
%
%
%
%
%
 Then
 $$
  \scri f(C)
  =
  \sup_{\alpha\in \overline C}\tau^{*}(\alpha)\,.
  $$
\item"(2)"
Let $\alpha\in \Bbb R^{M}$
and
let
$\,\,\scri f(\alpha)$
be the unique real number
such that
 $$
 \align
 \limsup_{r\searrow 0}
 \sigma_{\radius}
 \big(
 \,
 \zeta_{B(C,r)}^{\dyn,U}(\,\,\scri f(\alpha)\,\varphi;\cdot)
 \,
 \big)
 &=
 1\,.
 \endalign
 $$ 
 
 \noindent
 Assume that there is a point $q$ in $\Bbb R^{M}$ 
 such that
 $\tau$
 is differentiable
 at 
 $q$
 and
if
  $\alpha\in\overline C\cap U\big(\Cal P_{S}(\Sigma_{\smallG}^{\Bbb N})\big)$
and $\tau^{*}(\alpha)>-\infty$,
then
$\alpha\in-D_{\subdif} \tau (q)$.
%
%
%
%
%
%
%
%
%
%
%
 Then
 $$
  \scri f(\alpha)
  =
 \tau^{*}(\alpha)\,.
  $$

\endroster
(Observe that since $\tau$ is convex,
we conclude that $\tau$
is differentiable almost everywhere,
and the conclusion in Theorem 4.4
is therefore 
satisfied for 
\lq\lq many"
points $\alpha$.)  

\endproclaim

\bigskip

Theorem 4.3 is proved in Section 8
and Theorem 4.4 is proved in Section 9.
Theorem 4.4 shows that 
 if
 $\tau$ is differentiable,
 then
 the
  Multifractal Formalism
  for multifractal zeta-functions
  is satisfied.
  It is therefore of interest to 
  find conditions guaranteeing 
  the differentiability of $\tau$.
 Corollary 4.5 below provides 
 such a condition. 
In Corollary 4.5 we will use the following notation.
If
$\bold f=(f_{1},\ldots,f_{M}):\Sigma_{\smallG}^{\Bbb N}\to\Bbb R^{M}$
is a continuous function taking values in $\Bbb R^{M}$
and
$\mu\in\Cal P(\Sigma_{\smallG}^{\Bbb N})$, then we will
write
$$
\int \bold f\,d\mu
 =
 \Bigg(
 \int f_{1}\,d\mu,\ldots,\int f_{M}\,d\mu
 \Bigg)\,.
 $$
We can now state Corollary 4.5.

\bigskip

\proclaim{Corollary 4.5}
 Fix a H\"older continuous function  $\varphi:\Sigma_{\smallG}^{\Bbb N}\to\Bbb R$ 
with
$\varphi<0$,
and let the function $\tau$ be defined by (4.1).
Let 
$\psi:\Sigma_{\smallG}^{\Bbb N}\to\Bbb R^{M}$
be a H\"older continuous function
and
define
 $U:\Cal P(\Sigma_{\smallG}^{\Bbb N})\to \Bbb R^{M}$ by
  $$
  U\mu
  =
  \frac{\int\psi\,d\mu}{\int\varphi\,d\mu}\,.
  $$
Let $C\subseteq \Bbb R^{M}$ be a  subset of $\Bbb R^{M}$.
Let
$\,\,\scri f(C)$
be the unique real number
such that
 $$
 \align
 \limsup_{r\searrow 0}
 \sigma_{\radius}
 \big(
 \,
 \zeta_{B(C,r)}^{\dyn,U}(\,\,\scri f(C)\,\varphi;\cdot)
 \,
 \big)
 &=
 1\,.
 \endalign
 $$ 
Then
 $$
  \scri f(C)
  =
  \sup_{\alpha\in \overline C}\tau^{*}(\alpha)\,.
  $$

\endproclaim 

\noindent{\it Proof}\newline
\noindent
We conclude  from  Theorem 4.2  that $\tau $ is real analytic,
and 
the desired result therefore follows from Theorem 4.4.
\hfill$\square$
  
  \bigskip

  \newpage


\centerline{\smc 5. Applications: }

\centerline{\smc multifractal spectra of measures}

\centerline{\smc  and}

\centerline{\smc  
multifractal spectra of ergodic Birkhoff averages}

 \medskip

We will now consider several of applications of 
Theorems 4.1--4.3
to multifractal spectra of measures and ergodic averages.
In particular, we consider the following examples:

\medskip

\roster

\item"$\bullet$"
 Section 5.1: Multifractal spectra of graph-directed self-conformal measures.
 
\medskip

\item"$\bullet$"
Section 5.2: Multifractal spectra of ergodic Birkhoff averages
of continuous functions on 
graph-directed self-conformal sets.

\endroster

\medskip

{\bf 5.1.
Multifractal 
spectra of 
graph directed 
self-conformal measures.}
Let
$(\,
 \V,\,\allowmathbreak
 \E,\,\allowmathbreak
 (V_{\smallvertexi})_{\smallvertexi\in\smallV},\,\allowmathbreak
 (X_{\smallvertexi})_{\smallvertexi\in\smallV},\,\allowmathbreak
 (S_{\smalledge})_{\smalledge\in\smallE},\,\allowmathbreak
 (p_{\smalledge})_{\smalledge\in\smallE}\,
 )$
 be a
 graph-directed conformal iterated function system with probabilities
 (see Section 2)
 and let
   $(K_{\smallvertexi})_{\smallvertexi\in\smallV}$
   and 
   $(\mu_{\smallvertexi})_{\smallvertexi\in\smallV}$
   be the list of graph-directed self-conformal 
   sets 
   and the list of graph-direted self-conformal
   measures
   associated with the list 
 $(\,
 \V,\,\allowmathbreak
 \E,\,\allowmathbreak
 (V_{\smallvertexi})_{\smallvertexi\in\smallV},\,\allowmathbreak
 (X_{\smallvertexi})_{\smallvertexi\in\smallV},\,\allowmathbreak
 (S_{\smalledge})_{\smalledge\in\smallE},\,\allowmathbreak
 (p_{\smalledge})_{\smalledge\in\smallE}\,
 )$, respectively,
   i\.e\.
   the sets in the list
  $(K_{\smallvertexi})_{\smallvertexi\in\smallV}$
  are the unique non-empty
  compact sets
  satisfying (2.6)
   and 
   the measures in the list
   $(\mu_{\smallvertexi})_{\smallvertexi\in\smallV}$
   are the unique probability measures
   satisfying (2.7).
Recall that
the Hausdorff multifractal spectrum
$f_{\mu_{\smallsmallvertexi}}$ of $\mu_{\smallvertexi}$ 
is defined
by
 $$
 \align
 f_{\mu_{\smallsmallvertexi}}(\alpha)
&=
   \,\dim_{\Haus}
   \left\{x\in K_{\smallvertexi}
    \,\left|\,
     \lim_{r\searrow0}
     \frac
     {\log\mu_{\smallvertexi} (B(x,r))}{\log r}
     =
     \alpha
      \right.
   \right\}\,,
  \endalign
 $$
for $\alpha\in\Bbb R$, see Section 1.
If the OSC is satisfied,
then
the 
multifractal spectrum
$f_{\mu_{\smallsmallvertexi}}(\alpha)$
can be computed as follows.
Define
$\Phi:\Sigma_{\smallG}^{\Bbb N}\to\Bbb R$ by
 $$
 \Phi(\bold i)
 =
 \log p_{\ini(\bold i)}
 \tag5.1
 $$
for 
$\bold i=\edge_{1}\edge_{2}\ldots
\in\Sigma_{\smallG}^{\Bbb N}$
and recall that the map 
$\Lambda:\Sigma_{\smallG}^{\Bbb N}\to\Bbb R$
is defined in (2.10).
Next,
we define the function $\beta:\Bbb R\to\Bbb R$
by
 $$
 P\big(
 \,
 q\Phi
 +
 \beta(q)\Lambda
 \,
 \big)
 =
 0\,.
 \tag5.2
 $$
  If the OSC is satisfied, then it follows
  from [EdMa,Col1,Col2,Pa] that
   $$
   f_{\mu_{\smallsmallvertexi}}(\alpha)
   =
   \beta^{*}(\alpha)\,.
   \tag5.3
   $$

Of course, in general, the limit
 $\lim_{r\searrow0}
     \frac
     {\log\mu_{\smallsmallvertexi} (B(x,r))}{\log r}$
     may not exist.
     Indeed, recently 
     Barreira \& Schmeling [BaSc]
     (see also
     Olsen \& Winter [OlWi1,OlWi2],
     Xiao, Wu \& Gao [XiWuGa]
     and 
Moran [Mo])
have shown that 
in many cases
the set $D_{\smallvertexi}$ of divergence points of $\mu_{\smallvertexi}$, 
i\.e\. the set 
of points $x$ for which the 
limit
$\lim_{r\searrow0}
     \frac
     {\log \mu_{\smallsmallvertexi} (B(x,r))}{\log r}$
does not exist, 
is  highly
\lq\lq visible", namely it has full Hausdorff dimension,
i\.e\.
$\dim_{\Haus}D_{\smallvertexi}
=
\dim_{\Haus}K_{\smallvertexi}$.
This
 suggests that the set 
 of divergence points
 has a surprising rich 
 fractal structure.
 In order to explore this more
 carefully,
  Olsen \& Winter [OlWi1,OlWi2]
 introduced various
 generalised multifractal spectra functions designed to 
 \lq\lq see"
 different sets of divergence points.
To define these spectra 
 we introduce the following notation.
 If $M$ 
 is a  metric space
 and
 $f:(0,\infty)\to M$ is a function, then we write
 $\acc_{r\searrow 0}f(r)$
 for the set of accumulation 
 points of $f$ as $r\searrow 0$, i\.e\.
  $$
  \underset {r\searrow0}\to\acc\,\,f(r)
  =
  \Big\{x\in M
  \,\Big|\,
  \text{$x$ is an accumulation point of $f$ as $r\searrow 0$}
  \Big\}\,.
  $$
Olsen \& Winter [OlWi1]
introduced 
and investigated
 the generalised
 Hausdorff multifractal spectrum
 $F_{\mu_{\smallsmallvertexi}}$ of  $\mu_{\smallsmallvertexi}$
defined
by
 $$
 \align
F_{\mu_{\smallsmallvertexi}}(C)
&=
   \,\dim_{\Haus}
   \left\{x\in K_{\smallvertexi}
    \,\left|\,
     \,\,
     \underset {r\searrow0}\to\acc
 \,
     \frac
     {\log\mu_{\smallvertexi} (B(x,r))}{\log r}
  \subseteq
 C
      \right.
      \,
   \right\}
  \endalign
 $$
 for 
 $C\subseteq \Bbb R$.
There 
is a natural
divergence point analogue of (5.3).
Namely if the OSC is satisfied, then
 $$
 F_{\mu_{\smallsmallvertexi}}(C)
 =
 \sup_{\alpha\in C}\beta^{*}(\alpha)\,.
 $$
for all $C\subseteq \Bbb R$, see
[Mo,OlWi1,LiWuXi]
(see also [Ca,Vo]
for earlier
but related results in 
a
slightly different setting).

\newpage

As a first application of 
Corollary 4.5
we obtain a 
fine
dynamical multifractal zeta-function
with an associated
 Bowen equation whose solution
 $\scr f\,\,(C)$
equals the 
generalised
multifractal spectrum 
of the graph-directed 
self-conformal measure $\mu_{\smallvertexi}$
and a  family of 
coarse 
multifractal zeta-functions
with 
abscissae  of 
convergence 
whose Legendre transform
equals
$\scr f\,\,(C)$.
This is the content of Theorem 5.1 below.
Theorem 5.1
follows by applying 
Corollary 4.5
to $X=\Bbb R$
and
the function
$U:\Cal P(\Sigma_{\smallG}^{\Bbb N})\to\Bbb R$
defined by
$U\mu
  =
 \frac{\int\Phi\,d\mu}{\int\Lambda\,d\mu}$
 where $\Phi$  and $\Lambda$ are defined in (5.1) and (2.10), respectively.
 For this choice of $X$ and $U$, 
the
zeta-functions
$\zeta_{C}^{\dyn,U}(\varphi;z)$
and
$\zeta_{C}^{\co,U}(\varphi;z)$ 
are clearly given by
  $$
  \zeta_{C}^{\dyn,U}(\varphi;z)
  =
  \sum_{n}
  \,\,
  \frac{z^{n}}{n}
  \,\,
  \left(
  \sum
  \Sb
   \bold i\in\Sigma_{\smallsmallG}^{n}\\
   {}\\
   \forall
 \bold v\in\Sigma_{\smallsmallG}^{\Bbb N}
 \,\,
 \text{with}
 \,\,
 \termi(\bold i)=\ini(\bold v)
 \,\,:\,\,
 \frac{\log p_{\bold i}}{\log |DS_{\bold i}(\pi\bold v)|}
\in
C
  \endSb
  \sup_{\bold u\in[\bold i]}
  \,\,
  \exp
  \sum_{k=0}^{n-1}
  \varphi S^{k}\bold u
  \right)
  \tag5.5
  $$
and
  $$
  \zeta_{q}^{\co,U}(\varphi;s)
  =
  \sum_{n}
  \,
  \sum_{\bold i\in\Sigma_{\smallsmallG}^{n} }
  \,
  \exp
  \left(
  \,
  \left(
  \sup_{
   \bold v\in\Sigma_{\smallsmallG}^{\Bbb N}
 \,\,
 \text{with}
 \,\,
 \termi(\bold i)=\ini(\bold v)
  }
  q
  \,
 \frac{\log p_{\bold i}}{\log |DS_{\bold i}(\pi\bold v)|}
   +
  s
  \right)
  \sup_{\bold u\in[\bold i]}
  \sum_{k=0}^{n-1}
  \varphi S^{k}\bold u
  \right)\,.
  \tag5.6
  $$
It follows from the 
Principle of Bounded Variation that
$\log |DS_{\bold i}(\pi\bold v)|$
behaves asymptotically as
$\diam K_{\bold i}$
as $n$ tends to infinity for all
$\bold i\in\Sigma_{\smallG}^{n}$
and
 $\bold v\in\Sigma_{\smallsmallG}^{\Bbb N}$
 with
$\termi(\bold i)=\ini(\bold v)$;
this, in turn, implies that
the radius of convergence and the abscissa of convergence
of the 
zeta-functions in (5.5) and (5.6)
behave asymptotically
as the 
the radius of convergence and the abscissa of convergence
of the zeta-functions
  $$
  \zeta_{C}^{\dyn\text{-}\scon}(\varphi;z)
  =
  \sum_{n}
  \,\,
  \frac{z^{n}}{n}
  \,\,
  \left(
  \sum
  \Sb
   \frac{\log p_{\bold i}}{\log\diam_{\smallNormal} K_{\bold i}|}
\in
C
  \endSb
  \sup_{\bold u\in[\bold i]}
  \,\,
  \exp
  \sum_{k=0}^{n-1}
  \varphi S^{k}\bold u
  \right)
  \tag5.7
  $$
and
  $$
  \zeta_{q}^{\co\text{-}\scon}(\varphi;z)
  =
  \sum_{n}
  \,
  \sum_{\bold i\in\Sigma_{\smallsmallG}^{n} }
  \,
  \exp
  \left(
  \,
  \left(
  q
  \,
 \frac{\log p_{\bold i}}{\log \diam_{\Normal} K_{\bold i}}
   +
  s
  \right)
  \sup_{\bold u\in[\bold i]}
  \sum_{k=0}^{n-1}
  \varphi S^{k}\bold u
  \right)\,,
  \tag5.8
  $$
obtained from (5.5) and (5.6)
by replacing
$|DS_{\bold i}(\pi\bold v)|$
by
the normalised diameter
$\diam_{\Normal} K_{\bold i}$
of $K_{\bold i}$ defined by
 $$
 \diam_{\Normal} K_{\bold i}
 =
 \frac{\diam K_{\bold i}}{\diam K_{\smallvertexi}}
 $$
for each finite string $\bold i\in\Sigma_{\smallG}^{*}$
with initial vertex equal to $\vertexi$;
Proposition 5.3 provides a  precise 
statement of this.
For this reason we have decided to formulate 
Theorem 5.1
using
the more natural
 zeta-functions
(5.7) and (5.8) (instead of (5.5) and (5.6)).

\bigskip

\proclaim{Theorem 5.1. 
Fine dynamical multifractal zeta-functions
and coarse multifractal zeta-functions
 for 
multifractal spectra of graph-directed
self-conformal measures}
Let $(\mu_{\smallvertexi})_{\smallvertexi\in\smallV}$
be the list of graph-directed self-conformal measures
associated with the list
$\big(\,
 \V,\,
 \E,\,
 (V_{\smallvertexi})_{\smallvertexi\in\smallV},\,\allowmathbreak
 (X_{\smallvertexi})_{\smallvertexi\in\smallV},\,\allowmathbreak
 (S_{\smalledge})_{\smalledge\in\smallE},\,\allowmathbreak
 (p_{\smalledge})_{\smalledge\in\smallE}\,
 \big)$, i\.e\.
$\mu_{\smallvertexi}$ is the unique probability measure such that
$ \mu_{\smallvertexi}
 =
 \sum_{\smalledge\in E_{\smallsmallvertexi}}\,
 p_{\smalledge}\,\mu_{\termi(\smalledge)}\circ S_{\smalledge}^{-1}$.

For $C\subseteq\Bbb R$ and 
a continuous function 
$\varphi:\Sigma_{\smallG}^{\Bbb N}\to\Bbb R$,   we
define the fine dynamical
graph-directed
self-conformal multifractal zeta-function 
$\zeta_{C}^{\sdyncon}(\varphi;z)$ by
(5.7).

 For $q\in\Bbb R$ and 
a continuous function 
$\varphi:\Sigma_{\smallG}^{\Bbb N}\to\Bbb R$,   we
define the coarse
graph-directed
self-conformal multifractal zeta-function 
$\zeta_{q}^{\co\text{-}\scon}(\varphi;s)$ by
(5.8).

Let $\Lambda$ be defined by (2.10) and 
let $\beta$ be defined by (5.2),
 and write
  $$
  \tau(q)
  =
  \sigma_{\abs}
  \big(
  \,
   \zeta_{q}^{\co\text{-}\scon}(\Lambda;s)
  \,
  \big)\,.
  \tag5.9
  $$

\noindent
 {\rm (1)} For all $q$, we have
  $$
  \tau(q)
  =
  \beta(q)\,.
  $$

\bigskip

\noindent
 {\rm (2)} Assume 
 that $C\subseteq\Bbb R$ is closed.

 \roster

  \item"(2.1)"
 There is a unique real number 
 $\,\,\scr f\,\,(C)$ such that
  $$
  \lim_{r\searrow0}
  \sigma_{\radius}
  \big(
  \,
  \zeta_{B(C,r)}^{\sdyncon}(\,\,\scr f\,\,(C)\,\Lambda;\cdot)
  \,
  \big)
  =
  1\,.
  $$
It $\alpha\in\Bbb R$ and 
$C=\{\alpha\}$, then we will write
  $\,\,\scr f\,\,(\alpha)
  =
 \,\,\scr f\,\,(C)$.

\item"(2.2)"
We have
 $$
 \scr f\,\,(C)
 =
  \sup_{\alpha\in C}\tau^{*}(\alpha)\,.
 $$
  
   \item"(2.3)"
 If the OSC is satisfied, then we have
 $$
\align
\quad
  \scr f\,\,(C)
 &=
 F_{\mu_{\smallsmallvertexi}}(C)
 =
 \dim_{\Haus}
 \Bigg\{
 x\in K_{\smallvertexi}
 \,\Bigg|\,
 \,\underset{r\searrow 0}\to\acc
 \frac{\log\mu_{\smallvertexi}(B(x,r))}{\log r}
 \subseteq
 C
 \Bigg\}
 =
  \sup_{\alpha\in C}\tau^{*}(\alpha)\,.\\
 \endalign
 $$
 In particular, 
 if the OSC is satisfied and $\alpha\in\Bbb R$,
 then we have
 $$
\align
\quad
  \scr f\,\,(\alpha)
 &=
 f_{\mu_{\smallsmallvertexi}}(\alpha)
 =
 \dim_{\Haus}
 \Bigg\{
 x\in K_{\smallvertexi}
 \,\Bigg|\,
 \lim_{r\searrow 0}
  \frac{\log\mu_{\smallvertexi}(B(x,r))}{\log r}
  =
 \alpha
 \Bigg\}
 =
 \tau^{*}(\alpha)
 \,.\qquad
 \endalign
 $$

 \endroster
 
 \bigskip
 
 \noindent
 {\rm (3)} Assume that $C\subseteq\Bbb R$ is a closed 
interval
 with
 $\overset{\,\circ}\to{C}
  \cap
  \big(
  -\beta'(\Bbb R)
  \big)
  \not=
  \varnothing$.

\roster

 \item"(3.1)"
 There is a unique real number 
 $\scr F\,\,(C)$ such that
  $$
  \sigma_{\radius}
  \big(
  \,
  \zeta_{C}^{\sdyncon}(\scr F\,\,(C)\,\Lambda;\cdot)
  \,
  \big)
  =
  1\,.
  $$

\item"(3.2)"
We have
 $$
 \scr F\,\,(C)
 =
  \sup_{\alpha\in C}\tau^{*}(\alpha)\,.
 $$

 \item"(3.3)"
 If the OSC is satisfied then
 $$
\align
\quad
  \scr F\,\,(C)
 &=
 F_{\mu_{\smallsmallvertexi}}(C)
 =
 \dim_{\Haus}
 \Bigg\{
 x\in K_{\smallvertexi}
 \,\Bigg|\,
 \,\underset{r\searrow 0}\to\acc
 \frac{\log\mu_{\smallvertexi}(B(x,r))}{\log r}
 \subseteq
 C
 \Bigg\}
 =
  \sup_{\alpha\in C}\tau^{*}(\alpha)\,.\\
 \endalign
 $$
\endroster
 
\endproclaim

\bigskip

We will prove Theorem 5.1
below.
However, we first note that
if 
all the maps 
$S_{\smalledge}$ are similarities,
then
the coarse multifractal zeta-functions 
in (5.6) and (5.8)
can be computed explicitly;
this is the content of Theorem 5.2 below.
In order to state this result, we introduce the following notation.
Assuming that the maps 
$S_{\smalledge}$ are similarities, i\.e\.
if for 
each 
$\edge\in \E$
there is a number $r_{\smalledge}\in(0,1)$ such that
 $$
 |S_{\smalledge}(x)-S_{\smalledge}(y)|
 =
 r_{\smalledge}|x-y|
 $$
 for all
 $x,y\in X_{\termi(\smalledge)}$,
 then
we define the matrix $A(q,s)$
for $q\in\Bbb R$ and $s\in\Bbb C$ by
$$
A(q,s)=\big(\,a_{\smallvertexi,\smallvertexj}(q,s)\,\big)_{\smallvertexi,\smallvertexj\in\smallV}
\tag5.10
 $$
 where
 $$
 a_{\smallvertexi,\smallvertexj}(q,s)
 =
 \sum_{\smalledge\in\smallE_{\smallsmallvertexi,\smallsmallvertexj}}
 p_{\smalledge}^{q}r_{\smalledge}^{s}\,.
 $$
We can now state Theorem 5.2.

 \bigskip

 \proclaim{Theorem 5.2}
 Let $(\mu_{\smallvertexi})_{\smallvertexi\in\smallV}$
be the list of graph-directed self-conformal measures
associated with the list
$\big(\,
 \V,\,
 \E,\,
 (V_{\smallvertexi})_{\smallvertexi\in\smallV},\,
 (X_{\smallvertexi})_{\smallvertexi\in\smallV},\,
 (S_{\smalledge})_{\smalledge\in\smallE},\,
 (p_{\smalledge})_{\smalledge\in\smallE}\,
 \big)$, i\.e\.
$\mu_{\smallvertexi}$ is the unique probability measure such that
$ \mu_{\smallvertexi}
 =
 \sum_{\smalledge\in E_{\smallsmallvertexi}}\,
 p_{\smalledge}\,\mu_{\termi(\smalledge)}\circ S_{\smalledge}^{-1}$.
Assume that all the maps 
$S_{\smalledge}$ are similarities
and 
let the matrix $A(q,s)$
be defined my (5.10).
Let $\Lambda$ be defined by (2.10) and 
let the zeta-functions
$\zeta_{q}^{\co,U}(\Lambda;s)$
and
$\zeta_{q}^{\co\text{-}\scon}(\Lambda;s)$ be defined by (5.6) and (5.8), respectively.
Finally, let $\tau(q)$ be defined by (5.9).

\newpage

If $s\in\Bbb C$ with
 $\text{\rm Re}\,s>\tau(q)$,
 then $I-A(q,s)$ is invertible and
 we have
  $$
 \zeta_{q}^{\co,U}(\Lambda;s)
=
 \zeta_{q}^{\co\text{-}\scon}(\Lambda;s)
 =
  \bold 1^{\top}\,(I-A(q,s))^{-1}\,A(q,s)\,\bold 1
    $$
where
$\bold 1=(1)_{\smallvertexi\in\smallV}$
is the column vector 
in $\Bbb R^{\smallV}$
consisting of $1$'s
and 
$\bold 1^{\top}$ denotes the transpose of $\bold 1$.

In particular, if the graph $\G$ has only one vertex and 
$N$ edges labelled $1,2,\ldots,N$, then
we have
 $$
 \zeta_{q}^{\co,U}(\Lambda;s)
=
 \zeta_{q}^{\co\text{-}\scon}(\Lambda;s)
 =
  \frac{\sum_{i=1}^{N}p_{i}^{q}r_{i}^{s}}{1-\sum_{i=1}^{N}p_{i}^{q}r_{i}^{s}}\,.
  \tag5.11
    $$
 
 \endproclaim
 \noindent{\it Proof}\newline
 \noindent
 We fist note that, in this particular case,
 the zeta-functions
 in (5.6) and (5.8) coincide, i\.e\.
 $ \zeta_{q}^{\co,U}(\Lambda;s)
=
 \zeta_{q}^{\co\text{-}\scon}(\Lambda;s)$.

 Next, fix 
 $s\in\Bbb C$ with
 $\text{\rm Re}\,s>\tau(q)$
 and write
 $A=A(q,s)$ for brevity.
 Since  
 $\sum_{k=0}^{n-1}
  \Lambda S^{k}\bold u
  =
  \log r_{\bold i}$
  for all $\bold i\in\Sigma_{\smallG}^{*}$
  and
  all
 $\bold u\in[\bold i]$, we conclude that
  $$
  \align
  \zeta_{q}^{\co,U}(\Lambda;s)
 &=
 \sum_{n}
  \,
  \sum_{\bold i\in\Sigma_{\smallsmallG}^{n} }
  \,
  \exp
  \left(
  \,
  \left(
  \sup_{
   \bold v\in\Sigma_{\smallsmallG}^{\Bbb N}
 \,\,
 \text{with}
 \,\,
 \termi(\bold i)=\ini(\bold v)
  }
  q
  \,
 \frac{\log p_{\bold i}}{\log |DS_{\bold i}(\pi\bold v)|}
   +
  s
  \right)
  \sup_{\bold u\in[\bold i]}
  \sum_{k=0}^{n-1}
  \Lambda S^{k}\bold u
  \right)\\
   &=
 \sum_{n}
  \,
  \sum_{\bold i\in\Sigma_{\smallsmallG}^{n} }
  \,
  \exp
  \left(
  \,
  \left(
  \sup_{
   \bold v\in\Sigma_{\smallsmallG}^{\Bbb N}
 \,\,
 \text{with}
 \,\,
 \termi(\bold i)=\ini(\bold v)
  }
  q
  \,
 \frac{\log p_{\bold i}}{\log r_{\bold i}}
   +
  s
  \right)
 \log r_{\bold i}
  \right)\\
    &=
 \sum_{n}
  \,
  \sum_{\bold i\in\Sigma_{\smallsmallG}^{n} }
  \,
p_{\bold i}^{q}r_{\bold i}^{s}\,.
  \endalign
  $$
 Noticing that
  $  \sum_{\bold i\in\Sigma_{\smallsmallG}^{n} }
  \,
p_{\bold i}^{q}r_{\bold i}^{s}
=
\bold 1^{\top}\,A^{n}\,\bold 1$, it therefore follows that
 $$
 \zeta_{q}^{\co,U}(\Lambda;s)
 =
 \sum_{n}
  \,
  \bold 1^{\top}\,A^{n}\,\bold 1\,.
  \tag5.12
  $$

 Let $M_{\smallV}(\Bbb C)$ 
 denote the vector space of $\V\times\V$ square matrices
 with complex entries
 and let $\|\cdot\|$
 be a norm on 
 $M_{\smallV}(\Bbb C)$
 (since $M_{\smallV}(\Bbb C)$
 is finite dimensional, all norms
 on $M_{\smallV}(\Bbb C)$
 are equivalent and it is therefore not important which norm we use).
 We now claim that
 (1) the matrix $I-A$ is invertible,
 and (2)
 the series
 $\sum_{n}A^{n}$ converges with respect to the norm $\|\cdot\|$
 and
  $$
 \sum_{n}A^{n}
 =
(I-A)^{-1}A\,.
\tag5.13
$$ 
We will now prove (5.13).  
It follows from Theorem 5.1
that
$
 \sigma_{\abs}
 \big(
 \,
 \zeta_{q}^{\co,U}(\Lambda;\cdot)
 \,
 \big)
=
 \sigma_{\abs}
 \big(
 \,
 \zeta_{q}^{\co\text{-}\scon}(\Lambda;\cdot)
 \,
 \big)
=
\beta(q) $.
It also follows from
[EdMa] that
$\beta(q)$ is the unique real number such that
$\rho_{\specradsmall}\,A(q,\beta(q))
  =
  1$
  and 
  the function $t\to\rho_{\specradsmall}\,A(q,t)$
  is strictly decreasing;
 here and below we use the following notation, namely,
 if $M$ is a square matrix, then we will write
 $\rho_{\specradsmall} \,M$ for the spectral radius of $M$.
 Consequently, since
 $\text{Re}\,s>
 \sigma_{\abs}
 \big(
 \,
 \zeta_{q}^{\co,U}(\Lambda;\cdot)
 \,
 \big)
 =
 \beta(q)$,
 we conclude that
 $\rho_{\specradsmall} A=\rho_{\specradsmall}\,A(q,s)<1$.
 It also follows from the spectral formula that
 $\lim_{n}\|A^{n}\|^{\frac{1}{n}}
 =
\rho_{\specradsmall} A$, and we therefore deduce that
$\lim_{n}\|A^{n}\|^{\frac{1}{n}}
 =
\rho_{\specradsmall} A<1$.
This implies that we can find a positive constant $c_{0}$ with
$0<c_{0}<1$ and an integer $N_{0}$ such that
$\|A^{n}\|^{\frac{1}{n}}\le c_{0}$ for $n\ge N_{0}$, i\.e\.
$\|A^{n}\|\le c_{0}^{n}$ for $n\ge N_{0}$, 
whence
$\sum_{n\ge 0}\|A^{n}\|
\le
\sum_{n=0}^{N_{0}}\|A^{n}\|
 +
\sum_{n>N_{0}}c_{0}^{n}
<
\infty$.
We conclude from this and the completeness of $M_{\smallV}(\Bbb C)$
that the series $\sum_{n\ge 0}A^{n}$ converges with respect to $\|\cdot\|$.
Writing
$B=\sum_{n\ge 0}A^{n}$, we clearly have
$(I-A)B=I$ and $B(I-A)=I$. This shows that $I-A$ is invertible with
$(I-A)^{-1}=B=\sum_{n\ge 0}A^{n}$.
It is not difficult to see that this equality  implies (5.13).
This
 completes the proof of (5.13).

Next, we note that it 
follows easily from (5.13) that
 $$
 \align
  \sum_{n}
  \,
  \bold 1^{\top}\,A^{n}\,\bold 1
 &=
  \bold 1^{\top}\Bigg( \sum_{n}A^{n}\Bigg)\bold 1
  =
  \bold 1^{\top}\,(I-A)^{-1} A \,\bold 1\,.
  \tag5.14
 \endalign
 $$

  The desired result now follows immediately from (5.12) and (5.14).
 \hfill$\square$

 \newpage


The zeta-function in (5.11) 
and its connection to the 
Renyi dimensions $\tau_{\mu}(q)$
of self-similar measures $\mu$
has recently been investigated
in [Le-VeMe,Ol5,Ol6].

We  will now prove Theorem 5.1.
We  note that Part (2) and Part (3) follow from Part (1) and 
[MiOl2,Ol7];
hence it suffices to prove Part (1).
Recall  that the functions
$\Phi:$ 
and
 $\Lambda$
are defined in (5.1) and (2.10), respectively.
The proof of Theorem 5.1 now
follows by applying Corollary 4.5 to the function
$U:\Cal P(\Sigma_{\smallG}^{\Bbb N})\to\Bbb R$
defined by
 $$
 U\mu
  =
 \frac{\int\Phi\,d\mu}{\int\Lambda\,d\mu}\,,
 \tag5.15
 $$
noticing that in this particular  case
the zeta-function 
$\zeta_{q}^{\co,U}(\varphi;s)$
is given by (5.6).
We first prove
the following auxiliary
result
showing that the 
abscissa of convergence
of the 
zeta-functions
$\zeta_{q}^{\co,U}(\varphi;s)$
and
$ \zeta_{q}^{\co\text{-}\scon}(\varphi;s)$
are equal.

\bigskip

\proclaim{Proposition 5.3}
Let $U$
be defined by (5.15).
Fix a continuous function
$\varphi:\Sigma_{\smallG}^{\Bbb N}\to\Bbb R$
with $\varphi<0$.
\roster
\item"(1)"
There is a sequence $(\Delta_{n})_{n}$ with $\Delta_{n}>0$ 
and 
$\Delta_{n}\to 0$ such that
for
all
$n\in\Bbb N$, $\bold i\in\Sigma_{\smallG}^{n}$ and 
$\bold u\in\Sigma_{\smallG}^{\Bbb N}$
with
 $\termi(\bold i)=\ini(\bold u)$, we have
 $
 |
 \,
\frac{\log p_{\bold i}}{\log |DS_{\bold i}(\pi\bold u)|}
 -
  \frac{\log p_{\bold i}}{\log \diam_{\smallNormal} K_{\bold i}}
 \,
|
\le
\Delta_{n}$.

\item"(2)"
For all $q\in X$,
we have
 $
 \sigma_{\abs}
 \big(
 \,
 \zeta_{q}^{\co,U}(\varphi;\cdot)
 \,
 \big)
=
 \sigma_{\abs}
 \big(
 \,
 \zeta_{q}^{\co\text{-}\scon}(\varphi;\cdot)
 \,
 \big)
 $.

\endroster
\endproclaim
\noindent{\it Proof}\newline
\noindent
(1)
It is well-known and follows from the 
Principle of Bounded Distortion 
(see, for example, [Bar,Fa])
that there is a constant $c>0$ such that
for all integers $n$
and all $\bold i\in\Sigma_{\smallG}^{n}$
and all $\bold u,\bold v\in[\bold i]$,
 we have
$\frac{1}{c}
\le
\frac{|DS_{\bold i}(\pi S^{n}\bold u)|}{\diam_{\smallNormal} K_{\bold i}}
\le
c$
and
$\frac{1}{c}
\le
\frac{|DS_{\bold i}(\pi S^{n}\bold u)|}{|DS_{\bold i}(\pi S^{n}\bold v)|}
\le
c$.
It is not difficult to see that the desired result follows from this.

\noindent
(2)
It is not difficult to see that this statement follows
from Part 1 and we have therefore decided to omit the details.
\hfill$\square$

\bigskip

We can now prove Theorem 5.1.

\bigskip

\noindent{\it Proof of Theorem 5.1}\newline
\noindent
(1)
Let the map $U$ be defined by (5.15).
Proposition 5.3 and Theorem 4.1 clearly imply that
 $
 0
 =
 \sup_{\mu\in\Cal P_{S}(\Sigma_{\smallsmallG}^{\Bbb N})}
 (
 h(\mu)
 +
 (qU\mu+\tau(q))
 \int\Lambda\,d\mu
 )
  =
 \sup_{\mu\in\Cal P_{S}(\Sigma_{\smallsmallG}^{\Bbb N})}
 (
 h(\mu)
 +
 \int(q\Phi+\tau(q)\Lambda)\,d\mu
 $.
It follows from this and the variational principle that
$P(q\Phi+\tau(q)\Lambda)=0$, and (5.2) therefore
implies that
$\tau(q)=\beta(q)$.

  \noindent
  (2)--(3)
  These statements follow from (1) and 
  [MiOl2,Ol7].
  \hfill$\square$

\bigskip

{\bf 5.2. Multifractal spectra of 
ergodic Birkhoff averages.}
We first fix $\gamma\in(0,1)$ and define  the metric
 $\distance_{\gamma}$ on $\Sigma_{\smallG}^{\Bbb N}$
 as follows.
For $\bold i,\bold j\in\Sigma_{\smallG}^{\Bbb N}$ 
with $\bold i\not=\bold j$,
we will write
$\bold i\wedge\bold j$ for the longest common 
prefix of 
$\bold i$ and $\bold j$.
The metric
$\distance_{\gamma}$ is now defined by
$
 \distance_{\gamma}(\bold i,\bold j)
 =
 \gamma^{|\bold i\wedge\bold j|}$
 if $\bold i\not=\bold j$,
 and
 $ \distance_{\gamma}(\bold i,\bold j)
 =
 0$
if $\bold i=\bold j$.
for $\bold i,\bold j\in\Sigma_{\smallG}^{\Bbb N}$;
throughout this section, we equip
 $\Sigma_{\smallG}^{\Bbb N}$
with the metric  $\distance_{\gamma}$,
and continuity and Lipschitz properties of functions $f:\Sigma_{\smallG}^{\Bbb N}\to\Bbb R$
from $\Sigma_{\smallG}^{\Bbb N}$ to $\Bbb R$
will  always
refer to
the metric  $\distance_{\gamma}$.
Multifractal analysis of Birkhoff
averages has received significant interest
during
the past 10 years, see, for example,
[BaMe,FaFe,FaFeWu,FeLaWu,Oli,Ol2,OlWi1,OlWi2].
The 
multifractal
spectrum
$F_{f}^{\erg}$
of ergodic Birkhoff averages of a 
continuous function
$f:\Sigma_{\smallG}^{\Bbb N}\to\Bbb R$ is defined by
 $$
F_{f}^{\erg}(\alpha)
=
 \dim_{\Haus}
 \pi
 \Bigg\{
 \bold i\in\Sigma_{\smallG}^{\Bbb N}
 \,\Bigg|\,
 \lim_{n}
 \frac{1}{n}\sum_{k=0}^{n-1} f(S^{k}\bold i)
 =
 \alpha
 \Bigg\}
$$
for $\alpha\in\Bbb R$;
recall, that the map $\pi$ is defined in Section 2.
One of the main problems
in
multifractal analysis of Birkhoff
averages
is the detailed study of the multifractal
spectrum
$F_{f}^{\erg}$.
For example,
it  is
proved
(in different settings and at various levels of generality)
in [FaFe,FaFeWu,FeLaWu,Oli,Ol2,OlWi1]
that if $f:\Sigma_{\smallG}^{\Bbb N}\to\Bbb R$ is continuous
and $\alpha$ is a real number, then
$$
F_{f}^{\erg}(\alpha)
=
  \dim_{\Haus}
 \pi
 \Bigg\{
 \bold i\in\Sigma_{\smallG}^{\Bbb N} 
 \,\Bigg|\,
 \lim_{n}
 \frac{1}{n}\sum_{k=0}^{n-1}f(S^{k}\bold i)
 =
 \alpha
 \Bigg\}
 =
 \sup
 \Sb
 \mu\in\Cal P_{S}(\Sigma_{\smallsmallG}^{\Bbb N})\\
 \int  f\,d\mu
 =
 \alpha
 \endSb
 -
 \frac{h(\mu)}{\int\Lambda\,d\mu}\,.
 $$ 
where
$\Lambda$ be defined by (2.10).

As a second
 application of 
Theorems 4.1--4.4
we will now obtain
a fine dynamical
multifractal
zeta-function
with an associated Bowen equation
 whose solution
  $\,\,\scr f\,\,(C)$
equals the multifractal spectrum of ergodic Birkhoff averages
and a  family of coarse multifractal zeta-functions with 
abscissae of convergence whose Legendre transform equals 
  $\,\,\scr f\,\,(C)$.
We first 
state and prove a 
rather
general 
result, namely Theorem 5.4 below.
Using Theorem 5.4 we then 
deduce
analogous results
for
a number
of
different
multifractal 
spectra
of ergodic Birkhoff averages, including
$F_{f}^{\erg}\alpha)$; see
Theorem 5.5 and Theorem 5.6.
Recall that if
$\bold f:\Sigma_{\smallG}^{\Bbb N}\to\Bbb R^{M}$
is a continuous function with
$\bold f=(f_{1},\ldots, f_{M})$, then we will write
$
 \int\bold f\,d\mu
 =
 (
 \int f_{1}\,d\mu
 \,,\,
 \ldots
 \,,\,
 \int f_{M}\,d\mu
 )
 $
for $\mu\in\Cal P(\Sigma_{\smallG}^{\Bbb N})$.
Also, 
if $(x_{n})_{n}$ is a sequence of points in a metric space $M$, then we write
$\acc_{n} x_{n}$ for the set of accumulation points 
of the sequence  $(x_{n})_{n}$,
 i\.e\.
 $$
 \underset{n}\to{\acc} \,x_{n}
=
\Big\{
x\in M
\,\Big|\,
\text{
$x$ is an accumulation point of $(x_{n})_{n}$
}
\Big\}\,.
$$

\bigskip

\proclaim{Theorem 5.4.
Multifractal zeta-functions for
abstract 
multifractal spectra of ergodic Birkhoff averages}
Fix $\gamma\in(0,1)$
and
$W\subseteq \Bbb R^{I}$
and
let
$\pmb\Phi:\Sigma_{\smallG}^{\Bbb N}\to\Bbb R^{I}$ be a Lipschitz function
with respect to the metric $\distance_{\gamma}$
such that
 $\{
 \int
 \pmb\Phi\,d\mu
 \,|\,
 \mu
 \in\Cal P(\Sigma_{\smallG}^{\Bbb N})
 \}
 \subseteq
 W$.
Let
 $Q:W\to \Bbb R^{M}$
be a continuous function.

For $C\subseteq \Bbb R^{M}$ and 
a continuous function 
$\varphi:\Sigma_{\smallG}^{\Bbb N}\to\Bbb R$,   we
define the 
abstract fine dynamical
ergodic multifractal zeta-function associated with $Q$
by
  $$
  \zeta_{C}^{\dyn\text{-}\erg}(\varphi;z)
  =
  \sum_{n}
  \frac{z^{n}}{n}
  \left(
  \sum
  \Sb
 \bold i\in\Sigma_{\smallsmallG}^{n}\\
  {}\\
  \forall\bold v\in[\bold i]
  \,\,:\,\,
  Q\big(
 \frac{1}{n}\sum_{k=0}^{n-1}\pmb\Phi(S^{k}\bold v)
 \big)
 \in  C
  \endSb
  \sup_{\bold u\in[\bold i]}
  \exp
  \sum_{k=0}^{n-1}
  \varphi S^{k}\bold u
  \right)\,.
  $$

For $\bold q\in \Bbb R^{M}$ and 
a continuous function 
$\varphi:\Sigma_{\smallG}^{\Bbb N}\to\Bbb R$,   we
define the 
abstract coarse
ergodic multifractal zeta-function associated with $Q$
by
   $$
  \zeta_{\bold q}^{\co\text{-}\erg}(\varphi;s)
  =
  \sum_{n}
  \,
  \sum_{\bold i\in\Sigma_{\smallsmallG}^{n} }
  \,
  \exp
  \Bigg(
  \,
  \Bigg(
  \sup_{
  \bold v\in[\bold i]
  }
\Bigg\langle
\bold q
\Bigg|
Q
\Bigg(
 \frac{1}{n}\sum_{k=0}^{n-1}\pmb\Phi(S^{k}\bold v
\Bigg)
\Bigg\rangle
   +
  s
  \Bigg)
  \sup_{\bold u\in[\bold i]}
  \sum_{k=0}^{n-1}
  \varphi S^{k}\bold u
  \Bigg)\,.
  $$

Let 
 $\Lambda:\Sigma_{\smallG}^{\Bbb N}\to\Bbb R$ be defined by (2.10)
 and write
  $$
  \tau(\bold q)
  =
  \sigma_{\abs}
  \big(
  \,
  \zeta_{\bold q}^{\co\text{-}\erg}(\Lambda;\cdot)
  \,
  \big)\,.
  $$

\noindent
 {\rm (1)} For all $\bold q$, we have
  $$
  \tau(\bold q)
  =
  \sup_{\mu\in\Cal P_{S}(\Sigma_{\smallsmallG}^{\Bbb N})}
  \Bigg(
  -\frac{h(\mu)}{\int\Lambda\,d\mu}
  -
 \Bigg\langle
\bold q
\Bigg|
Q
\Bigg(
\int\pmb\Phi\,d\mu
\Bigg)
\Bigg\rangle 
\,\,
  \Bigg)\,.
 $$

\noindent
 {\rm (2)} Assume 
 that $C\subseteq \Bbb R^{M}$ is closed.

 \roster
  \item"(2.1)"
 There is a unique real number 
 $\,\,\scr f\,\,(C)$ such that
  $$
  \lim_{r\searrow0}
  \sigma_{\radius}
  \big(
  \,
  \zeta_{B(C,r)}^{\dyn\text{-}\erg}(\,\,\scr f\,\,(C)\,\Lambda;\cdot)
  \,
  \big)
  =
  1\,.
  $$

\item"(2.2)"
If for each $\pmb\alpha\in C$, there is $\bold q\in\Bbb R^{M}$ such that
$\pmb\alpha=-\nabla\tau(\bold q)$,
then
we have
$$
\,\,\scr f\,\,(C)
 =
 \sup_{\pmb\alpha\in C}
 \,
 \tau^{*}(\pmb\alpha)\,.
 $$

\item"(2.3)"
If the OSC is satisfied
and
if
for each $\pmb\alpha\in C$, there is $\bold q\in\Bbb R^{M}$ such that
$\pmb\alpha=-\nabla\tau(\bold q)$,
then we have
$$
\align
\quad
\,\,\scr f\,\,(C)
 &=
 \dim_{\Haus}
 \pi
 \Bigg\{
 \bold i\in\Sigma_{\smallG}^{\Bbb N}
 \,\Bigg|\,
 \,\underset{n}\to\acc
 \,\,
 Q
 \Bigg(
 \frac{1}{n}\sum_{k=0}^{n-1}\pmb\Phi(S^{k}\bold i)
 \Bigg)
\subseteq
C
 \Bigg\}
 =
  \sup_{\pmb\alpha\in C}
 \,
 \tau^{*}(\pmb\alpha)\,.\\
 \endalign
 $$
\endroster

\endproclaim
\noindent{\it Proof}\newline
\noindent
(1) This follows by applying Theorem 4.1 to the map
$U:\Cal P_{S}(\Sigma_{\smallG}^{\Bbb N})\to\Bbb R^{M}$
defined by
$U\mu=Q(\int\pmb\Phi\,d\mu)$.

\noindent
(2) This follows from Part (1), Theorem 4.4 and 
[MiOl2,Ol7].
\hfill$\square$

\bigskip

\newpage

\noindent
Below we present two
corollaries of Theorem 5.4.
Theorem 5.5 studies 
the multifractal spectrum 
$F_{f}^{\erg}(\alpha)$
of  continuos functions $f$.
The Legendre transform
representation of 
the spectrum $F_{f}^{\erg}(\alpha)$
in Part (2.3) of Theorem 5.4
seems to be new.

\bigskip

\proclaim{Theorem 5.5.
Multifractal zeta-functinons for
multifractal spectra of vector valued ergodic Birkhoff averages}
Fix $\gamma\in(0,1)$
and
let $\bold f:\Sigma_{\smallG}^{\Bbb N}\to\Bbb R^{M}$ be a Lipschitz function
with respect to the metric $\distance_{\gamma}$.
For $C\subseteq\Bbb R^{M}$ and 
a continuous function 
$\varphi:\Sigma_{\smallG}^{\Bbb N}\to\Bbb R$,   we
define the dynamical
ergodic multifractal zeta-function by
  $$
  \zeta_{C}^{\dyn\text{-}\vector\text{-}\erg}(\varphi;z)
  =
  \sum_{n}
  \frac{z^{n}}{n}
  \left(
  \sum
  \Sb
 \bold i\in\Sigma_{\smallsmallG}^{n}\\
  {}\\
  \forall\bold v\in[\bold i]
  \,\,:\,\,
 \frac{1}{n}\sum_{k=0}^{n-1}\bold f(S^{k}\bold v)
 \in  C
  \endSb
  \sup_{\bold u\in[\bold i]}
  \exp
  \sum_{k=0}^{n-1}
  \varphi S^{k}\bold u
  \right)\,.
  $$

For $\bold q\in \Bbb R^{M}$ and 
a continuous function 
$\varphi:\Sigma_{\smallG}^{\Bbb N}\to\Bbb R$,   we
define the 
coarse
ergodic multifractal zeta-function 
by
   $$
  \zeta_{\bold q}^{\co\text{-}\vector\text{-}\erg}(\varphi;s)
  =
  \sum_{n}
  \,
  \sum_{\bold i\in\Sigma_{\smallsmallG}^{n} }
  \,
  \exp
  \Bigg(
  \,
  \Bigg(
  \sup_{
  \bold v\in[\bold i]
  }
\Bigg\langle
\bold q
\Bigg|
 \frac{1}{n}\sum_{k=0}^{n-1}\bold f(S^{k}\bold v)
\Bigg\rangle
   +
  s
  \Bigg)
  \sup_{\bold u\in[\bold i]}
  \sum_{k=0}^{n-1}
  \varphi S^{k}\bold u
  \Bigg)\,.
  $$

Let 
 $\Lambda:\Sigma_{\smallG}^{\Bbb N}\to\Bbb R$ be defined by (2.10)
 and write
  $$
  \tau(\bold q)
  =
  \sigma_{\abs}
  \big(
  \,
  \zeta_{\bold q}^{\co\text{-}\vector\text{-}\erg}(\Lambda;s)
  \,
  \big)\,.
  $$

\noindent
 {\rm (1)} For all $\bold q$, we have
  $$
  \tau(\bold q)
  =
  \sup_{\mu\in\Cal P_{S}(\Sigma_{\smallsmallG}^{\Bbb N})}
  \Bigg(
  -\frac{h(\mu)}{\int\Lambda\,d\mu}
  -
 \Bigg\langle
\bold q
\Bigg|
\int
\bold f\,d\mu
\Bigg\rangle 
\,\,
  \Bigg)\,.
 $$

\noindent
 {\rm (2)} Assume 
 that $C\subseteq\Bbb R^{M}$ is closed.

 \roster
  \item"(2.1)"
 There is a unique real number 
 $\,\,\scr f\,\,(C)$ such that
  $$
  \lim_{r\searrow0}
  \sigma_{\radius}
  \big(
  \,
  \zeta_{B(C,r)}^{\dyn\text{-}\vector\text{-}\erg}(\,\,\scr f\,\,(C)\,\Lambda;\cdot)
  \,
  \big)
  =
  1\,.
  $$
If $\pmb\alpha\in\Bbb R^{M}$ and 
$C=\{\pmb\alpha\}$, then we will write
  $\,\,\scr f\,\,(\pmb\alpha)
  =
 \,\,\scr f\,\,(C)$.

\item"(2.2)"
If for each $\pmb\alpha\in C$, there is $\bold q\in\Bbb R^{M}$ such that
$\pmb\alpha=-\nabla\tau(\bold q)$,
then
we have
$$
\,\,\scr f\,\,(C)
 =
 \sup_{\pmb\alpha\in C}
 \,
 \tau^{*}(\pmb\alpha)\,.
  $$

\item"(2.3)"
If the OSC is satisfied
and
if for each $\pmb\alpha\in C$, there is $\bold q\in\Bbb R^{M}$ such that
$\pmb\alpha=-\nabla\tau(\bold q)$, then we have
$$
\align
\qquad
\,\,\scr f\,\,(C)
 &=
 \dim_{\Haus}
 \pi
 \Bigg\{
 \bold i\in\Sigma_{\smallG}^{\Bbb N}
 \,\Bigg|\,
 \,\underset{n}\to\acc
 \frac{1}{n}\sum_{k=0}^{n-1}\bold f(S^{k}\bold i)
\subseteq
C
 \Bigg\}
  =
 \sup_{\pmb\alpha\in C}
 \,
 \tau^{*}(\pmb\alpha)\,.\\
 \endalign
 $$
In particular, if the OSC is satisfied and $\pmb\alpha\in\Bbb R^{M}$
and there is $\bold q\in\Bbb R^{M}$ such that
$\pmb\alpha=-\nabla\tau(\bold q)$,
then we have 
$$
\align
\,\,\scr f\,\,(\pmb\alpha)
 &=
 \dim_{\Haus}
 \pi
 \Bigg\{
 \bold i\in\Sigma_{\smallG}^{\Bbb N}
 \,\Bigg|\,
 \lim_{n}
 \frac{1}{n}\sum_{k=0}^{n-1}\bold f(S^{k}\bold i)
 =
 \pmb\alpha
 \Bigg\}
  =
 \tau^{*}(\pmb\alpha)\,.
 \endalign
 $$
(Observe that since $\tau$ is convex,
we conclude that $\tau$
is differentiable almost everywhere,
and the conclusion in Part (2.3)
is therefore 
satisfied for 
\lq\lq many"
points $\pmb\alpha$.)   
\endroster

\endproclaim

\noindent{\it Proof}\newline
\noindent
This follows immediately
by applying
Theorem 5.4
to $W=\Bbb R^{M}$
and the maps
$\pmb\Phi:\Sigma_{\smallG}^{\Bbb N}\to\Bbb R^{M}$
and
 $Q:W\to\Bbb R^{M}$ 
defined by
$\pmb\Phi=\bold f$
and
$Q(\bold x)=\bold x$.
\hfill$\square$

 \bigskip

\newpage

\noindent
As a final application of Theorem 5.4 we
consider
a type of relative ergodic 
multifractal spectra involving
quantities 
similar to those 
appearing in H\"older's inequality;
for this reason we have decided
to refer to these multifractal spectra as
\lq\lq H\"older-like relative ergodic Birkhoff averages".

\bigskip

\proclaim{Theorem 5.6.
Multifractal zeta-functinons for
multifractal spectra of 
H\"older-like relative ergodic Birkhoff averages}
Fix $\gamma\in(0,1)$
and let $f_{1},\ldots,f_{M},g_{1},\ldots,g_{M}:\Sigma_{\smallG}^{\Bbb N}\to\Bbb R$
be Lipschitz functions
with respect to the metric
$\distance_{\gamma}$
and
assume that 
$f_{i}(\bold i)>0$
for all $i$ and all $\bold i$,
and that
$g_{i}(\bold i)>0$
for all $i$ and all $\bold i$.
Fix
 $s_{1},\ldots,s_{M},t_{1},\ldots,t_{M}>0$.
For $C\subseteq\Bbb R$ and 
a continuous function 
$\varphi:\Sigma_{\smallG}^{\Bbb N}\to\Bbb R$,   we
define the dynamical
H\"older-like
relative ergodic multifractal zeta-function by
  $$
  \zeta_{C}^{\dyn\text{-}\Hol\text{-}\erg}(\varphi;z)
  =
  \sum_{n}
  \frac{z^{n}}{n}
  \left(
  \sum
  \Sb
 \bold i\in\Sigma_{\smallsmallG}^{n}\\
  {}\\
  \forall\bold v\in[\bold i]
  \,\,:\,\,
 \frac
 {\prod_{i=1}^{M}\big(\frac{1}{n}\sum_{k=0}^{n-1}f_{i}(S^{k}\bold v)\big)^{s_{i}}}
 {\prod_{i=1}^{M}\big(\frac{1}{n}\sum_{k=0}^{n-1}g_{i}(S^{k}\bold v)\big)^{t_{i}}}
 \in  C
  \endSb
  \sup_{\bold u\in[\bold i]}
  \exp
  \sum_{k=0}^{n-1}
  \varphi S^{k}\bold u
  \right)\,.
  $$

For $q\in \Bbb R$ and 
a continuous function 
$\varphi:\Sigma_{\smallG}^{\Bbb N}\to\Bbb R$,   we
define the 
coarse
H\"older-like
ergodic multifractal zeta-function 
by
   $$
  \zeta_{q}^{\co\text{-}\Hol\text{-}\erg}(\varphi;s)
  =
  \sum_{n}
  \,
  \sum_{\bold i\in\Sigma_{\smallsmallG}^{n} }
  \,
  \exp
  \Bigg(
  \,
  \Bigg(
  \sup_{
  \bold v\in [\bold i]
  }
q
\,
\frac
{\prod_{i=1}^{M}\big(\frac{1}{n}\sum_{k=0}^{n-1}f_{i}(S^{k}\bold v)\big)^{s_{i}}}
 {\prod_{i=1}^{M}\big(\frac{1}{n}\sum_{k=0}^{n-1}g_{i}(S^{k}\bold v)\big)^{t_{i}}}
   +
  s
  \Bigg)
  \sup_{\bold u\in[\bold i]}
  \sum_{k=0}^{n-1}
  \varphi S^{k}\bold u
  \Bigg)\,.
  $$

Let 
 $\Lambda:\Sigma_{\smallG}^{\Bbb N}\to\Bbb R$ be defined by (2.10)
 and write
  $$
  \tau(q)
  =
  \sigma_{\abs}
  \big(
  \,
  \zeta_{q}^{\co\text{-}\Hol\text{-}\erg}(\Lambda;s)
  \,
  \big)\,.
  $$

\noindent
 {\rm (1)} For all $q$, we have
  $$
  \tau(q)
  =
  \sup_{\mu\in\Cal P_{S}(\Sigma_{\smallsmallG}^{\Bbb N})}
  \Bigg(
  -\frac{h(\mu)}{\int\Lambda\,d\mu}
  -
q
\,
\frac
{\prod_{i=1}^{M}\big(\int f_{i}\,d\mu\big)^{s_{i}}}
 {\prod_{i=1}^{M}\big(\int g_{i}\,d\mu\big)^{t_{i}}}
\,\,
  \Bigg)\,.
 $$

\noindent
 {\rm (2)} Assume 
 that $C\subseteq\Bbb R^{M}$ is closed.

 \roster
  \item"(2.1)"
 There is a unique real number 
 $\,\,\scr f\,\,(C)$ such that
  $$
  \lim_{r\searrow0}
  \sigma_{\radius}
  \big(
  \,
  \zeta_{B(C,r)}^{\dyn\text{-}\Hol\text{-}\erg}(\,\,\scr f\,\,(C)\,\Lambda;\cdot)
  \,
  \big)
  =
  1\,.
  $$
If $\alpha\in\Bbb R$ and 
$C=\{\alpha\}$, then we will write
  $\,\,\scr f\,\,(\alpha)
  =
 \,\,\scr f\,\,(C)$.

\item"(2.2)"
If for each $\alpha\in C$, there is $q\in\Bbb R$ such that
$\alpha=-\tau'(q)$,
then
we have
$$
\,\,\scr f\,\,(C)
 =
 \sup_{\alpha\in C}
 \,
 \tau^{*}(\alpha)\,.
 $$

 \item"(2.3)"
If the OSC is satisfied
and
if for each $\alpha\in C$, there is $q\in\Bbb R$ such that
$\alpha=-\tau'(q)$, then we have
$$
\align
\qquad\quad\,\,\,\,
\,\,\scr f\,\,(C)
 &=
 \dim_{\Haus}
 \pi
 \Bigg\{
 \bold i\in\Sigma_{\smallG}^{\Bbb N}
 \,\Bigg|\,
 \,\underset{n}\to\acc
 \,\,
 \frac
 {\prod_{i=1}^{M}\big(\frac{1}{n}\sum_{k=0}^{n-1}f_{i}(S^{k}\bold i)\big)^{s_{i}}}
 {\prod_{i=1}^{M}\big(\frac{1}{n}\sum_{k=0}^{n-1}g_{i}(S^{k}\bold i)\big)^{t_{i}}}
\subseteq
C
 \Bigg\}
 =
  \sup_{\alpha\in C}
 \,
 \tau^{*}(\alpha)\,.\\
 \endalign
 $$
In particular, if the OSC is satisfied and $\alpha\in\Bbb R$
and
there is $q\in\Bbb R$ such that
$\alpha=-\tau'(q)$, 
, then we have 
$$
\align
\quad
\,\,\scr f\,\,(\alpha)
 &=
 \dim_{\Haus}
 \pi
 \Bigg\{
 \bold i\in\Sigma_{\smallG}^{\Bbb N}
 \,\Bigg|\,
 \lim_{n}
 \frac
 {\prod_{i=1}^{M}\big(\frac{1}{n}\sum_{k=0}^{n-1}f_{i}(S^{k}\bold i)\big)^{s_{i}}}
 {\prod_{i=1}^{M}\big(\frac{1}{n}\sum_{k=0}^{n-1}g_{i}(S^{k}\bold i)\big)^{t_{i}}}
 =
 \alpha
 \Bigg\}
  =
 \tau^{*}(\alpha)\,.
 \endalign
 $$
(Observe that since $\tau$ is convex,
we conclude that $\tau$
is differentiable almost everywhere,
and the conclusion in Part (2.3)
is therefore 
satisfied for 
\lq\lq many"
points $\alpha$.)   
\endroster

\endproclaim

\noindent{\it Proof}\newline
\noindent
This follows immediately
by applying
Theorem 5.4
to $W=\Bbb R^{M}\times (\Bbb R\setminus\{0\})^{M}$
and the maps
$\pmb\Phi:\Sigma_{\smallG}^{\Bbb N}\to\Bbb R^{2M}$
and
 $Q:W\to\Bbb R$ 
defined by
$\pmb\Phi=(f_{1},\ldots,f_{M},g_{1},\ldots,g_{M})$
and
$Q(x_{1},\ldots,x_{M},y_{1},\ldots,y_{M})
=
\frac
{\prod_{i=1}^{M}x_{i}^{s_{i}}}
{\prod_{i=1}^{M}y_{i}^{t_{i}}}$.
\hfill$\square$


\newpage

\heading
{
6. Proof of Theorem 4.1 and Theorem 4.2, Part 1: The map $M_{n}$.}
\endheading

{\bf The map $M_{n}$.}
Since the graph
$\G=(\V,\E)$ is strongly connected,
it follows that
for
each $\bold i\in\Sigma_{\smallG}^{*}$, 
 we 
can choose
$\widehat{\bold i}\in\Sigma_{\smallG}^{*}$ 
with
$|\,\widehat{\bold i}\,|\le|\V|$
such that
$\termi(\bold i)=\ini(\,\widehat{\bold i}\,)$
and
$\termi(\,\widehat{\bold i}\,)=\ini(\bold i)$,
and we now define
$\overline{\bold i}\in\Sigma_{\smallG}^{\Bbb N}$ by
 $$
 \overline{\bold i}
 =
 \bold i\,\,\widehat{\bold i}\,\,
 \bold i\,\,\widehat{\bold i}\,\,
 \bold i\,\,\widehat{\bold i}\ldots\,.
 $$
 For  a positive integer $n$,
 we define
$M_{n}:\Sigma_{\smallG}^{\Bbb N}\to\Cal P_{S}(\Sigma_{\smallG}^{\Bbb N})$ by
 $$
 \align
 M_{n}\bold i
&=
 L_{n+|\,\widehat{\,\bold i|n\,}\,|}\left(\,\overline{\bold i|n}\,\right)\\
&= 
 \frac{1}{n+|\,\widehat{\,\bold i|n\,}\,|}
 \sum_{k=0}^{n+|\,\widehat{\,\bold i|n\,}\,|-1}
 \delta_{S^{k}(\,\overline{\bold i|n}\,)}
 \tag6.1
 \endalign
 $$
for $\bold i\in\Sigma_{\smallG}^{\Bbb N}$;
recall, that
the map
$L_{n}:\Sigma_{\smallG}^{\Bbb N}\to\Cal P(\Sigma_{\smallG}^{\Bbb N})$ is defined in 
(3.5).

\bigskip

{\bf Why the map $M_{n}$?}
The main reason for introduction the map $M_{n}$ is the following.
In order to prove Theorem 4.1 we will
use results from large deviation theory.
In particular, we
will use
Varadhan's
integral lemma
which says that if $X$
is a complete separable metric space
and $(P_{n})_{n}$ is a sequence of probability measures on
$X$ satisfying the large deviation property
with
rate constants $a_{n}\in(0,\infty)$ for $n\in\Bbb N$
and rate function $I:\Bbb R\to[-\infty,\infty]$
(this terminology will be explained in Section 7),
then
$$
\lim_{n}
\frac{1}{a_{n}}
\log
\int
\exp(a_{n}F)\,dP_{n}
=
-
\inf_{x\in X}(\,I(x)-F(x))
$$
for any bounded continuous function $F:X\to \Bbb R$
(see Section 7 for more a detailed and precise discussion of this result).

More precisely, in Section 7
we intend to use 
Varadhan's
integral lemma
to analyse the asymptotic
behaviour of the integral
$$
\frac{1}{n}
\log
\int
\exp(nF(L_{n}(\overline{\bold i|n})))\,d\Pi(\bold i)
\tag6.2
$$
as $n\to\infty$
where $\Pi$ is the Parry measure on $\Sigma_{\smallG}^{\Bbb N}$
(the Parry measure will be defined in Section 7)
 and the function
 $F:\Cal P\big(\Sigma_{\smallG}^{\Bbb N}\big)\to\Bbb R$
is defined by
 $
 F(\mu)
 =
 (\langle q|U\mu\rangle+s)\int\varphi\,d\mu$.
Defining
$\Lambda_{n}: \Cal P\big(\Sigma_{\smallG}^{\Bbb N}\big)\to\Bbb R$
by
$\Lambda_{n}(\bold i)
 =
 L_{n}(\,\overline{\bold i|n}\,)$,
 then (6.2) can be written as
$$
\frac{1}{n}
\log
\int
\exp(nF)\,d(\Pi\circ\Lambda_{n}^{-1})\,.
\tag6.3
$$
Consequently,
if the sequence $(\Pi\circ\Lambda_{n}^{-1})_{n}$
satisfied the large deviation property
with rate constants $a_{n}=n$,
then
Varadhan's
integral lemma
could be applied to analyse the asymptotic behaviour of the sequence of integrals in (6.3). 
 However, 
 it follows from
 results by Orey \& Pelikan [OrPe1,OrPe2]
 that the sequence
 $(\Pi\circ L_{n})_{n}$
satisfies the large deviation property
with rate constants $a_{n}=n$
and
Varadhan's
integral lemma can therefore be applied to 
provide information about the asymptotic behaviour of the sequence of integral
defined by
$$
\frac{1}{n}
\log
\int
\exp(nF)\,d(\Pi\circ L_{n}^{-1})\,.
\tag6.4
$$
In order to utilise the knowledge of the
asymptotic behaviour of (6.4)
for analysing the asymptotic behaviour of (6.3),
we must therefore show that the measures
 $$
 \Pi\circ\Lambda_{n}^{-1}
 $$
 and
  $$
  \Pi\circ L_{n}^{-1}
  $$
  are
  \lq\lq close".
 In fact, for technical reasons we will prove and use a similar statement involving the measures
  $$
  \Pi\circ M_{n}^{-1}
  $$
  and
  $$
  \Pi\circ L_{n}^{-1}\,.
  $$

   Indeed, below we prove a number
   of
   results
showing that the maps
$M_{n}$
  and
$L_{n}$
(and therefore also the
measures
$\Pi\circ M_{n}^{-1}$
  and
$\Pi\circ L_{n}^{-1}$)
are 
\lq\lq close".
These results
play an important role in
Section 7.
In particular, they
allow us to:
(1)
use 
Orey \& Pelikan's
result from [OrPe1,OrPe2] saying that
the sequence
$(\Pi\circ L_{n}^{-1})_{n}$
satisfies the large deviation property
to
prove
that
the sequence
$(\Pi\circ M_{n}^{-1})_{n}$
also satisfies the large deviation property (see Theorem 7.2),
and (2)
replace all occurencies
of 
$L_{n}(\overline{\bold i|n})$
in the formulas in Section 7 by
$M_{n}\bold i$
allowing us to use
the large deviation property of the sequence
$(\Pi\circ M_{n}^{-1})_{n}$.
 The above discussion 
 explains the mean reason for introducing the map $M_{n}$
 and the associated measure $\Pi\circ M_{n}^{-1}$.

\bigskip

{\bf Comparing $M_{n}$ and $L_{n}$.}
 We now prove various 
continuity statements saying that the maps
$M_{n}$
  and
$L_{n}$
are 
\lq\lq close".
These statements play an important role in Section 7
where we
apply
Varadhan's
integral lemma 
to prove Theorem 4.1.
We first introduce the
metric $\LDistance$
on 
$\Cal P(\Sigma_{\smallG}^{\Bbb N})$.
Fix $\gamma\in(0,1)$
and let $\distance_{\gamma}$
denote the metric on $\Sigma_{\smallG}^{\Bbb N}$
introduced in Section 5.2.
For a function $f:\Sigma_{\smallG}^{\Bbb N}\to\Bbb R$, we let
$\Lip_{\gamma}(f)$ denote the Lipschitz constant of $f$
with respect to the metric
$\distance_{\gamma}$, 
i\.e\.
$\Lip_{\gamma}(f)
=
\sup_{
\bold i,\bold j\in\Sigma_{\smallsmallG}^{\Bbb N},
\bold i\not=\bold j}
\frac
{|f(\bold i)-f(\bold j)|}
{\distance_{\gamma}(\bold i,\bold j)}$
and
we
define the metric $\LDistance$ in $\Cal P(\Sigma_{\smallG}^{\Bbb N})$ by
$$
\LDistance(\mu,\nu)
=
\sup
\Sb
f:\Sigma_{\smallsmallG}^{\Bbb N}\to\Bbb R\\
\Lip_{\gamma}(f)\le 1
\endSb
\Bigg|
\int f\,d\mu-\int f\,d\nu
\Bigg|;
\tag6.5
$$
we note that it is well-known that $\LDistance$ is a metric
and  that $\LDistance$ induces the weak topology.
Below we will always equip 
the space
$\Cal P(\Sigma_{\smallG}^{\Bbb N})$ 
with the metric $\LDistance$.
We can now state and prove the main results in this section.

\bigskip

\proclaim{Lemma 6.1}
Let $(X,\distance)$ be a metric space
and let $U:\Cal P(\Sigma_{\smallG}^{\Bbb N})\to X$ be continuous with respect to the weak topology.

\roster
\item"(1)"
We have
 $$
 \sup_{\bold u\in\Sigma_{\smallsmallG}^{n}}
 \sup
 \Sb
 \bold k,\bold l\in\Sigma_{\smallsmallG}^{\Bbb N}\\
 \termi(\bold u)=\ini(\bold k)\\
 \termi(\bold u)=\ini(\bold l)
 \endSb
  \distance
 \big(
 \,
 UL_{n}(\bold u\bold k)
 \,,\,
 UM_{n}(\bold u\bold l)
 \,
 \big)
 \to
 0
 \,\,\,\,
 \text{as $n\to\infty$.}
 $$
 
\item"(2)"
If $\varphi:\Sigma_{\smallG}^{\Bbb N}\to\Bbb R$ is a H\"older continuous function, then we have
 $$
 \sup_{\bold u\in\Sigma_{\smallsmallG}^{n}}
 \sup
 \Sb
 \bold k,\bold l\in\Sigma_{\smallsmallG}^{\Bbb N}\\
 \termi(\bold u)=\ini(\bold k)\\
 \termi(\bold u)=\ini(\bold l)
 \endSb
\Bigg|
\int\varphi\,d(L_{n}(\bold u\bold k))
-
\int\varphi\,d(M_{n}(\bold u\bold l))
\Bigg|
 \to
 0
 \,\,\,\,
 \text{as $n\to\infty$.}
 $$ 
\endroster

\endproclaim
\noindent{\it Proof}\newline
\noindent
We first note that
if
$\bold u\in\Sigma_{\smallG}^{n}$
and 
$\bold k\in\Sigma_{\smallG}^{\Bbb N}$
with
$\termi(\bold u)=\ini(\bold k)$,
and
$f:\Sigma_{\smallG}^{\Bbb N}\to\Bbb R$ is a continuous function, 
then
 $$
\align
\Bigg|
\int f\,d(L_{n}(\bold u\bold k))-\int f\,d(L_{n+|\,\widehat{\bold u}\,|}\overline{\bold u})
\Bigg|
&=
\Bigg|
\frac{1}{n}\sum_{i=0}^{n-1} f(S^{i}(\bold u\bold k))
-
\frac{1}{n+|\,\widehat{\bold u}\,|}\sum_{i=0}^{n+|\,\widehat{\bold u}\,|-1} 
f(S^{i}\overline{\bold u})
\Bigg|\\
&\le
\Bigg|
\frac{1}{n}\sum_{i=0}^{n-1} f(S^{i}(\bold u\bold k))
-
\frac{1}{n+|\,\widehat{\bold u}\,|}\sum_{i=0}^{n-1} f(S^{i}\overline{\bold u})
\Bigg|\\
&\qquad\qquad
\qquad\qquad
+
\Bigg|
\frac{1}{n+|\,\widehat{\bold u}\,|}
\sum_{i=n}^{n+|\,\widehat{\bold u}\,|-1} f(S^{i}\overline{\bold u})
\Bigg|
\\
&\le
\Bigg|
\frac{1}{n}\sum_{i=0}^{n-1} f(S^{i}(\bold u\bold k))
-
\frac{1}{n}\sum_{i=0}^{n-1} f(S^{i}\overline{\bold u})
\Bigg|\\
&\qquad\qquad
\qquad\qquad
+
\Bigg|
\frac{1}{n}\sum_{i=0}^{n-1} f(S^{i}\overline{\bold u})
-
\frac{1}{n+|\,\widehat{\bold u}\,|}\sum_{i=0}^{n-1} f(S^{i}\overline{\bold u})
\Bigg|\\
&\qquad\qquad
\qquad\qquad
+
\Bigg|
\frac{1}{n+|\,\widehat{\bold u}\,|}
\sum_{i=n}^{n+|\,\widehat{\bold u}\,|-1} f(S^{i}\overline{\bold u})
\Bigg|
\\
&\le
\frac{1}{n}
\sum_{i=0}^{n-1} 
|
f(S^{i}(\bold u\bold k))
-
f(S^{i}\overline{\bold u})
|
\\
&\qquad\qquad
\qquad\qquad
+
\frac{|\,\widehat{\bold u}\,|}{n+|\,\widehat{\bold u}\,|}
\frac{1}{n}\sum_{i=0}^{n-1}\|f\|_{\infty}\\
&\qquad\qquad
\qquad\qquad
+
\frac{1}{n+|\,\widehat{\bold u}\,|}
\sum_{i=n}^{n+|\,\widehat{\bold u}\,|-1}\|f\|_{\infty}
\\
&=
\frac{1}{n}
\sum_{i=0}^{n-1} 
|
f(S^{i}(\bold u\bold k))
-
f(S^{i}\overline{\bold u})
|
+
2
\frac{|\,\widehat{\bold u}\,|}{n+|\,\widehat{\bold u}\,|}
\|f\|_{\infty}\\
&\le
\frac{1}{n}
\sum_{i=0}^{n-1} 
|
f(S^{i}(\bold u\bold k))
-
f(S^{i}\overline{\bold u})
|
+
2
\frac{|\V| }{n}
\|f\|_{\infty}\,.
\quad
\text{[since $|\,\widehat{\bold u}\,|\le|\V|$]}
\tag6.6
\endalign
$$

\noindent (1).
Let $\varepsilon>0$.
Fix $\gamma\in(0,1)$ and let $\LDistance$ be the metric on 
$\Cal P(\Sigma_{\smallG}^{\Bbb N})$
defined in (6.5).
Since 
$U:\Cal P(\Sigma_{\smallG}^{\Bbb N})\to X$ is continuous and 
$\Cal P(\Sigma_{\smallG}^{\Bbb N})$ is compact, we
conclude that
$U:\Cal P(\Sigma_{\smallG}^{\Bbb N})\to X$
is uniformly continuous.
This implies that we can choose $\delta>0$
such that
all measures
$\mu,\nu\in\Cal P(\Sigma_{\smallG}^{\Bbb N})$
satisfy the following implication:
 $$
 \LDistance(\mu,\nu)<\delta
 \,\,\,\,
 \Rightarrow
 \,\,\,\,
 \distance(U\mu,U\nu)<\frac{\varepsilon}{2}\,.
 \tag6.7
 $$

Next,
choose a positive integer $N_{0}$ such that
$\frac{1}{N_{0}}(2 |\V|+\frac{1}{1-\gamma})
 <
\delta$\,.

Now, fix
$n\ge N_{0}$,
$\bold u\in\Sigma_{\smallG}^{n}$ and $\bold k,\bold l\in\Sigma_{\smallG}^{\Bbb N}$
with
$\termi(\bold u)=\ini(\bold k)$
and
$\termi(\bold u)=\ini(\bold l)$.
It follows that
$$
 \align
\LDistance
 \big(
 \,
 L_{n}(\bold u\bold k)
 \,,\,
 M_{n}(\bold u\bold l)
 \,
 \big)
&\le
 \LDistance
 \big(
 \,
 L_{n}(\bold u\bold k)
 \,,\,
 L_{n}(\bold u\bold l)
 \,
 \big)
+
\LDistance
 \big(
 \,
 L_{n}(\bold u\bold l)
 \,,\,
 M_{n}(\bold u\bold l)
 \,
 \big)\\
&\le
 \LDistance
 \big(
 \,
 L_{n}(\bold u\bold k)
 \,,\,
 L_{n}(\bold u\bold l)
 \,
 \big)
+
\LDistance
 \big(
 \,
 L_{n}(\bold u\bold l)
 \,,\,
 L_{n+|\,\widehat{\bold u}\,|}\overline{\bold u}
 \,
 \big) 
\endalign
$$

 We first estimate the distance
 $\LDistance
 \big(
 \,
 L_{n}(\bold u\bold k)
 \,,\,
 L_{n}(\bold u\bold l)
 \,
 \big)$.
 Indeed,
 since 
 $\frac{1}{N_{r}(1-\gamma)}
 \le
 \frac{1}{N_{0}}(2 |\V|+\frac{1}{1-\gamma})
 <
 \delta$,
 it follows that
 $$
 \align
\LDistance
 \big(
 \,
 L_{n}(\bold u\bold k)
 \,,\,
 L_{n}(\bold u\bold l)
 \,
 \big)
&=
\sup
\Sb
f:\Sigma_{\smallsmallG}^{\Bbb N}\to\Bbb R\\
\Lip_{\gamma}(f)\le 1
\endSb
\Bigg|
\int f\,d(L_{n}(\bold u\bold k))-\int f\,d(L_{n}(\bold u\bold l))
\Bigg|\\
&=
\sup
\Sb
f:\Sigma_{\smallsmallG}^{\Bbb N}\to\Bbb R\\
\Lip_{\gamma}(f)\le 1
\endSb
\Bigg|
\frac{1}{n}\sum_{i=0}^{n-1} f(S^{i}(\bold u\bold k))
-
\frac{1}{n}\sum_{i=0}^{n-1} f(S^{i}(\bold u\bold l))
\Bigg|\\
&\le
\sup
\Sb
f:\Sigma_{\smallsmallG}^{\Bbb N}\to\Bbb R\\
\Lip_{\gamma}(f)\le 1
\endSb
\frac{1}{n}\sum_{i=0}^{n-1} 
|
f(S^{i}(\bold u\bold k))
-
 f(S^{i}(\bold u\bold l))
|\\
&\le
\frac{1}{n}\sum_{i=0}^{n-1} 
\distance_{\gamma}
\big(
\,
S^{i}(\bold u\bold k)
\,,\,
S^{i}(\bold u\bold l)
\,
\big)\\
&=
\frac{1}{n}\sum_{i=0}^{n-1} 
\gamma^{
|S^{i}(\bold u\bold k)
\wedge
S^{i}(\bold u\bold l)|
}\\
&\le
\frac{1}{N_{0}}\sum_{i=0}^{n-1} 
\gamma^{n-i}\\
&\le
 \frac{1}{N_{0}(1-\gamma)}\\
&<
\delta\,,
\endalign
$$
and we therefore conclude from (6.7) that
 $$
 \distance
 \big(
 \,
 UL_{n}(\bold u\bold k)
 \,,\,
 UL_{n}(\bold u\bold l)
 \,
 \big)
 <
\frac{\varepsilon}{2}\,.
\tag6.8
$$

 Next, we estimate the distance
 $\LDistance
 \big(
 \,
 L_{n}(\bold u\bold l)
 \,,\,
 L_{n+|\,\widehat{\bold u}\,|}\overline{\bold u}
 \,
 \big)$.
We start by observing that
if we fix $\bold i_{0}\in\Sigma_{\smallG}^{\Bbb N}$, then
$$
\align
\LDistance(\mu,\nu)
&=
\sup
\Sb
f:\Sigma_{\smallsmallG}^{\Bbb N}\to\Bbb R\\
\Lip_{\gamma}(f)\le 1
\endSb
\Bigg|
\int f\,d\mu-\int f\,d\nu
\Bigg|\\
&=
\sup
\Sb
f:\Sigma_{\smallsmallG}^{\Bbb N}\to\Bbb R\\
\Lip_{\gamma}(f)\le 1
\endSb
\Bigg|
\int (f-f(\bold i_{0}))\,d\mu-\int (f-f(\bold i_{0}))\,d\nu
\Bigg|\\
&=
\sup
\Sb
g:\Sigma_{\smallsmallG}^{\Bbb N}\to\Bbb R\\
\Lip_{\gamma}(g)\le 1\\
g(\bold i_{0})=0
\endSb
\Bigg|
\int g\,d\mu-\int g\,d\nu
\Bigg|
\tag6.9
\endalign
$$
for all $\mu,\nu\in\Cal P(\Sigma_{\smallG}^{\Bbb N})$.
It  follows from 
(6.9) and (6.6)  that
 $$
 \align
\LDistance
 \big(
 \,
 L_{n}(\bold u\bold l)
 \,,\,
 L_{n+|\,\widehat{\bold u}\,|}\overline{\bold u}
 \,
 \big)
&=
\sup
\Sb
g:\Sigma_{\smallsmallG}^{\Bbb N}\to\Bbb R\\
\Lip_{\gamma}(g)\le 1\\
g(\bold i_{0})=0
\endSb
\Bigg|
\int g\,d(L_{n}(\bold u\bold l))-\int g\,d(L_{n+|\,\widehat{\bold u}\,|}\overline{\bold u})
\Bigg|\\
\allowdisplaybreak
&\le
\sup
\Sb
g:\Sigma_{\smallsmallG}^{\Bbb N}\to\Bbb R\\
\Lip_{\gamma}(g)\le 1\\
g(\bold i_{0})=0
\endSb
\Bigg(
\frac{1}{n}
\sum_{i=0}^{n-1} 
|
g(S^{i}(\bold u\bold l))
-
g(S^{i}\overline{\bold u})
|
+
2
\frac{|\V| }{n}
\|g\|_{\infty}
\Bigg)
\\
\allowdisplaybreak
&\le
\sup
\Sb
g:\Sigma_{\smallsmallG}^{\Bbb N}\to\Bbb R\\
\Lip_{\gamma}(g)\le 1\\
g(\bold i_{0})=0
\endSb
\Bigg(
\frac{1}{n}
\sum_{i=0}^{n-1} 
\Lip_{\gamma}(g)
\,\,
\distance_{\gamma}
\big(
\,
S^{i}(\bold u\bold l))
\,,\,
S^{i}\overline{\bold u}
\,
\big)
+
2
\frac{|\V| }{n}
\|g\|_{\infty}
\Bigg)
\\
\allowdisplaybreak
&\le
\sup
\Sb
g:\Sigma_{\smallsmallG}^{\Bbb N}\to\Bbb R\\
\Lip_{\gamma}(g)\le 1\\
g(\bold i_{0})=0
\endSb
\Bigg(
\frac{1}{n}
\sum_{i=0}^{n-1} 
\gamma^{
|S^{i}(\bold u\bold l)
\wedge
S^{i}\overline{\bold u}|
}
+
2
\frac{|\V| }{n}
\|g\|_{\infty}
\Bigg)
\\
\allowdisplaybreak
&\le
\sup
\Sb
g:\Sigma_{\smallsmallG}^{\Bbb N}\to\Bbb R\\
\Lip_{\gamma}(g)\le 1\\
g(\bold i_{0})=0
\endSb
\Bigg(
\frac{1}{n}
\sum_{i=0}^{n-1} 
\gamma^{n-i}
+
2
\frac{|\V| }{n}
\|g\|_{\infty}
\Bigg)
\\
\allowdisplaybreak
&\le
\sup
\Sb
g:\Sigma_{\smallsmallG}^{\Bbb N}\to\Bbb R\\
\Lip_{\gamma}(g)\le 1\\
g(\bold i_{0})=0
\endSb
\Bigg(
\frac{1}{n(1-\gamma)}
+
2
\frac{|\V| }{n}
\|g\|_{\infty}
\Bigg)\,.
\tag6.10
\endalign
$$
However, 
if
$g:\Sigma_{\smallG}^{\Bbb N}\to\Bbb R$ satisfies
$\Lip_{\gamma}(g)\le 1$ and
$g(\bold i_{0})=0$,
then
$|g(\bold i)|
=
|g(\bold i)-g(\bold i_{0})|
\le
\distance_{\gamma}(\bold i,\bold i_{0})
\le
1$
for all $\bold i\in\Sigma_{\smallG}^{\Bbb N}$, whence
$\|g\|_{\infty}\le 1$.
It therefore follows from (6.10) that
if 
$n\ge N_{0}$,
$\bold u\in\Sigma_{\smallG}^{n}$ and $\bold l\in\Sigma_{\smallG}^{\Bbb N}$
with
$\termi(\bold u)=\ini(\bold l)$, then
 $\LDistance
 \big(
 \,
 L_{n}(\bold u\bold l)
 \,,\,
 L_{n+|\,\widehat{\bold u}\,|}\overline{\bold u}
 \,
 \big)
 \le
\frac{1}{n} (\frac{1}{1-\gamma}
 +
 2|\V|)
 \le
\frac{1}{N_{0}} (\frac{1}{1-\gamma}
 +
 2|\V|)
 <
 \delta$,
 and we therefore conclude from (6.7) that
 $$
 \distance
 \big(
 \,
 UL_{n}(\bold u\bold l)
 \,,\,
 L_{n+|\,\widehat{\bold u}\,|}\overline{\bold u}
 \,
 \big)
 <
 \frac{\varepsilon}{2}\,.
 \tag6.11
 $$

Finally, combining (6.8) and (6.11) shows 
that
$\distance
 \big(
 \,
 UL_{n}(\bold u\bold k)
 \,,\,
 UM_{n}(\bold u\bold l)
 \,
 \big)
\le
 \distance
 \big(
 \,
 UL_{n}(\bold u\bold k)
 \,,\,
 UL_{n}(\bold u\bold l)
 \,
 \big)
+
\distance
 \big(
 \,
 UL_{n}(\bold u\bold l)
 \,,\,
 UL_{n+|\,\widehat{\bold u}\,|}\overline{\bold u}
 \,
 \big)
 <
 \frac{\varepsilon}{2}
 +
 \frac{\varepsilon}{2}
 =
 \varepsilon$.

 \noindent (2)
 Since $\varphi$ is H\"older continuous 
 there are constant $c$ and $a$ with $c,a>0$ such that
$|\varphi(\bold i)-\varphi(\bold j)|
\le
c\,\distance_{\gamma}(\bold i,\bold j)^{a}$ for all 
$\bold i,\bold j\in\Sigma_{\smallG}^{\Bbb N}$.
It follows from this and (6.6) that
 $$
 \align
 \sup_{\bold u\in\Sigma_{\smallsmallG}^{n}}
 \sup
 \Sb
 \bold k,\bold l\in\Sigma_{\smallsmallG}^{\Bbb N}\\
 \termi(\bold u)=\ini(\bold k)\\
 \termi(\bold u)=\ini(\bold l)
 \endSb
&\Bigg|
\int\varphi\,d(L_{n}(\bold u\bold k))
-
\int\varphi\,d(M_{n}(\bold u\bold l))
\Bigg|\\
&=
 \sup_{\bold u\in\Sigma_{\smallsmallG}^{n}}
 \sup
 \Sb
 \bold k,\bold l\in\Sigma_{\smallsmallG}^{\Bbb N}\\
 \termi(\bold u)=\ini(\bold k)\\
 \termi(\bold u)=\ini(\bold l)
 \endSb
\Bigg|
\int\varphi\,d(L_{n}(\bold u\bold k))
-
\int\varphi\,d(L_{n+|\,\widehat{\bold u}\,|}\overline{\bold u})
\Bigg|\\
\allowdisplaybreak
&\le
 \sup_{\bold u\in\Sigma_{\smallsmallG}^{n}}
 \sup
 \Sb
 \bold k,\bold l\in\Sigma_{\smallsmallG}^{\Bbb N}\\
 \termi(\bold u)=\ini(\bold k)\\
 \termi(\bold u)=\ini(\bold l)
 \endSb
\Bigg(
\frac{1}{n}
\sum_{i=0}^{n-1} 
|
\varphi(S^{i}(\bold u\bold k))
-
\varphi(S^{i}\overline{\bold u})
|
+
2
\frac{|\V| }{n}
\|\varphi\|_{\infty}
\Bigg)\\
\allowdisplaybreak
&\le
 \sup_{\bold u\in\Sigma_{\smallsmallG}^{n}}
 \sup
 \Sb
 \bold k,\bold l\in\Sigma_{\smallsmallG}^{\Bbb N}\\
 \termi(\bold u)=\ini(\bold k)\\
 \termi(\bold u)=\ini(\bold l)
 \endSb
\Bigg(
\frac{1}{n}
\sum_{i=0}^{n-1} 
c
\,
\distance_{\gamma}
\big(\,
S^{i}(\bold u\bold k))
\,,\,
S^{i}\overline{\bold u}
\,
\big)^{a}
+
2
\frac{|\V| }{n}
\|\varphi\|_{\infty}
\Bigg)\\
\allowdisplaybreak
&\le
 \sup_{\bold u\in\Sigma_{\smallsmallG}^{n}}
 \sup
 \Sb
 \bold k,\bold l\in\Sigma_{\smallsmallG}^{\Bbb N}\\
 \termi(\bold u)=\ini(\bold k)\\
 \termi(\bold u)=\ini(\bold l)
 \endSb
\Bigg(
\frac{1}{n}
\sum_{i=0}^{n-1} 
c
\,
\gamma^{
a
|S^{i}(\bold u\bold l)
\wedge
S^{i}\overline{\bold u}|
}
+
2
\frac{|\V| }{n}
\|\varphi\|_{\infty}
\Bigg)\\
\allowdisplaybreak
&\le
 \sup_{\bold u\in\Sigma_{\smallsmallG}^{n}}
 \sup
 \Sb
 \bold k,\bold l\in\Sigma_{\smallsmallG}^{\Bbb N}\\
 \termi(\bold u)=\ini(\bold k)\\
 \termi(\bold u)=\ini(\bold l)
 \endSb
\Bigg(
\frac{1}{n}
\sum_{i=0}^{n-1} 
c
\,
\gamma^{
a(n-i)
}
+
2
\frac{|\V| }{n}
\|\varphi\|_{\infty}
\Bigg)\\
&\le
\frac{c}{n}
\,
\frac{1}{1-\gamma^{a}}
+
2
\frac{|\V| }{n}
\|\varphi\|_{\infty}
\to
0
\,\,\,\,\text{as $n\to\infty$.}
\endalign
 $$
This completes the proof. 
  \hfill$\square$

 \bigskip

\proclaim{Lemma 6.2}
Let $X$ be an inner product space with inner product 
$\langle\cdot|\cdot\rangle$
and let $U:\Cal P(\Sigma_{\smallG}^{\Bbb N})\to X$ be continuous with respect to the weak topology.
Fix a H\"older continuous function $\varphi:\Sigma_{\smallG}^{\Bbb N}\to\Bbb R$
and
let $s$ be a real number and $q\in X$.
We have
 $$
 \sup_{\bold i\in\Sigma_{\smallsmallG}^{\Bbb N}}
 \Bigg|
 \Bigg(
 \sup_{\bold u\in[\overline{\bold i|n}]}
 \langle q|UL_{n}\bold u\rangle+s
 \Bigg)
 \int\varphi\,d\big(L_{n}\big(\overline{\bold i|n}\big)\big)
 -
 (
 \langle q|UM_{n}\bold i\rangle+s
 )
 \int\varphi\,d(M_{n}\bold i)
 \Bigg|
  \to
 0
 \,\,\,\,
 \text{as $n\to\infty$.}
 $$

\endproclaim
\noindent{\it Proof}\newline
\noindent
We note that if
$\bold u\in\Sigma_{\smallsmallG}^{n}$
and
$\bold k,\bold l,\bold m\in\Sigma_{\smallsmallG}^{\Bbb N}$
with
$\termi(\bold u)=\ini(\bold k)$,
$\termi(\bold u)=\ini(\bold l)$
and
$\termi(\bold u)=\ini(\bold m)$
then
Cauchy--Schwarz inequality implies that
$$
\align
 \Bigg|
 (
 \langle q|UL_{n}(\bold u\bold m)\rangle
& +s
 \big)
 \int\varphi\,d(L_{n}(\bold u\bold k))
 -
 (
 \langle q|UM_{n}(\bold u\bold l)\rangle
 +s
 )
 \int\varphi\,d(M_{n}(\bold u\bold l))
 \Bigg|\\
 &\le
 \Bigg|
 (
 \langle q|UL_{n}(\bold u\bold m)\rangle
 +s
 )
 \int\varphi\,d(L_{n}(\bold u\bold k))
 -
 (
 \langle q|UM_{n}(\bold u\bold l)\rangle
 +s
 )
 \int\varphi\,d(L_{n}(\bold u\bold k))
 \Bigg|\\
 &\qquad
 +
  \Bigg|
 (
 \langle q|UM_{n}(\bold u\bold l)\rangle
 +s
 )
 \int\varphi\,d(L_{n}(\bold u\bold k))
 -
 (
 \langle q|UM_{n}(\bold u\bold l)\rangle
 +s
 )
 \int\varphi\,d(M_{n}(\bold u\bold l))
 \Bigg|\\
&\le
 (
 |
 \langle q|UL_{n}(\bold u\bold m)-UM_{n}(\bold u\bold l)\rangle
 |
 +
 |s|
 )
 \|\varphi\|_{\infty}\\
 &\qquad
 +
 (
 |
 \langle q|UM_{n}(\bold u\bold l)\rangle
 |
 +
 |s|
 )
 \Bigg|
 \int\varphi\,d(L_{n}(\bold u\bold k))
 -
 \int\varphi\,d(M_{n}(\bold u\bold l))
 \Bigg|\\
 &\le
 (
\|q\|
\,
\|UL_{n}(\bold u\bold m)-UM_{n}(\bold u\bold l)\|
 +
 |s|
 )
 \|\varphi\|_{\infty}\\
 &\qquad
 +
 (
 \|q\|
 \,
 \|U\|_{\infty}
 +
 |s|
 )
 \Bigg|
 \int\varphi\,d(L_{n}(\bold u\bold k))
 -
 \int\varphi\,d(M_{n}(\bold u\bold l))
 \Bigg|\,,
 \tag6.12
\endalign
$$
where $\|U\|_{\infty}<\infty$
since $U:\Cal P(\Sigma_{\smallG}^{\Bbb N})\to X$ is continuous 
and 
$\Cal P(\Sigma_{\smallG}^{\Bbb N})$
is compact with respect to the weak topology.
The desired result follows immediately from 
Lemma 6.1 and 
(6.12) .
\hfill$\square$

 \bigskip

\proclaim{Lemma 6.3}
Let $X$ be an inner product space with inner product 
$\langle\cdot|\cdot\rangle$
and let $U:\Cal P(\Sigma_{\smallG}^{\Bbb N})\to X$ be continuous with respect to the weak topology.
Let $q\in X$.
We have
$$
 \sup_{\bold u\in\Sigma_{\smallsmallG}^{n}}
 \sup
 \Sb
 \bold k,\bold l\in\Sigma_{\smallsmallG}^{\Bbb N}\\
 \termi(\bold u)=\ini(\bold k)\\
 \termi(\bold u)=\ini(\bold l)
 \endSb
|
 \langle q|UL_{n}(\bold u\bold k)\rangle
 -
  \langle q|UL_{n}(\bold u\bold l)\rangle
|
 \to
 0
 \,\,\,\,
 \text{as $n\to\infty$.}
 $$

\endproclaim
\noindent{\it Proof}\newline
\noindent
We note that if
$\bold u\in\Sigma_{\smallsmallG}^{n}$
and
$\bold k,\bold l\in\Sigma_{\smallsmallG}^{\Bbb N}$
with
$\termi(\bold u)=\ini(\bold k)$
and
$\termi(\bold u)=\ini(\bold l)$,
then
$M_{n}(\bold u\bold k)=M_{n}(\bold u\bold l)$,
and
Cauchy--Schwarz inequality therefore implies that
$$
\align
 |
 \langle q|UL_{n}(\bold u\bold k)\rangle
 -
  \langle q|UL_{n}(\bold u\bold l)\rangle
|
 &\le
 \|q\|
 \,
 \|UL_{n}(\bold u\bold k)-UL_{n}(\bold u\bold l)\|\\
 &\le
 \|q\|
 \,
\Big(
 \|UL_{n}(\bold u\bold k)-UM_{n}(\bold u\bold k)\|
 +
 \|UM_{n}(\bold u\bold k)-UM_{n}(\bold u\bold l)\|\\
 &\qquad\qquad
 \qquad\qquad
 \qquad\qquad
 \qquad\qquad
 +
 \|UM_{n}(\bold u\bold l)-UL_{n}(\bold u\bold l)\|
 \Big)\\
&\le
 \|q\|
 \,
\Big(
 \|UL_{n}(\bold u\bold k)-UM_{n}(\bold u\bold k)\|
 +
 \|UM_{n}(\bold u\bold l)-UL_{n}(\bold u\bold l)\|
 \Big)\,.
\endalign
$$
The desired result follows immediately from 
this and
Lemma 6.1.
\hfill$\square$

  \bigskip


\heading
{
7. Proof of Theorems 4.1 and Theorem 4.2, Part 2: The proofs.}
\endheading

The purpose of this section is to use the continuity results in  Section 6
and
Varadhan's integral lemma from large
deviation theory to prove 
Theorem 4.1.
We first introduce the sequence of measures
(see (7.3) below)
 that
we will
apply Varadhan's integral lemma
to.

\bigskip

{\bf The measure $\Pi$.}
Let
$B=(b_{\,\smallvertexi,\smallvertexj})_{\smallvertexi,\smallvertexj\in\smallV}$
denote the matrix defined by
 $$
 \align
 b_{\,\smallvertexi,\smallvertexj}
&=
 |\E_{\smallvertexi,\smallvertexj}|\,;
 \endalign
  $$
  recall, that
  $\E_{\smallvertexi,\smallvertexj}$ denotes the set of edges from $\vertexi$
  to
  $\vertexj$.
We denote the spectral radius of $B$ by
$\lambda$.
Since $\G=(\V,\E)$ is strongly connected,
we conclude that the matrix $B$ is irreducible,
and it therefore follows from the Perron-Frobenius theorem that there 
is a unique
right eigen-vector 
$\bold u=(u_{\smallvertexi})_{\smallvertexi\in\smallV}$ 
of $B$ with 
eigen-value $\lambda$
and a 
unique
left eigen-vector
$\bold v=(v_{\smallvertexi})_{\smallvertexi\in\smallV}$ 
of $B$ with 
eigen-value $\lambda$, i\.e\.
$$
 \aligned
 \bold u B
&=
 \lambda\bold u\,,\\
 B\bold v
&=\lambda \bold v\,,
 \endaligned
 \tag7.1
 $$
such that
$u_{\smallvertexi},v_{\smallvertexi}>0$
for all $\vertexi$,
$\sum_{\smallvertexi}u_{\smallvertexi}v_{\smallvertexi}=1$
and
$\sum_{\smallvertexi}u_{\smallvertexi}=1$.
For $\edge\in\V$, write
 $p_{\smalledge}
 =
 v_{\ini(\smalledge)}^{-1}
 \,
 v_{\termi(\smalledge)}
 \,
 \lambda^{-1}$.
A simple calculation shows that
$\sum_{
 \smalledge\in\smallE_{\smallsmallvertexi}
 }
 p_{\smalledge} =
 1$
 for all $\vertexi$
and
that
$\sum_{\smallvertexi}
\sum_{
 \smalledge\in\smallE_{\smallsmallvertexi,\smallsmallvertexj}
 }
 u_{\smallvertexi}
 \,
 v_{\smallvertexi}
 \,
 p_{\smalledge}
 =
 u_{\smallvertexj}
 \,
 v_{\smallvertexj}$
  for all $\vertexj$.
It follows from this
that there is
 a unique Borel probability measure
 $\Pi\in\Cal P(\Sigma_{\smallG}^{\Bbb N})$
 such that
 $$
 \align
 \Pi[\bold i]
 &=
  u_{\ini(\smalledge_{1})}
 \,
 v_{\ini(\smalledge_{1})}
 \,
  p_{\smalledge_{1}}
  \,
  \cdots
  \,
 p_{\smalledge_{n}} \\
&=
u_{\ini(\smalledge_{1})}
 \,
 v_{\termi(\smalledge_{n})}
 \,
 \lambda^{-n}\\
&=
u_{\ini(\bold i)}
 \,
 v_{\termi(\bold i)}
 \,
 \lambda^{-n}
 \tag7.2
\endalign 
 $$
for
all
$\bold i
=
\edge_{1}\ldots \edge_{n}\in\Sigma_{\smallG}^{*}$.

\bigskip

{\bf The measure $\Pi_{n}$.}
Finally, for a positive integer $n$, we define the probability measures
$\Pi_{n}\in\Cal P(\Cal P(\Sigma_{\smallG}^{\Bbb N}))$ by
 $$
 \Pi_{n}
 =
 \Pi\circ M_{n}^{-1}\,;
 \tag7.3
 $$
recall, that the map
$M_{n}:\Sigma_{\smallG}^{\Bbb N}
\to
\Cal P(\Sigma_{\smallG}^{\Bbb N})$
is defined in (6.1).

\bigskip

We now turn towards the proof of the main result in
  this section, namely, Theorem 4.1.
  The proof
  of Theorem 4.1 uses
  Varadhan's integral lemma from large
deviation theory 
and the fact that the sequence
  $(\Pi_{n})_{n}$ has the large deviation property.
  We now state the definition of the
  large deviation property
  and
  Varadhan's integral lemma.

\bigskip

\proclaim{Definition}
Let $X$
be a complete separable metric space
and
let $(P_{n})_{n}$ be a sequence of probability measures on $X$. 
Let
$(a_{n})_{n}$ be a sequence of positive numbers with
$a_{n}\to\infty$ 
and let
$I:X\to[0,\infty]$ be a lower semi-continuous function
with compact level sets.
The sequence $(P_{n})_{n}$ is said to have the large deviation property
with rate constants $(a_{n})_{n}$ and rate function $I$ if the following 
two condistions
hold:
\roster
\item"(i)" For each closed subset $K$ of $X$, we have
 $$
 \limsup_{n}\frac{1}{a_{n}}\log P_{n}(K)\le-\inf_{x\in K}I(x)\,;
 $$
\item"(ii)" For each open subset $G$ of $X$, we have
 $$
 \liminf_{n}\frac{1}{a_{n}}\log P_{n}(G)\ge-\inf_{x\in G}I(x)\,.
 $$
\endroster
\endproclaim

%
%
%
%

\bigskip

\proclaim{Theorem 7.1. 
Varadhan's integral lemma}
Let $X$ be a complete separable metric space
and
let $(P_{n})_{n}$ be a sequence of probability measures on $X$.
Assume that the sequence $(P_{n})_{n}$
has the large deviation property
with rate constants $(a_{n})_{n}$ and rate function $I$.
Let $F:X\to\Bbb R$ be a continuous function
satisfying the following two conditions:
\roster
\item"(i)"
For all $n$, we have
 $$
 \int\exp(a_{n}F)\,dP_{n}
 <
 \infty\,.
 $$
\item"(ii)"
We have
 $$
 \lim_{M\to\infty}\,\,
 \limsup_{n}\,\,
 \frac{1}{a_{n}}
 \log
 \int_{\{M\le F\}}
 \exp(a_{n}F)\,dP_{n}
 =
 -\infty\,.
 $$
 \endroster
(Observe that the  Conditions (i)--(ii)
 are satisfied if $F$ is bounded.)
Then we have
 $$
 \lim_{n}\,\,
 \frac{1}{a_{n}}
 \log
 \int
 \exp(a_{n}F)\,dP_{n}
 =
 -\inf_{x\in X}(I(x)-F(x))\,.
 $$
\endproclaim 
\noindent{\it Proof}\newline
\noindent 
This  follows from [El, Theorem II.7.1] or [DeZe, Theorem 4.3.1].
\hfill$\square$

 \bigskip

\noindent
The next result says that the
sequence
$(\Pi_{n})_{n}$ has the large deviation property.

\bigskip

\proclaim{Theorem 7.2}
Define $I:\Cal P(\Sigma_{\smallG}^{\Bbb N})\to[0,\infty]$ by
 $$
 I(\mu)
 =
 \cases
 \log\lambda-h(\mu)
\quad
 &\text{for $\mu\in \Cal P_{S}(\Sigma_{\smallG}^{\Bbb N})$;}\\
 \infty
\quad
 &\text{for 
 $\mu\in \Cal P(\Sigma_{\smallG}^{\Bbb N})\setminus \Cal P_{S}(\Sigma_{\smallG}^{\Bbb N})$.
 }
 \endcases
$$
Then the 
sequence
$(\Pi_{n})_{n}$
has the large deviation deviation property with respect to the sequence
$(n)_{n}$ and rate function $I$.
\endproclaim
\noindent{\it Proof}\newline
\noindent
Let
 $$
 \Gamma_{n}
 =
 \Pi\circ L_{n}^{-1}\,.
 $$
It follows from
 Orey \& Pelikan [OrPe1,OrPe2]
 (see also [JiQiQi, Remark 7.2.2])
 that the sequence $(\Gamma_{n})_{n}$
 has the large deviation 
 property
 with respect to the sequence
$(n)_{n}$ and rate function $I$.
Using the results from Section 6
(saying that the measures
$\Pi_{n}=\Pi\circ M_{n}^{-1}$
and
$\Gamma_{n}=\Pi\circ L_{n}^{-1}$
are \lq\lq close")
it can be shown that
that this implies that the 
sequence $(\Pi_{n})_{n}$
also
 has the large deviation 
 property
 with respect to the sequence
$(n)_{n}$ and rate function $I$;
the reader is 
referred to [MiOl2]
for a detailled
proof of this.
\hfill$\square$

\bigskip

We now turn towards the proof of Theorem 4.1.

\bigskip

\proclaim{Theorem 7.3}
Let $X$ be an inner product space
 with inner product $\langle\cdot|\cdot\rangle$
  and let $U:\Cal P(\Sigma_{\smallG}^{\Bbb N})\to X$ be 
continuous with respect to the weak topology.
Fix a H\"older continuous function $\varphi:\Sigma_{\smallG}^{\Bbb N}\to\Bbb R$
and let $s$ be a real number
and $q\in X$.
For each positive integer $n$ and each 
$ \bold i\in\Sigma_{\smallG}^{n}$, 
we write
 $$
 \align
 u_{\bold i}
&=
\sup_{\bold u\in[\bold i]}\langle q|UL_{n}\bold u\,\rangle\,,\\
s_{\bold i}
&=
\sup_{\bold u\in[\bold i]}
\exp
 \sum_{k=0}^{n-1}\varphi S^{k}\bold u\,.
 \endalign
 $$
 Then
 $$
 \align
  \,
  \liminf_{n}
 \,
  \,\,
  \frac{1}{n}
 \,\,
 \log
   \sum
  \Sb
 \bold i\in\Sigma_{\smallsmallG}^{n}
  \endSb
 \exp
((u_{\bold i}+s)s_{\bold i})k
&=
\sup
\Sb
\mu\in\Cal P_{S}(\Sigma_{\smallsmallG}^{\Bbb N})\\
\endSb
\Bigg(
h(\mu)+(\langle q|U\mu\rangle+s)\int\varphi\,d\mu
\Bigg)\,,\\
  \limsup_{n}
  \,\,
  \frac{1}{n}
 \,\,
 \log
   \sum
  \Sb
 \bold i\in\Sigma_{\smallsmallG}^{n}
  \endSb
 \exp
((u_{\bold i}+s)s_{\bold i})
&=
\sup
\Sb
\mu\in\Cal P_{S}(\Sigma_{\smallsmallG}^{\Bbb N})\\
\endSb
\Bigg(
h(\mu)+(\langle q|U\mu\rangle+s)\int\varphi\,d\mu
\Bigg)\,.
\endalign
$$ 
\endproclaim 
\noindent{\it Proof}\newline
\noindent
For brevity write
$t_{\bold i}=\exp((u_{i}+s)s_{\bold i})$.
We also define $F:\Cal P(\Sigma_{\smallG}^{\Bbb N})\to\Bbb R$ by
 $$
 F(\mu)
 =
 (\langle q|U\mu\rangle+s)
 \int\varphi\,d\mu\,.
 $$
Observe that since $\varphi$ and $U$ are bounded, 
i\.e\. $\|\varphi\|_{\infty}<\infty$ and $\|U\|_{\infty}<\infty$,
we conclude that
$\|F\|_{\infty}
\le
(\|q\|\,\|U\|_{\infty}+|s|)\|\varphi\|_{\infty}<\infty$.

Next, we prove the following claim.

\bigskip

{\it Claim 1.
Then there is a constant $c$
such that
the following holds:
for all $\varepsilon>0$,
there is a positive integer $N_{\varepsilon}$,
such that
if $n\ge N_{\varepsilon}$, then
$$
 \align
 \sum
   \Sb
   \bold i\in\Sigma_{\smallsmallG}^{n}\\
   	 \endSb
	t_{\bold i}
	&\le
c
\,\,
\lambda^{n}
\,\,\,\,\,\,
e^{n\varepsilon}
\,\,\,\,
\int\exp(nF)\,d\Pi_{n}\,,
\tag7.4\\
 \sum
   \Sb
   \bold i\in\Sigma_{\smallsmallG}^{n}\\
   	 \endSb
t_{\bold i}
&\ge
\frac{1}{c}
\,\,
\lambda^{n}
\,\,
e^{-n\varepsilon}
\,\,
\int\exp(nF)\,d\Pi_{n}\,.
\tag7.5
\endalign
 $$
}

{\it Proof of Claim 1.}
For each positive integer $n$, we clearly have
 $$
 \align
\int
 t_{\bold i|n}
 \,d\Pi(\bold i)
&=
 \sum
   \Sb
  \bold k\in\Sigma_{\smallsmallG}^{n}\\
 	 \endSb
	 \,\,
\int\limits_{[\bold k]}
 t_{\bold i|n}
 \,d\Pi(\bold i)\\
&=
 \sum
   \Sb
  \bold k\in\Sigma_{\smallsmallG}^{n}\\
	 \endSb
 t_{\bold k}
 \,\,
 \Pi\big(\,
 [\bold k]
 \,\big)\\
&=
 \sum
   \Sb
  \bold k\in\Sigma_{\smallsmallG}^{n}\\
  	 \endSb
 t_{\bold k}
 \,
 u_{\ini(\bold k)}
 \,
 v_{\termi(\bold k)}
 \,
 \lambda^{-n}\,.
 \tag7.6
 \endalign
 $$

Next,
it follows from the 
Principle of
Bounded Distortion 
(see, for example, [Bar,Fa])
that there is a constant $C>0$
such that
if $n\in\Bbb N$, $\bold i\in\Sigma_{\smallG}^{n}$
and $\bold u,\bold v\in[\bold i]$, then
$|
\sum_{k=0}^{n-1}\varphi S^{k}\bold u
-
\sum_{k=0}^{n-1}\varphi S^{k}\bold v
|
\le
C$.
In particular, this implies that
for all positive integers $n$
and all
$\bold i\in\Sigma_{\smallG}^{\Bbb N}$, we have
$|
(u_{\bold i}+s)s_{\bold i}
-
(u_{\bold i}+s)\sum_{k=0}^{n-1}\varphi S^{k}\overline{\bold i}
|
\le
(|u_{\bold i}|+|s|)C
=
(|\sup_{\bold u\in[\bold i]}\langle q|UL_{n}\bold u\,\rangle|+|s|)c_{0}
\le
(\|q\|\,\|U\|_{\infty}+|s|)C$,
whence
 $$
 \gathered
 \frac{1}{c_{0}}
 \exp
 \Bigg(
  \Bigg(
  \sup_{\bold u\in[\bold i]}
	 \langle q|UL_{n}\bold u\,\rangle
	 +
	 s
	 \Bigg)
 \,\,
 \sum_{k=0}^{n-1}\varphi S^{k}\overline{\bold i}
 \Bigg)
 \le
 t_{\bold i}\,,\\
 t_{\bold i
 }\le
 c_{0}
  \exp
 \Bigg(
  \Bigg(
  \sup_{\bold u\in[\bold i]}
	 \langle q|UL_{n}\bold u\,\rangle
	 +
	 s
	 \Bigg)
 \,\,
 \sum_{k=0}^{n-1}\varphi S^{k}\overline{\bold i}
 \Bigg)\,,
 \endgathered
 \tag7.7
 $$
 for all positive integers $n$
and all
$\bold i\in\Sigma_{\smallG}^{\Bbb N}$
where
$c_{0}
=
\exp(\|q\|\,\|U\|_{\infty}+|s|)C$.

We also note that it  follows from Lemma 6.2 that there is a positive integer $N_{\varepsilon}$
such that
 $$
 \gathered
 \Bigg(
 \sup_{\bold u\in[\overline{\bold i|n}]}
 \langle q|UL_{n}\bold u\rangle+s
 \Bigg)
 \int\varphi\,d\big(L_{n}\big(\overline{\bold i|n}\big)\big)
 \le
 \varepsilon
 +
 (
 \langle q|UM_{n}\bold i\rangle+s
 )
 \int\varphi\,d(M_{n}\bold i)\,,\\
 (
 \langle q|UM_{n}\bold i\rangle+s
 )
 \int\varphi\,d(M_{n}\bold i)
 \le
 \varepsilon
 +
 \Bigg(
 \sup_{\bold u\in[\overline{\bold i|n}]}
 \langle q|UL_{n}\bold u\rangle+s
 \Bigg)
 \int\varphi\,d\big(L_{n}\big(\overline{\bold i|n}\big)\big)\,,
 \endgathered
 \tag7.8
  $$
 for all integers $n$ with $n\ge N_{\varepsilon}$
and all
$\bold i\in\Sigma_{\smallG}^{\Bbb N}$.

Finally, we can
find a constant $w>0$
such that
$\frac{1}{w}
\le
u_{\smallvertexi}
v_{\smallvertexi}
\le
w$
for all $\vertexi$.
Now put $c=c_{0}w$.

It follows from (7.6), (7.7) and (7.8) that
if $n$ is a positive integer with $n\ge N_{\varepsilon}$, then
we have
 $$
 \align
  \sum
   \Sb
\bold k\in\Sigma_{\smallsmallG}^{n}\\
	 \endSb
 t_{\bold k}
 &\le
 w
 \,
 \lambda^{n}
 \,
\int
 t_{\bold i|n}
 \,d\Pi(\bold i)\\
&\le 
 c
 \,
 \lambda^{n}
 \,
 \int
  \exp
 \Bigg(
  \Bigg(
 \sup_{\bold u\in[\overline{\bold i|n}]}
 \langle q|UL_{n}\bold u\rangle+s
 \Bigg)
 \,\,
 \sum_{k=0}^{n-1}\varphi S^{k}\big(\overline{\bold i|n}\big)
 \Bigg)
  \,d\Pi(\bold i)\\
 &=
 c
 \,
 \lambda^{n}
 \,
 \int
  \exp
 \Bigg(
 n
 \,
 \Bigg(
 \sup_{\bold u\in[\overline{\bold i|n}]}
 \langle q|UL_{n}\bold u\rangle+s
 \Bigg)
 \int\varphi\,d\big(L_{n}\big(\overline{\bold i|n}\big)\big)
 \Bigg)
  \,d\Pi(\bold i)\\ 
\allowdisplaybreak 
 &\le
 c
 \,
 \lambda^{n}
 \,
 \int
  \exp
 \Bigg(
 n
 \,
 \Bigg(
  \varepsilon
 +
 (
 \langle q|UM_{n}\bold i\rangle+s
 )
 \int\varphi\,d(M_{n}\bold i)
 \Bigg)
 \Bigg)
  \,d\Pi(\bold i)\\ 
\allowdisplaybreak 
&=
 c
 \,
 \lambda^{n}
 \,
 e^{n\varepsilon}
 \,
 \int
  \exp\left(n F(M_{n}\bold i)\right)
 \,d\Pi(\bold i)\,.
 \endalign
 $$
This proves inequality (7.4).
Inequality (7.5) is proved similarly.
This completes the proof of Claim 1.

\bigskip

We now turn towards the proof the statement
in the theorem.
Let $\varepsilon>0$.
We first observe 
that it follows immediately from Claim 1 that
 $$
 \aligned
 \liminf_{n}
 \frac{1}{n}
 \log
 \sum
   \Sb
\bold i\in\Sigma_{\smallsmallG}^{n}\\
	 \endSb
 t_{\bold i}
&\ge
 \log \lambda
 \,-\,
 \varepsilon
 \,+\,
\, \liminf_{n}\,
 \frac{1}{n}
 \log 
 \int
 \exp
 \left(
 nF
 \right)
 \,d\Pi_{n}\,,\\
\limsup_{n}
 \frac{1}{n}
 \log
 \sum
   \Sb
\bold i\in\Sigma_{\smallsmallG}^{n}\\
	 \endSb
 t_{\bold i}
&\le
 \log \lambda
 \,+\,
 \varepsilon
 \,+\,
  \limsup_{n}
 \frac{1}{n}
 \log 
 \int
 \exp
 \left(
 nF
 \right)
 \,d\Pi_{n}\,.
 \endaligned
 \tag7.9
 $$
Since inequalities (7.9) hold for all
$\varepsilon>0$, we now conclude that
 $$
 \aligned
 \liminf_{n}
 \frac{1}{n}
 \log
 \sum
   \Sb
\bold i\in\Sigma_{\smallsmallG}^{n}\\
	 \endSb
 t_{\bold i}
&\ge
 \log \lambda
 \,+\,
\, \liminf_{n}\,
 \frac{1}{n}
 \log 
 \int
 \exp
 \left(
 nF
 \right)
 \,d\Pi_{n}\,,\\
\limsup_{n}
 \frac{1}{n}
 \log
 \sum
   \Sb
\bold i\in\Sigma_{\smallsmallG}^{n}\\
	 \endSb
 t_{\bold i}
&\le
 \log \lambda
 \,+\,
  \limsup_{n}
 \frac{1}{n}
 \log 
 \int
 \exp
 \left(
 nF
 \right)
 \,d\Pi_{n}\,.
 \endaligned
 \tag7.10
 $$

Next,
define $I:\Cal P(\Sigma_{\smallG}^{\Bbb N})\to[0,\infty]$ by
 $$
 I(\mu)
 =
 \cases
 \log\lambda-h(\mu)
\quad
 &\text{for $\mu\in \Cal P_{S}(\Sigma_{\smallG}^{\Bbb N})$;}\\
 \infty
\quad
 &\text{for 
 $\mu\in \Cal P(\Sigma_{\smallG}^{\Bbb N})\setminus \Cal P_{S}(\Sigma_{\smallG}^{\Bbb N})$.
 }
 \endcases
$$
It follows from Theorem 7.2
that the sequence 
$(\Pi_{n})_{n}
\subseteq
\Cal P\big(\,\Cal P(\Sigma_{\smallG}^{\Bbb N})\,\big)$ 
has the large deviation property with respect to
the sequence
$(n)_{n}$ and rate function $I$,
and
we therefore 
conclude from Theorem 7.1 that
 $$
 \align
  \lim_{n}
 \frac{1}{n}
 \log 
 \int
 \exp
 \left(
 nF
 \right)
 \,d\Pi_{n}
&=
 -
 \,
  \inf_{\nu\in\Cal P(\Sigma_{\smallsmallG}^{\Bbb N})}(I(\nu)-F(\nu))\\
&=
 -
  \inf_{\nu\in\Cal P_{S}(\Sigma_{\smallsmallG}^{\Bbb N})}(I(\nu)-F(\nu))\,.
  \tag7.11
  \endalign
 $$

Finally, we introduce the following notation and prove 
the two claims below.
Write
 $$
 \align
 \underline Q_{q}^{U}(\varphi)
&=
  \,
  \liminf_{n}
 \,
  \,\,
  \frac{1}{n}
 \,\,
 \log
   \sum
  \Sb
 \bold i\in\Sigma_{\smallsmallG}^{n}
  \endSb
 \exp
((u_{\bold i}+s)s_{\bold i})\,,\\
 \overline Q_{q}^{U}(\varphi)
&=
  \limsup_{n}
  \,\,
  \frac{1}{n}
 \,\,
 \log
   \sum
  \Sb
 \bold i\in\Sigma_{\smallsmallG}^{n}
  \endSb
 \exp
((u_{\bold i}+s)s_{\bold i})\,.
\endalign
$$

\bigskip

{\it Claim 2. We have
$ \underline Q_{q}^{U}(\varphi)
\ge
 \sup_{
	\mu\in\Cal P_{S}(\Sigma_{\smallsmallG}^{\Bbb N})}
(h(\mu)+(\langle q|U\mu\rangle+s)\int\varphi\,d\mu)$.
}

{\it Proof of Claim 2.}
Combining (7.10) and (7.11) now yields
 $$
 \align
 \underline Q_{q}^{U}(\varphi)
 &=
 \liminf_{n}
 \frac{1}{n}
 \log
 \sum
   \Sb
   \bold i\in\Sigma_{\smallsmallG}^{n}\\
	 \endSb
 t_{\bold i}\\
&\ge
 \log \lambda
 \,+\,
 \liminf_{n}
 \frac{1}{n}
 \log 
 \int
 \exp
 \left(
 nF
 \right)
 \,d\Pi_{n}\\
&\ge
 \log \lambda
 \,-\,
 \inf_{\nu\in\Cal P_{S}(\Sigma_{\smallsmallG}^{\Bbb N})}(I(\nu)-F(\nu))\\
&=
 \log \lambda
 \,+\,
 \sup
  \Sb
	\mu\in\Cal P_{S}(\Sigma_{\smallsmallG}^{\Bbb N})\\
	\endSb
 (F(\mu)-I(\mu))\\
&=
 \sup
  \Sb
	\mu\in\Cal P_{S}(\Sigma_{\smallsmallG}^{\Bbb N})\\
	\endSb
 \left(h(\mu)+(\langle q|U\mu\rangle+s)\int\varphi\,d\mu\right)\,. 
 \endalign
 $$
 This completes the proof of Claim 2.

 \bigskip

{\it Claim 3. We have
 $\overline Q_{q}^{U}(\varphi)
\le
 \sup_{
	\mu\in\Cal P_{S}(\Sigma_{\smallsmallG}^{\Bbb N})}
(h(\mu)+(\langle q|U\mu\rangle+s)\int\varphi\,d\mu)$.
}

{\it Proof of Claim 3.}
The proof of Claim 3 is similar to the proof of Claim 2
and is therefore omitted.
This completes the proof of Claim 3.

\bigskip

The statement in Theorem 7.3
now follows immediately by combining Claim 2 and Claim 3.
\hfill$\square$

\bigskip

We can now prove Theorem 4.1
and 
Theorem 4.2.
We first prove Theorem 4.1;
for the benefit of the reader we have decided to repeat the statement of Theorem 4.1.

\bigskip

\proclaim{Theorem 4.1. The variational principle for $\tau$}
Let $X$ be an inner product space
with
inner product
$\langle\cdot|\cdot\rangle$
and let $U:\Cal P(\Sigma_{\smallG}^{\Bbb N})\to X$ be 
continuous with respect to the weak topology.
Fix a H\"older continuous function $\varphi:\Sigma_{\smallG}^{\Bbb N}\to\Bbb R$ with $\varphi<0$
and
let $q\in X$.

\roster
\item"(1)"
We have
 $$
 \tau(q)
 =
\sup
\Sb
\mu\in\Cal P_{S}(\Sigma_{\smallsmallG}^{\Bbb N})
\endSb
\Bigg(
\,
-
\frac{h(\mu)}{\int\varphi\,d\mu}
-
\big\langle q\big|U\mu \big\rangle
\,
\Bigg)\,.
$$ 
\item"(2)"
We have
  $$
\sup
\Sb
\mu\in\Cal P_{S}(\Sigma_{\smallsmallG}^{\Bbb N})
\endSb
\Bigg(
\,
h(\mu)
+
\big(
\big\langle q\big|U\mu \big\rangle
+
\tau(q)
\big)
\int\varphi\,d\mu
\,
\Bigg)
=
0\,.
$$ 
\endroster
\endproclaim

\noindent{\it Proof}\newline

\noindent
(1) For brevity write
 $$
 u
 =
  \sup
  \Sb
  \mu\in\Cal P_{S}(\Sigma^{\Bbb N})\\
  \endSb
  \Bigg(
 -\frac{h(\mu)}{\int \varphi\,d\mu}-\langle q|U\mu\rangle
 \Bigg)\,.
 $$
 We will also use the same notation as in Theorem 7.3.
Namely,
for each positive integer $n$ and each 
$ \bold i\in\Sigma_{\smallG}^{n}$, 
we write
$u_{\bold i}
=
\sup_{\bold u\in[\bold i]}
\langle q|UL_{n}\bold u\,\rangle$
and
$s_{\bold i}
=
\sup_{\bold u\in[\bold i]}
\exp
 \sum_{k=0}^{n-1}\varphi S^{k}\bold u$.
 Since $\varphi:\Sigma_{\smallG}^{\Bbb N}\to\Bbb N$ is continuous with $\varphi<0$,
 we 
 also conclude that there is a constant $c>0$, such that
  $$
\varphi\le -c\,.
  $$
We must now prove the following two inequalities
 $$
 \gather
 \tau(q)\le u\,,\tag7.12\\
 u\le \tau(q)\,.\tag7.13
 \endgather
 $$

{\it Proof of (7.12)}.
We must prove that if $s>u$, then 
  $$
 \sum
  \Sb
  \bold i\\
    \endSb
 \exp((u_{\bold i}+s)s_{\bold i})
 <
 \infty\,.
 $$

Let $s>u$ and write $\varepsilon=\frac{s-u}{3}>0$.
It follows from the 
definition
of $u$ that if
$\mu\in\Cal P_{S}(\Sigma^{\Bbb N})$,
then
 we have
$-\frac{h(\mu)}{\int\varphi\,d\mu}
-\langle q|U\mu\rangle
<
u+\varepsilon$, 
whence
$h(\mu)+\langle q|U\mu\rangle\int\varphi\,d\mu
<
-(u+\varepsilon)\int\varphi\,d\mu$
where
we have used the fact that
$\int\varphi\,d\mu<0$
because 
$\varphi<0$. This implies that
 if 
$\mu\in\Cal P_{S}(\Sigma^{\Bbb N})$,
then
 $$
 \align
 h(\mu)
 +
 (\langle q|U\mu\rangle+s)\int\varphi\,d\mu
 &=
 h(\mu)
 +
 (\langle q|U\mu\rangle+(u+3\varepsilon))\int\varphi\,d\mu\\
 &\le
 2\varepsilon\int\varphi\,d\mu\\
&\le
 -2\varepsilon c\,.
 \endalign
 $$
We deduce from this inequality and Theorem 7.3 that
 $$
 \align
 \limsup_{n}
\frac{1}{n}
 \log
&\sum
   \Sb
  \bold i\in\Sigma_{\smallsmallG}^{n}\\
	 \endSb
  \exp((u_{\bold i}+s)s_{\bold i})\\
&=
  \sup
  \Sb
  \mu\in\Cal P_{S}(\Sigma_{\smallsmallG}^{\Bbb N})\\
  \endSb
 \Bigg(
 h(\mu)
 +
 (\langle q|U\mu\rangle +s)\int\varphi\,d\mu
 \Bigg)
  \qquad\qquad
  \text{[by Theorem 7.3]}\\
 &\le
	-2\varepsilon c\\
&<
 -\varepsilon c\,.	
 \tag7.14
 \endalign	
 $$
Inequality (7.14) shows
that
there is an  integer $N_{0}$ such that 
$ \frac{1}{n}
 \log
 \sum_{
   \bold i\in\Sigma_{\smallsmallG}^{n}
   }
 \exp((u_{\bold i}+s)s_{\bold i})
 \le
  -\varepsilon c$
for all $n\ge N_{0}$, whence
 $$
  \sum
   \Sb
  \bold i\in\Sigma_{\smallsmallG}^{n}	 \endSb
 \exp((u_{\bold i}+s)s_{\bold i})
 \le
 e^{-\varepsilon cn}
 \tag7.15
 $$
for all $n\ge N_{0}$. 
Using (7.15) we now conclude that
 $$
 \align
 \sum
  \Sb
 \bold i
   \endSb
\exp((u_{\bold i}+s)s_{\bold i})
&=
 \sum_{n<N_{0}}
 \sum
   \Sb
 \bold i\in\Sigma_{\smallsmallG}^{n}
	 \endSb
 \exp((u_{\bold i}+s)s_{\bold i})
 \,\,
 +
 \,\,
  \sum_{n\ge N_{0}}
 \sum
   \Sb
  \bold i\in\Sigma_{\smallsmallG}^{n}
	 \endSb
\exp((u_{\bold i}+s)s_{\bold i})\\
&{}\\ 
&\le
 \sum_{n<N_{0}}
 \sum
   \Sb
   \bold i\in\Sigma_{\smallsmallG}^{n}
	 \endSb
 \exp((u_{\bold i}+s)s_{\bold i})
 \,\,
 +
 \,\,
  \sum_{n\ge N_{0}}
e^{-\varepsilon cn}\\
&{}\\
&<
\infty\,.
\endalign
$$
This completes the proof of (7.12).

{\it Proof of 7.13).}
We must prove that if $s<u$, then 
  $$
 \sum
  \Sb
  \bold i\\
    \endSb
 \exp((u_{\bold i}+s)s_{\bold i})
 =
 \infty\,.
 $$

Let $s<u$ and write $\varepsilon=\frac{u-s}{3}>0$.
It follows from the 
definition
of $u$ that
there is a measure 
$\mu_{\varepsilon}\in\Cal P_{S}(\Sigma^{\Bbb N})$,
such that
$u-\varepsilon
\le
-\frac{h(\mu_{\varepsilon})}{\int\varphi\,d\mu_{\varepsilon}}
-\langle q|U\mu_{\varepsilon}\rangle$, 
whence
$h(\mu_{\varepsilon})+\langle q|U\mu_{\varepsilon}\rangle\int\varphi\,d\mu_{\varepsilon}
\ge
-(u-\varepsilon)\int\varphi\,d\mu_{\varepsilon}$
where
we have used the fact that
$\int\varphi\,d\mu_{\varepsilon}<0$
because 
$\varphi<0$. This implies that
 $$
 \align
 h(\mu_{\varepsilon})
 +
 (\langle q|U\mu_{\varepsilon}\rangle+s)\int\varphi\,d\mu_{\varepsilon}
 &=
 h(\mu_{\varepsilon})
 +
 (\langle q|U\mu_{\varepsilon}\rangle+(u-3\varepsilon))\int\varphi\,d\mu_{\varepsilon}\\
 &\ge
 -2\varepsilon\int\varphi\,d\mu_{\varepsilon}\\
&\ge
 2\varepsilon c\,.
 \endalign
 $$
We deduce from this inequality and Theorem 7.3 that
 $$
 \align
 \limsup_{n}
\frac{1}{n}
 \log
&\sum
   \Sb
  \bold i\in\Sigma_{\smallsmallG}^{n}\\
	 \endSb
  \exp((u_{\bold i}+s)s_{\bold i})\\
&=
  \sup
  \Sb
  \mu\in\Cal P_{S}(\Sigma_{\smallsmallG}^{\Bbb N})\\
  \endSb
 \Bigg(
 h(\mu)
 +
 (\langle q|U\mu\rangle +s)\int\varphi\,d\mu
 \Bigg)
  \qquad\qquad
  \text{[by Theorem 7.3]}\\
 &\ge
 h(\mu_{\varepsilon})
 +
 (\langle q|U\mu_{\varepsilon}\rangle+s)\int\varphi\,d\mu_{\varepsilon} \\
 &\ge
	2\varepsilon c\\
&>
 \varepsilon c\,.	
 \tag7.16
 \endalign	
 $$
Inequality (7.16) shows
that
there is an  integer $N_{0}$ such that 
$ \frac{1}{n}
 \log
 \sum_{
   \bold i\in\Sigma_{\smallsmallG}^{n}
   }
 \exp((u_{\bold i}+s)s_{\bold i})
 \ge
  \varepsilon c$
for all $n\ge N_{0}$, whence
 $$
  \sum
   \Sb
  \bold i\in\Sigma_{\smallsmallG}^{n}	 \endSb
 \exp((u_{\bold i}+s)s_{\bold i})
 \ge
 e^{\varepsilon cn}
 \tag7.17
 $$
for all $n\ge N_{0}$. 
Using (7.17) we now conclude that
 $$
 \align
 \sum
  \Sb
 \bold i
   \endSb
\exp((u_{\bold i}+s)s_{\bold i})
&=
 \sum_{n<N_{0}}
 \sum
   \Sb
 \bold i\in\Sigma_{\smallsmallG}^{n}
	 \endSb
 \exp((u_{\bold i}+s)s_{\bold i})
 \,\,
 +
 \,\,
  \sum_{n\ge N_{0}}
 \sum
   \Sb
  \bold i\in\Sigma_{\smallsmallG}^{n}
	 \endSb
\exp((u_{\bold i}+s)s_{\bold i})\\
&{}\\ 
&\ge
 \sum_{n<N_{0}}
 \sum
   \Sb
   \bold i\in\Sigma_{\smallsmallG}^{n}
	 \endSb
 \exp((u_{\bold i}+s)s_{\bold i})
 \,\,
 +
 \,\,
  \sum_{n\ge N_{0}}
e^{\varepsilon cn}\\
&{}\\
&=
\infty\,.
\endalign
$$
This completes the proof of (7.13).

\noindent
(2) This statement follows easily from (1).
\hfill$\square$

\bigskip

Finally, we prove Theorem 4.2.

\bigskip

%
%
\noindent{\it Proof of theorem 4.2}\newline
\noindent
We denote the inner product in $X$ by
$\langle\cdot|\cdot\rangle$.

\noindent
(1)
Fix $q,p\in X$ and 
$s,t\in[0,1]$ with
$s+t=1$
and let $\varepsilon>0$.
Note that it follows from Lemma 6.3
that we can find a positive integer $N_{\varepsilon}$ such that
if 
$n\ge N_{\varepsilon}$ and
$\bold i\in\Sigma_{\smallG}^{n}$, then
$\langle q|UL_{n}\bold u\rangle
\ge
\langle q|UL_{n}\bold v\rangle-\frac{\varepsilon}{2}$
and
$\langle p|UL_{n}\bold u\rangle
\ge
\langle p|UL_{n}\bold w\rangle-\frac{\varepsilon}{2}$
for all
$\bold u,\bold v,\bold w\in[\bold i]$.
Hence, if
$n\ge N_{\varepsilon}$ and
$\bold i\in\Sigma_{\smallG}^{n}$, then
$\langle sq+tp|UL_{n}\bold u\rangle
=
s\langle q|UL_{n}\bold u\rangle
+
t\langle p|UL_{n}\bold u\rangle
\ge
s(\langle q|UL_{n}\bold v\rangle
-
\frac{\varepsilon}{2})
+
t(\langle p|UL_{n}\bold w\rangle
-
\frac{\varepsilon}{2})
=
s \langle q|UL_{n}\bold v\rangle
+
t \langle p|UL_{n}\bold w\rangle
-
\frac{\varepsilon}{2}$
for all
$\bold u,\bold v,\bold w\in[\bold i]$.
This implies that if
$n\ge N_{\varepsilon}$
and $\bold i\in\Sigma_{\smallG}^{n}$, then

 $$
 \align
 \sup_{\bold u\in[\bold i]}\langle sq+tp|UL_{|\bold i|}\bold u\rangle
&\ge
s \sup_{\bold u\in[\bold i]}\langle q|UL_{|\bold i|}\bold u\rangle
+
t \sup_{\bold u\in[\bold i]}\langle p|UL_{|\bold i|}\bold u\rangle
-
\frac{\varepsilon}{2}\\
\endalign
$$ 
Writing
$$
M_{\varepsilon}
=
 \sum_{n<N_{\varepsilon}}
  \sum_{\bold i\in\Sigma_{\smallsmallG}^{n}}
 \exp
 \Bigg(
 \Bigg(
 \sup_{\bold u\in[\bold i]}
\langle sq+tp|UL_{|\bold i|}\bold u \rangle
+
s\tau(q)+t\tau(p)+\varepsilon
\Bigg)
 \sup_{\bold u\in[i]}
 \sum_{k=0}^{n-1}\varphi(S^{k}\bold u)
 \Bigg)
 $$
and
using H\"older's inequality
together
with the fact that $\varphi<0$,
we therefore conclude that
 $$
 \align
 \zeta_{sq+tp}^{\co,U}
 &(\varphi;
 s\tau(q)+t\tau(p)+\varepsilon)\\
 &=
  \sum_{\bold i}
 \exp
 \Bigg(
 \Bigg(
 \sup_{\bold u\in[\bold i]}
\langle sq+tp|UL_{|\bold i|}\bold u \rangle
+
s\tau(q)+t\tau(p)+\varepsilon
\Bigg)
 \sup_{\bold u\in[i]}
 \sum_{k=0}^{|\bold i|-1}\varphi(S^{k}\bold u)
 \Bigg)\\
 &=
 M_{\varepsilon}
  +
 \sum_{n\ge N_{\varepsilon}}
  \sum_{\bold i\in\Sigma_{\smallsmallG}^{n}}
 \exp
 \Bigg(
 \Bigg(
 \sup_{\bold u\in[\bold i]}
\langle sq+tp|UL_{|\bold i|}\bold u \rangle
+
s\tau(q)+t\tau(p)+\varepsilon
\Bigg)
 \sup_{\bold u\in[i]}
 \sum_{k=0}^{n-1}\varphi(S^{k}\bold u)
 \Bigg)\\ 
 &\le
M_{\varepsilon}
  +
 \sum_{n\ge N_{\varepsilon}}
  \sum_{\bold i\in\Sigma_{\smallsmallG}^{n}}
 \exp
 \Bigg(
 \Bigg(
s \sup_{\bold u\in[\bold i]}\langle q|UL_{|\bold i|}\bold u\rangle
+
t \sup_{\bold u\in[\bold i]}\langle p|UL_{|\bold i|}\bold u\rangle\\
&\qquad\qquad
\qquad\qquad
\qquad\qquad
+
s\tau(q)+t\tau(p)+
\tfrac{\varepsilon}{2}
\Bigg)
 \sup_{\bold u\in[i]}
 \sum_{k=0}^{n-1}\varphi(S^{k}\bold u)
 \Bigg)\\ 
 &=
 M_{\varepsilon}
  +
 \sum_{n\ge N_{\varepsilon}}
  \sum_{\bold i\in\Sigma_{\smallsmallG}^{n}}
  \Bigg(
  \,\,
 \exp
 \Bigg(
 \Bigg(
\sup_{\bold u\in[\bold i]}\langle q|UL_{|\bold i|}\bold u\rangle
+
\tau(q)+
\tfrac{\varepsilon}{2}
\Bigg)
 \sup_{\bold u\in[i]}
 \sum_{k=0}^{|\bold i|-1}\varphi(S^{k}\bold u)
 \Bigg)^{s}\\
&\qquad\qquad 
  \qquad\qquad
  \qquad\qquad
  \times
   \exp
 \Bigg(
 \Bigg(
\sup_{\bold u\in[\bold i]}\langle p|UL_{|\bold i|}\bold u\rangle
+
\tau(p)+
\tfrac{\varepsilon}{2}
\Bigg)
 \sup_{\bold u\in[i]}
 \sum_{k=0}^{|\bold i|-1}\varphi(S^{k}\bold u)
 \Bigg)^{t}
 \,\,
 \Bigg)\\
  &\le
  M_{\varepsilon}
  +
  \Bigg(
  \,\,
  \sum_{n\ge N_{\varepsilon}}
  \sum_{\bold i\in\Sigma_{\smallsmallG}^{n}}
 \exp
 \Bigg(
 \Bigg(
\sup_{\bold u\in[\bold i]}\langle q|UL_{|\bold i|}\bold u\rangle
+
\tau(q)
+
\tfrac{\varepsilon}{2}
\Bigg)
 \sup_{\bold u\in[i]}
 \sum_{k=0}^{|\bold i|-1}\varphi(S^{k}\bold u)
 \Bigg)
 \,\,
 \Bigg)^{s}\\
&\qquad\qquad 
  \qquad\qquad
  \times
  \Bigg(
  \,\,
  \sum_{n\ge N_{\varepsilon}}
  \sum_{\bold i\in\Sigma_{\smallsmallG}^{n}}
   \exp
 \Bigg(
 \Bigg(
\sup_{\bold u\in[\bold i]}\langle p|UL_{|\bold i|}\bold u\rangle
+
\tau(p)
+
\tfrac{\varepsilon}{2}
\Bigg)
 \sup_{\bold u\in[i]}
 \sum_{k=0}^{|\bold i|-1}\varphi(S^{k}\bold u)
 \Bigg)
 \,\,
 \Bigg)^{t}\\
  &\le
 M_{\varepsilon}
 +
 \zeta_{q}^{\co,U}
 (\varphi;
 \tau(q)+\tfrac{\varepsilon}{2})^{s}
 \,\,
 \zeta_{p}^{\co,U}
 (\varphi;
 \tau(p)+\tfrac{\varepsilon}{2})^{t}\\
 &{}
  \tag7.18
 \endalign
  $$
 Since $\varepsilon>0$, 
it follows from the definition of $\tau(q)$ and $\tau(p)$
that
$ \zeta_{q}^{\co,U}
 (\varphi;
 \tau(q)+\frac{\varepsilon}{2})<\infty$
and
$\zeta_{p}^{\co,U}
 (\varphi;
 \tau(p)+\frac{\varepsilon}{2})<\infty$,
 and we therefore conclude from 
 (7.18) that
 $\zeta_{sq+tp}^{\co,U}
 (\varphi;s\tau(q)+t\tau(p)+\varepsilon)
\le
M_{\varepsilon}
+
 \zeta_{q}^{\co,U}
 (\varphi;
 \tau(q)+\frac{\varepsilon}{2})^{s}
 \,\,
 \zeta_{p}^{\co,U}
 (\varphi;
 \tau(p)+\frac{\varepsilon}{2})^{t}
 <\infty$.
This inequality shows that
$\tau(sq+tp)
\le
s\tau(q)+t\tau(p)+\varepsilon$ for all $\varepsilon>0$.
Finally,
letting $\varepsilon$ tend to $0$
gives
$\tau(sq+tp)
\le
s\tau(q)+t\tau(p)$.

\noindent
(2)
For a H\"older continuous function 
$f:\Sigma_{\smallG}^{\Bbb N}\to\Bbb R$, we write $P(f)$ for the pressure of
$f$.
Next, define $\Phi:\Bbb R^{M+1}\to\Bbb R$ by
  $$
  \Phi(q,s)
  =
  \sup
  _{\mu\in\Cal P_{S}(\Sigma_{\smallsmallG}^{\Bbb N})}
  \Bigg(
  h(\mu)
  +
  (
  \langle q|U\mu\rangle
  +
  s
  )
  \int\varphi\,d\mu
  \Bigg)
  $$
  for $q\in\Bbb R^{M}$ and $s\in\Bbb R$.
Observe, that
since $U\mu=\frac{\int\psi\,d\mu}{\int\varphi\,d\mu}$ for 
$\mu\in\Cal P_{S}(\Sigma_{\smallG}^{\Bbb N})$, 
we conclude that
 $$
 \align
 \Phi(q,s)
 &=
 \sup_{\mu\in\Cal P_{S}(\Sigma_{\smallsmallG}^{\Bbb N})}
 \Bigg(
 h(\mu)
 +
 \Bigg(
 \Bigg\langle q\Bigg|\frac{\int\psi\,d\mu}{\int\varphi\,d\mu} \Bigg\rangle+s
 \Bigg)
 \int\varphi\,d\mu
  \Bigg)\\  
 &=
  \sup_{\mu\in\Cal P_{S}(\Sigma_{\smallsmallG}^{\Bbb N})}
 \Bigg(
 h(\mu)
 +
 \int
 (\langle q|\psi\rangle+s\varphi)\,d\mu
 \Bigg)\,.
 \tag7.19
 \endalign
 $$
 It follows from (7.19)
 and
 the Variational Principle (see [Wa]) that
 $\Phi(q,s)$ equals the pressure of 
 $\langle q|\psi\rangle+s\varphi$, i\.e\.
 $\Phi(q,s)
 =
  P(\langle q|\psi\rangle+s\varphi)$,
  and we therefore conclude
  from [Rue1]
  that $\Phi$ is  real analytic.
  Since also Theorem 4.1 implies that
  $\Phi(q,\tau(q))=0$ for all $q$, we now conclude from the implicit function theorem
  that $\tau $ is real analytic.

\noindent
(3)
It follows from
Theorem 4.1 that
$\tau(0)
=
\sup_{\mu\in\Cal P_{S}(\Sigma_{\smallsmallG}^{\Bbb N})}-\frac{h(\mu)}{\int\Lambda\,d\mu}$.
The desired result follows from this since it is well-known that if the OSC is satisfied,
then
$\dim_{\Haus}K_{\smallvertexi}
=
\sup_{\mu\in\Cal P_{S}(\Sigma_{\smallsmallG}^{\Bbb N})}-\frac{h(\mu)}{\int\Lambda\,d\mu}$
for all $\vertexi\in\V$
(indeed, if we write
$t=\dim_{\Haus}K_{\smallvertexi}$,
then it follows 
from Bowen's formula (see [Bar,Fa])
that
$0
=
P(t\Lambda)
=
\sup_{\mu\in\Cal P_{S}(\Sigma_{\smallsmallG}^{\Bbb N})}(h(\mu)+t\int\Lambda\,d\mu)$,
whence
$t
=
\sup_{\mu\in\Cal P_{S}(\Sigma_{\smallsmallG}^{\Bbb N})}-\frac{h(\mu)}{\int\Lambda\,d\mu}$).
\hfill$\square$

  \bigskip


\heading
{
8. Proof of Theorem 4.3}
\endheading

The purpose of this section is to prove 
Theorem 4.3.
We start by proving the following
auxiliary
result.

\bigskip

\proclaim{Proposition 8.1}
Let $X$ be an inner product space and let $U:\Cal P(\Sigma_{\smallG}^{\Bbb N})\to X$ be 
continuous with respect to the weak topology.
Fix a continuous function  $\varphi:\Sigma_{\smallG}^{\Bbb N}\to\Bbb R$ 
with
$\varphi<0$.
Let $\alpha\in X$ 
and let
$\,\,\scri f(\alpha)$
be the unique real number
such that
 $$
 \align
 \limsup_{r\searrow 0}
 \sigma_{\radius}
 \big(
 \,
 \zeta_{B(\alpha,r)}^{\dyn,U}(\,\,\scri f(\alpha)\,\varphi;\cdot)
 \,
 \big)
 &=
 1\,.
 \endalign
 $$ 
Then
 $$
  \scri f(\alpha)
  \le
  \tau^{*}(\alpha)\,.
  $$

\endproclaim 

\noindent{\it Proof}\newline
\noindent
Let $\langle\cdot|\cdot\rangle$
denote the inner product in $X$.
Fix $q\in X$.
Next, let $\varepsilon>0$ and fix 
$r>0$ with
$0<r<\frac{\varepsilon}{\|q\|}$
where we write
$\frac{\varepsilon}{\|q\|}=\infty$ if $q=0$.
We now have
 $$
 \align
 \zeta_{q}^{\co,U}
 (\varphi;-\langle q|\alpha\rangle+
 &\,\,\scri f(\alpha)-\varepsilon)\\
 &=
 \qquad\,
  \sum_{\bold i}
 \qquad\,
 \exp
 \Bigg(
 \Bigg(
\sup_{\bold u\in[\bold i]} \langle q|UL_{|\bold i|}\bold u \rangle
-\langle q|\alpha\rangle+\,\,\scri f(\alpha)-\varepsilon
\Bigg)
 \sup_{\bold u\in[i]}
 \sum_{k=0}^{|\bold i|-1}\varphi(S^{k}\bold u)
 \Bigg)\\
 &\ge
  \sum
  \Sb
  \bold i\\
  {}\\
  UL_{|\bold i|}[\bold i]\subseteq B(\alpha,r)
  \endSb
 \exp
 \Bigg(
  \Bigg(
\sup_{\bold u\in[\bold i]} \langle q|UL_{|\bold i|}\bold u \rangle
-\langle q|\alpha\rangle+\,\,\scri f(\alpha)-\varepsilon
\Bigg)
 \sup_{\bold u\in[\bold i]}
 \sum_{k=0}^{|\bold i|-1}\varphi(S^{k}\bold u)
 \Bigg)\,.\\
 &{}
 \tag8.1
 \endalign
  $$
Next, note that if
$\bold i\in\Sigma_{\smallG}^{*}$ with
$UL_{|\bold i|}[\bold i]\subseteq B(\alpha,r)$ and $\bold u\in[\bold i]$, then
$UL_{|\bold i|}\bold u\in UL_{|\bold i|}[\bold i]\subseteq B(\alpha,r)$,
whence
$\|UL_{|\bold i|}\bold u-\alpha\|\le r$, 
and so
$| \langle q|UL_{|\bold i|}\bold u\rangle
-
\langle q|\alpha\rangle|
=
|
 \langle q|UL_{|\bold i|}\bold u-\alpha\rangle
|
\le
\|q\|\,\|UL_{|\bold i|}\bold u -\alpha\|
\le
\|q\|\,r$.
Hence,
 if
$\bold i\in\Sigma_{\smallG}^{*}$ with
$UL_{|\bold i|}[\bold i]\subseteq B(\alpha,r)$, then
(using the fact that $\varphi<0$)
 $$
 \align
  \exp
 \Bigg(
 \Bigg(
\sup_{\bold u\in[\bold i]} \langle q|UL_{|\bold i|}\bold u \rangle
-\langle q|\alpha\rangle
&+
\,\,\scri f(\alpha)-\varepsilon
\Bigg)
 \sup_{\bold u\in[\bold i]}
 \sum_{k=0}^{|\bold i|-1}\varphi(S^{k}\bold u)
 \Bigg)\\
&\ge
   \exp
 \Bigg(
 \big(
\|q\|r 
+
\,\,\scri f(\alpha)-\varepsilon
\big)
 \sup_{\bold u\in[\bold i]}
 \sum_{k=0}^{|\bold i|-1}\varphi(S^{k}\bold u)
 \Bigg)\\
&=
   \exp
 \Bigg(
 \big(
\,\,\scri f(\alpha)-(\varepsilon-\|q\|r)
\big)
 \sup_{\bold u\in[\bold i]}
 \sum_{k=0}^{|\bold i|-1}\varphi(S^{k}\bold u)
 \Bigg)\,.
 \tag8.2
 \endalign
 $$
Combining (8.1) and (8.2) now shows that
$$
 \zeta_{q}^{\co,U}
 (\varphi;-\langle q|\alpha\rangle+
 \,\,\scri f(\alpha)-\varepsilon)
 \ge
  \sum
  \Sb
  \bold i\\
  {}\\
  UL_{|\bold i|}[\bold i]\subseteq B(\alpha,r)
  \endSb
  \exp
 \Bigg(
 \big(
\,\,\scri f(\alpha)-(\varepsilon-\|q\|r)
\big)
 \sup_{\bold u\in[\bold i]}
 \sum_{k=0}^{|\bold i|-1}\varphi(S^{k}\bold u)
 \Bigg)\,.
 \tag8.3
  $$

Next,
we observe
that if $\eta>0$, then
it follows 
immediately 
from the definition
of $\,\,\scri f(\alpha)$ that
$\limsup_{r\searrow 0}
 \sigma_{\radius}
 \big(
 \,
 \zeta_{B(\alpha,r)}^{\dyn,U}((\,\,\scri f(\alpha)-\eta)\,\varphi;\cdot)
 \,
 \big)
 <
 1$,
 and we can therefore find
positive real numbers
$\delta_{\eta}$
and
 $\rho_{\eta}$
such that
if $0<\rho<\rho_{\eta}$, then
$\sigma_{\radius}
 \big(
 \,
 \zeta_{B(\alpha,\rho)}^{\dyn,U}((\,\,\scri f(\alpha)-\eta)\,\varphi;\cdot)
 \,
 \big)
 <
 1-\delta_{\eta}$.
 In particular, this
shows that 
if $0<\rho<\rho_{\eta}$ and $t>1-\delta_{\eta}$, then
 $$
 \sum_{n}
\frac{t^{n}}{n!}
  \sum
  \Sb
  \bold i\in\Sigma_{\smallsmallG}^{n}\\
  {}\\
  UL_{|\bold i|}[\bold i]\subseteq B(\alpha,\rho)
  \endSb
  \exp
 \Bigg(
 \big(
\,\,\scri f(\alpha)-\eta
\big)
 \sup_{\bold u\in[\bold i]}
 \sum_{k=0}^{|\bold i|-1}\Phi(S^{k}\bold u)
 \Bigg) 
 =
 \infty\,.
 \tag8.4
  $$

Since $\varepsilon-\|q\|r>0$, we
can
choset a positive real number $\rho$ with $0<\rho<\min(\rho_{\varepsilon-\|q\|r},r)$.
It now follows from (8.3)  and (8.4) that
$$
 \align
 \zeta_{q}^{\co,U}
 (\varphi;-\langle q|\alpha\rangle+
 &\,\,\scri f(\alpha)-\varepsilon)\\
 &\ge
  \sum
  \Sb
  \bold i\\
  {}\\
  UL_{|\bold i|}[\bold i]\subseteq B(\alpha,r)
  \endSb
  \exp
 \Bigg(
 \big(
\,\,\scri f(\alpha)-(\varepsilon-\|q\|r)
\big)
 \sup_{\bold u\in[\bold i]}
 \sum_{k=0}^{|\bold i|-1}\varphi(S^{k}\bold u)
 \Bigg)\\
  &\ge
  \sum
  \Sb
  \bold i\\
  {}\\
  UL_{|\bold i|}[\bold i]\subseteq B(\alpha,\rho)
  \endSb
  \exp
 \Bigg(
 \big(
\,\,\scri f(\alpha)-(\varepsilon-\|q\|r)
\big)
 \sup_{\bold u\in[\bold i]}
 \sum_{k=0}^{|\bold i|-1}\varphi(S^{k}\bold u)
 \Bigg)\\
 &{}\\
 &\qquad\qquad
   \qquad\qquad
   \qquad\qquad
   \qquad\qquad
   \text{[since $B(\alpha,\rho)\subseteq B(\alpha,r)$ because $\rho<r$]} \\
 &{}\\
  &=
  \sum_{n}
  \sum
  \Sb
  \bold i\in\Sigma_{\smallsmallG}^{n}\\
  {}\\
  UL_{|\bold i|}[\bold i]\subseteq B(\alpha,\rho)
  \endSb
  \exp
 \Bigg(
 \big(
\,\,\scri f(\alpha)-(\varepsilon-\|q\|r)
\big)
 \sup_{\bold u\in[\bold i]}
 \sum_{k=0}^{n-1}\varphi(S^{k}\bold u)
 \Bigg)\\
 &\ge
  \sum_{n}
  \frac{1}{n!}
  \sum
  \Sb
  \bold i\in\Sigma_{\smallsmallG}^{n}\\
  {}\\
  UL_{|\bold i|}[\bold i]\subseteq B(\alpha,\rho)
  \endSb
  \exp
 \Bigg(
 \big(
\,\,\scri f(\alpha)-(\varepsilon-\|q\|r)
\big)
 \sup_{\bold u\in[\bold i]}
 \sum_{k=0}^{n-1}\varphi(S^{k}\bold u)
 \Bigg)\\
 &=
 \infty\,.
 \qquad\qquad
 \qquad\qquad
 \qquad\qquad
 \quad\,
 \text{[by (8.4) because $\rho<\rho_{\varepsilon-\|q\|r}$]}
 \tag8.5
  \endalign
  $$
It follows from (8.5) that
$-\langle q|\alpha\rangle+\,\,\scri f(\alpha)-\varepsilon
\le
 \sigma_{\abs}
 \big(
 \,
 \zeta_{q}^{\co,U}(\varphi;\cdot)
 \,
 \big)
=
\tau(q)$, and so
$\,\,\scri f(\alpha)\le \langle q|\alpha\rangle+\tau(q)+\varepsilon$.
Letting $\varepsilon$ tend to $0$ and taking infimum over all $q$, 
now shows that
$\,\,\scri f(\alpha)\le\inf_{q} (\langle q|\alpha\rangle+\tau(q))
=
\tau^{*}(\alpha)$.
This completes the proof.
\hfill$\square$

\bigskip

\noindent
We can now prove Theorem 4.3.

\bigskip

\noindent{\it Proof of Theorem 4.3}\newline
\noindent
Using Theorem C and Proposition 8.1 we conclude that
 $$
 \align
 \scri f(C)
&=
\sup
\Sb
\mu\in\Cal P_{S}(\Sigma^{\Bbb N})\\
{}\\
U\mu\in \overline C
\endSb
-
\frac{h(\mu)}{\int\varphi\,d\mu}
\qquad\qquad
\qquad\qquad
\text{[by Theorem C]}\\
&=
\sup
_{\alpha\in\overline C}
\,\,
\sup
\Sb
\mu\in\Cal P_{S}(\Sigma^{\Bbb N})\\
{}\\
U\mu=\alpha
\endSb
-
\frac{h(\mu)}{\int\varphi\,d\mu}\\
&=
\sup
_{\alpha\in\overline C}
\,\,
 \,\,\,\scri f(\alpha)
 \qquad\qquad
\qquad\qquad
\qquad\quad\,\,\,
\text{[by Theorem C]}\\
 &\le
\sup
_{\alpha\in\overline C}
\,\,
\tau^{*}(\alpha)\,.
\qquad\qquad
\qquad\qquad
\qquad\,\,\,\,
\text{[by Proposition 8.1]}
\endalign
$$
 This completes the proof.
 \hfill$\square$

  \bigskip


\heading
{
9. Proof of Theorem 4.4}
\endheading

The purpose of this section 
is to prove Theorem 4.4.
We start by
proving a simple , but useful,
auxiliary result, namely, Lemma 9.1 below.
Lemma 9.1 says 
that if $C\subset X$, then
$\,\,\scri f$
and 
$\tau$
satisfy the multifractal formalism
at $C$, 
i\.e\.
$  \scri f(C)
  =
  \sup_{\alpha\in \overline C}\tau^{*}(\alpha)$,
  provided 
  each
$\alpha\in\overline C\cap U\big(\,\,\Cal P_{S}(\Sigma_{\smallG}^{\Bbb N})\big)$
with $\tau^{*}(\alpha)>-\infty$
is the 
image 
of a \lq\lq Gibbs like" measure,
i\.e\.
there is a
measure $\mu_{\alpha}$
with a
certain 
 \lq\lq Gibbs like" property
 (namely, property (ii) in Lemma 9.1)
such that
$\alpha=U\mu_{\alpha}$.
This lemma plays a key role in the proof of Theorem 4.4.

\bigskip

\proclaim{Lemma 9.1}
Let $X$ be an inner product space 
with inner product $\langle\cdot |\cdot\rangle$
and let $U:\Cal P(\Sigma_{\smallG}^{\Bbb N})\to X$ be 
continuous with respect to the weak topology.
Fix a H\"older continuous 
function  $\varphi:\Sigma_{\smallG}^{\Bbb N}\to\Bbb R$ 
with
$\varphi<0$.
Let $C\subseteq X$ be a  subset of $X$
and
let
$\,\,\scri f(C)$
be the unique real number
such that
 $$
 \align
 \limsup_{r\searrow 0}
 \sigma_{\radius}
 \big(
 \,
 \zeta_{B(C,r)}^{\dyn,U}(\,\,\scri f(C)\,\varphi;\cdot)
 \,
 \big)
 &=
 1\,.
 \endalign
 $$ 
 
 \noindent
 Assume that for each
$\alpha\in\overline C\cap U\big(\,\,\Cal P_{S}(\Sigma_{\smallG}^{\Bbb N})\big)$
with $\tau^{*}(\alpha)>-\infty$
there 
is a point
$q_{\alpha}\in X$
and a measure
$\mu_{\alpha}\in \Cal P_{S}(\Sigma_{\smallG}^{\Bbb N})$
such that
 \roster
 \medskip
 
 \item"(i)"
 $\alpha=U\mu_{\alpha}$;
 
 \medskip
 
 \item"(ii)"
 $\dsize{
 \sup_{\mu\in\Cal P_{S}(\Sigma_{\smallsmallG}^{\Bbb N})}
 \Bigg(
 -\frac{h(\mu)}{\int\varphi\,d\mu}-\langle q_{\alpha}|U\mu\rangle
 \Bigg)
 =
 -\frac{h(\mu_{\alpha})}{\int\varphi\,d\mu_{\alpha}}-\langle q_{\alpha}|U\mu_{\alpha}\rangle
 }$.
 
 \medskip
 \endroster

\noindent
Then
 $$
  \scri f(C)
  =
  \sup_{\alpha\in \overline C}\tau^{*}(\alpha)\,.
  $$

\endproclaim

\noindent{\it Proof}\newline
\noindent
For 
$\alpha\in\overline C\cap U\big(\,\,\Cal P_{S}(\Sigma_{\smallG}^{\Bbb N})\big)$
with $\tau^{*}(\alpha)>-\infty$
we have
(using Theorem 4.1 and Theorem C)
 $$
 \align
 \tau^{*}(\alpha)
&=
\inf_{q}(\langle q|\alpha\rangle+\tau(q))\\
&\le
\langle q_{\alpha}|\alpha\rangle+\tau(q_{\alpha})\\
&=
\langle q_{\alpha}|\alpha\rangle
+
\sup_{\mu\in\Cal P_{S}(\Sigma_{\smallsmallG}^{\Bbb N})}
 \Bigg(
 -\frac{h(\mu)}{\int\varphi\,d\mu}-\langle q_{\alpha}|U\mu\rangle
 \Bigg)
 \qquad\qquad
 \text{[by Theorem 4.1]}\\
&=
\langle q_{\alpha}|U\mu_{\alpha}\rangle
 -
 \frac{h(\mu_{\alpha})}{\int\varphi\,d\mu_{\alpha}}
 -
 \langle q_{\alpha}|U\mu_{\alpha}\rangle
 \qquad\qquad
 \qquad\qquad
 \,
 \text{[by (i) and (ii)]}\\ 
 &\le
 -
 \frac{h(\mu_{\alpha})}{\int\varphi\,d\mu_{\alpha}}\\
 &\le
 \sup
 \Sb
 \mu\in\Cal P_{S}(\Sigma_{\smallsmallG}^{\Bbb N})\\
 {}\\
 U\mu\in\overline C
 \endSb
 -\frac{h(\mu)}{\int\varphi\,d\mu}\\
 &=
  \,\, \scri f(C)\,.
  \qquad\qquad
 \qquad\qquad
 \qquad\qquad
 \qquad\qquad
 \qquad\qquad
 \text{[by Theorem C]}
 \tag9.1
 \endalign
 $$
Of course, if
$\alpha\in\overline C\cap U\big(\,\,\Cal P_{S}(\Sigma_{\smallG}^{\Bbb N})\big)$
with
$  \tau^{*}(\alpha)=-\infty$, then
it is clear that
 $\tau^{*}(\alpha)
\le
  \,\, \scri f(C)$.
 It follows immediately from this and (9.1)
 that
  $$\sup_{
  \alpha\in
  \overline C\cap U(\Cal P_{S}(\Sigma_{\smallsmallG}^{\Bbb N}))
  }
  \tau^{*}(\alpha)
 \le
  \,\, \scri f(C)\,.
  $$
  Since it also follows from Theorem 4.3 that
  $\,\, \scri f(C)
  \le
  \sup_{\alpha\in\overline C}\tau^{*}(\alpha)
  =
  \sup_{
  \alpha\in
  \overline C\cap U(\Cal P_{S}(\Sigma_{\smallsmallG}^{\Bbb N}))
  }
  \tau^{*}(\alpha)$, we therefore conclude that
  $\,\, \scri f(C)
 =
  \sup_{\alpha\in\overline C}\tau^{*}(\alpha)$.
  \hfill$\square$

\bigskip

We can now prove Theorem 4.4.

\bigskip

\noindent{\it Proof of Theorem 4.4}\newline
\noindent
(1)
By Lemma  9.1 it suffices to show that for each 
$\alpha\in\overline C\cap U\big(\,\,\Cal P_{S}(\Sigma_{\smallG}^{\Bbb N})\big)$
with $\tau^{*}(\alpha)>-\infty$
there is a point
$q_{\alpha}\in \Bbb R^{M}$
and
a measure
$\mu_{\alpha}\in \Cal P_{S}(\Sigma_{\smallG}^{\Bbb N})$
such that
 \roster
 \medskip
 
 \item"(i)"
 $\alpha=U\mu_{\alpha}$;
 
 \medskip
 
 \item"(ii)"
 $\dsize{
 \sup_{\mu\in\Cal P_{S}(\Sigma_{\smallsmallG}^{\Bbb N})}
 \Bigg(
 -\frac{h(\mu)}{\int\varphi\,d\mu}-\langle q_{\alpha}|U\mu\rangle
 \Bigg)
 =
 -\frac{h(\mu_{\alpha})}{\int\varphi\,d\mu_{\alpha}}-\langle q_{\alpha}|U\mu_{\alpha}\rangle
 }$.
 
 \medskip
 \endroster

Define
$\Phi:\Bbb R^{M+1}\to\Bbb R$ by
  $$
  \Phi(q,s)
  =
  \sup
  _{\mu\in\Cal P_{S}(\Sigma_{\smallsmallG}^{\Bbb N})}
  \Bigg(
  h(\mu)
  +
  (
  \langle q|U\mu\rangle
  +
  s
  )
  \int\varphi\,d\mu
  \Bigg)
  $$
  for $q\in\Bbb R^{M}$ and $s\in\Bbb R$.
Also,
 for $q\in\Bbb R^{M}$ and $s\in\Bbb R$, we
write
  $$
  \Cal E_{q,s}
  =
  \Bigg\{
  \mu\in\Cal P_{S}(\Sigma_{\smallG}^{\Bbb N})
  \,\Bigg|\,
  \Phi(q,s)
  =
   h(\mu)
  +
  (
  \langle q|U\mu\rangle
  +
  s
  )
  \int\varphi\,d\mu
  \Bigg\}\,.
  $$
 We first show that  $\Cal E_{q,s}\not=\varnothing$.
This is the statement of Claim 1 below.
  
 \bigskip

{\it Claim 1. 
 For all $q\in\Bbb R^{M}$ and $s\in\Bbb R$,
 we have
$\Cal E_{q,s}\not=\varnothing$.
}

{\it Proof of Claim 1.} 
Fix $q\in\Bbb R^{M}$ and $s\in\Bbb R$
and
define the map
$H: \Cal P_{S}(\Sigma_{\smallG}^{\Bbb N})\to\Bbb R$
by
$H(\mu)
=
 h(\mu)
  +
  (
  \langle q|U\mu\rangle
  +
  s
  )
  \int\varphi\,d\mu$
  for $\mu\in \Cal P_{S}(\Sigma_{\smallG}^{\Bbb N})$.
 Since the entropy map $\mu\to h(\mu)$
  is upper semi-continuous (see,for example, [Wa, Theorem 8.2])
  and $U$ is continuous, we conclude that
  $H$ is upper semi-continuous.
  We therefore deduce from the compactness of 
$\Cal P_{S}(\Sigma_{\smallG}^{\Bbb N})$
with respect to the weak topology that the map 
$H$ has a maximum, i\.e\.
$\Cal E_{q,s}\not=\varnothing$.
This completes the proof of Claim 1.

\bigskip

\noindent
Since $\Cal E_{q,s}\not=\varnothing$, we can choose 
$\mu_{q,s}\in \Cal E_{q,s}$.

\bigskip

{\it Claim 2. For all $q$, we have
 $$
 \sup_{\mu\in\Cal P_{S}(\Sigma_{\smallsmallG}^{\Bbb N})}
 \Bigg(
 -\frac{h(\mu)}{\int\varphi\,d\mu}-\langle q|U\mu\rangle
 \Bigg)
 =
 -\frac{h(\mu_{q,\tau(q)})}{\int\varphi\,d\mu_{q,\tau(q)}}-\langle q|U\mu_{q,\tau(q)}\rangle\,.
$$
}

\noindent
{\it Proof of Claim 2.}
Since 
$\Phi(q,\tau(q))=0$
(by Theorem 4.1),
we deduce  that
if $\mu\in\Cal P_{S}(\Sigma_{\smallG}^{\Bbb N})$, then
 $
 h(\mu)
 +
(\langle q|U\mu\rangle+\tau(q))\int\varphi\,d\mu
\le
 \sup_{\nu\in\Cal P(\Sigma_{\smallsmallG}^{\Bbb N})}
 (
 h(\nu)
 +
 (\langle q|U\nu\rangle+\tau(q))\int\varphi\,d\nu
 )
= 
\Phi(q,\tau(q))
=
0$.
Dividing this inequality by $-\int\varphi\,d\mu$,
 clearly implies that
  $$
  \align
   -\frac{h(\mu)}{\int\varphi\,d\mu}-\langle q|U\mu\rangle-\tau(q)
   &
    \le
    0\,.
    \tag9.2
   \endalign
   $$
 In addition,
 since
 $\mu_{q,\tau(q)}\in\Cal E_{q,\tau(q)}$
and
$\Phi(q,\tau(q))=0$,
we conclude that
$
0
=
\Phi(q,\tau(q))
=
h(\mu_{q,\tau(q)})
 +
(\langle q|U\mu_{q,\tau(q)}\rangle+\tau(q))\int\varphi\,d\mu_{q,\tau(q)}$.
Dividing this equality
by $-\int\varphi\,d\mu_{q,\tau(q)}$,
 shows that
 $$
  \align
  0
  &    =
   -\frac{h(\mu_{q,\tau(q)})}{\int\varphi\,d\mu_{q,\tau(q)}}-\langle q|U\mu_{q,\tau(q)}\rangle-\tau(q)\,.
   \tag9.3
       \endalign
   $$
Finally, combining (9.2) and (9.3), we
deduce that 
if $\mu\in\Cal P_{S}(\Sigma_{\smallG}^{\Bbb N})$, then
  $$
   -\frac{h(\mu)}{\int\varphi\,d\mu}-\langle q|U\mu\rangle
  \le
-\frac{h(\mu_{q,\tau(q)})}{\int\varphi\,d\mu_{q,\tau(q)}}-\langle q|U\mu_{q,\tau(q)}\rangle\,.
   $$
The desired conclusion follows immediately from this
inequality.
This completes the proof of Claim 2.

\bigskip

\noindent
Next, we 
compute the 
gradient $\nabla\tau(q)$
of $\tau$
at points $q$ at which $\tau$ is differentiable.
This is the statement of Claim 3.
In Claim 3 we use the 
following notation, namely,
for $i=1,\ldots,M$,
we write
$D_{i}^{-}\tau(q)$
and 
$D_{i}^{+}\tau(q)$
for the left
and right  partial
derivative
of $\tau$ with respect to the $i$'th variable at the point $q$, respectively;
note
that since $\tau$ is convex (by Theorem 4.2), 
we conclude that
the left
and right  partial
derivatives
of $\tau$ with respect to the $i$'th variable exist at all points.
Also, if $\tau$
is partially differentiable with respect to the $i$'th variable at $q$,
then   we write
$D_{i}\tau(q)$
for the  partial
derivative
of $\tau$ with respect to the $i$'th variable at the point $q$.

\bigskip

{\it Claim 3. For $i=1,\ldots,M$ and all $q\in\Bbb R^{M}$, we have
 $$
 D_{i}^{-}\tau(q)
 \le
-
  \langle e_{i}|U\mu_{q,\tau(q)}\rangle
\le
D_{i}^{+}\tau(q)
\tag9.4
 $$
where $e_{i}=(0,\ldots,0,1,0\ldots,0)$ 
is the canonical basis vector in $\Bbb R^{M}$ with
$1$ in the $i$'th coordinate and $0$'s elsewhere.
In particular, if $\tau$ is  partially differentiable at $q$,
then (9.4) implies that
$
 D_{i}\tau(q)
 =
 -
  \langle e_{i}|U\mu_{q,\tau(q)}\rangle
 $,
 and so
 $$
 \nabla\tau(q)
 =
 -U\mu_{q,\tau(q)}\,.
 $$
 }

{\it Proof of Claim 3.}
It follows Theorem 4.1 and Claim 2 that
if $h\in \Bbb R$, then we
have
 $$
 \align
 \tau(q+he_{i})
 -
 \tau(q)
&=
 \sup_{\mu\in\Cal P_{S}(\Sigma_{\smallsmallG}^{\Bbb N})}
 \Bigg(
 -\frac{h(\mu)}{\int\varphi\,d\mu}-\langle q+he_{i}|U\mu\rangle
 \Bigg)
-
 \sup_{\mu\in\Cal P_{S}(\Sigma_{\smallsmallG}^{\Bbb N})}
 \Bigg(
 -\frac{h(\mu)}{\int\varphi\,d\mu}-\langle q|U\mu\rangle
 \Bigg)\\
&\ge
\Bigg(
  -\frac{h(\mu_{q,\tau(q)})}{\int\varphi\,d\mu_{q,\tau(q)}}-\langle q+he_{i}|U\mu_{q,\tau(q)}\rangle
  \Bigg)
  -
  \Bigg(
 -\frac{h(\mu_{q,\tau(q)})}{\int\varphi\,d\mu_{q,\tau(q)}}-\langle q|U\mu_{q,\tau(q)}\rangle
  \Bigg)\\
 &=
   -h\langle e_{i}|U\mu_{q,\tau(q)}\rangle\,.
    \tag9.5
\endalign
 $$
For $h<0$, (9.5) implies that
$\frac{\tau(q+he_{i})
 -
 \tau(q)
}
{h}
\le
 -\langle e_{i}|U\mu_{q,\tau(q)}\rangle$
  and for $h>0$, (9.5) implies that
$ -\langle e_{i}|U\mu_{q,\tau(q)}\rangle
 \le
\frac{\tau(q+he_{i})
 -
 \tau(q)
}
{h}$.
The result follows immediately from these inequalities.
This 
completes the proof of Claim 3.

\bigskip

 Let $\alpha\in\overline C\cap U\big(\,\Cal P_{S}(\Sigma_{\smallG}^{\Bbb N})\,\big)$
 with 
  $ \tau^{*}(\alpha)>-\infty$.
  Since 
  $ \tau^{*}(\alpha)>-\infty$,
   it follows from the convexity of $\tau$
that there is a point $q_{\alpha}\in\Bbb R^{M}\cap Q$
such that
$\alpha=-\nabla\tau(q_{\alpha})$, see [Ro].
Now put
$\mu_{\alpha}=\mu_{q_{\alpha},\tau(q_{\alpha})}$.
We claim that $\mu_{\alpha}$ satisfies (i) and (ii).
Indeed,
it follows from Claim 2 and Claim 3 that
 $$
 \alpha
 =
 -\nabla\tau(q_{\alpha})
=
U\mu_{q_{\alpha},\tau(q_{\alpha})}
=
U\mu_{\alpha}
$$
and
$$
 \sup_{\mu\in\Cal P_{S}(\Sigma_{\smallsmallG}^{\Bbb N})}
 \Bigg(
 -\frac{h(\mu)}{\int\varphi\,d\mu}-\langle q_{\alpha}|U\mu\rangle
 \Bigg)
 =
 -\frac{h(\mu_{q_{\alpha},\tau(q_{\alpha})})}{\int\varphi\,d\mu_{q_{\alpha},\tau(q_{\alpha})}}
 -
 \langle q_{\alpha}|U\mu_{q_{\alpha},\tau(q_{\alpha})}\rangle
 =
 -\frac{h(\mu_{\alpha})}{\int\varphi\,d\mu_{\alpha}}
 -
 \langle q_{\alpha}|U\mu_{\alpha}\rangle\,.
$$
This completes the proof
of
(1).

\noindent
(2)
This follows immediately by applying
the statement in Part (1)
to the set $C=\{\alpha\}$.
\hfill$\square$

\end